\documentclass[reqno,11pt]{amsart}
\usepackage{mathrsfs}
\usepackage{verbatim}
\usepackage{upgreek}
\usepackage{amsmath}
\usepackage{amsxtra}
\usepackage{amscd}
\usepackage{amsthm}
\usepackage{amsfonts}
\usepackage{amssymb}
\usepackage{eucal}
\usepackage{tikz-cd}
\usepackage[matrix,arrow,curve]{xy}
\usepackage{young}
\usepackage[vcentermath]{youngtab}
\usepackage{mathtools}

\theoremstyle{plain}
\newtheorem{Thm}{Theorem}[section]
\newtheorem{Cor}[Thm]{Corollary}
\newtheorem{Lem}[Thm]{Lemma}
\newtheorem{Prop}[Thm]{Proposition}
\newtheorem{Conj}[Thm]{Conjecture}
\newtheorem{Warning}[Thm]{Warning}
\newtheorem{Assumption}[Thm]{Assumption}

\newtheorem{thmx}{Theorem}

\theoremstyle{definition}
\newtheorem{Def}[Thm]{Definition}

\theoremstyle{definition}
\newtheorem{Ex}[Thm]{Example}

\theoremstyle{remark}
\newtheorem{Rem}[Thm]{Remark}

\newtheorem{conj}{Conjecture}

\numberwithin{equation}{section}
\renewcommand{\rm}{\normalshape}

\newif\ifShowLabels
\ShowLabelstrue
\newdimen\theight
\def\TeXref#1{%
	\leavevmode\vadjust{\setbox0=\hbox{{\tt
				\quad\quad  {\small \rm #1}}}%
		\theight=\ht0
		\advance\theight by \lineskip
		\kern -\theight \vbox to
		\theight{\rightline{\rlap{\box0}}%
			\vss}%
}}%

\ShowLabelsfalse

\newcommand{\ssec}[2]{\subsection{#2}\label{SS:#1}%
	\ifShowLabels \TeXref{{SS:#1}} \fi}

\newcommand{\sssec}[2]{\subsubsection{#2}\label{SSS:#1}%
	\ifShowLabels \TeXref{{SSS:#1}} \fi}

\newenvironment{thm}[1]%
{ \begin{Thm} \label{T:#1}  \ifShowLabels \TeXref{T:#1} \fi }%
	{ \end{Thm} }

\renewcommand{\th}[1]{\begin{thm}{#1} \sl }
	\renewcommand{\eth}{\end{thm} }

\newenvironment{lemma}[1]%
{ \begin{Lem} \label{L:#1}  \ifShowLabels \TeXref{L:#1} \fi }%
	{ \end{Lem} }
\newcommand{\lem}[1]{\begin{lemma}{#1} \sl}
	\newcommand{\elem}{\end{lemma}}

\newenvironment{warning}[1]%
{ \begin{Warning} \label{W:#1}  \ifShowLabels \TeXref{W:#1} \fi }%
	{ \end{Warning} }
\newcommand{\war}[1]{\begin{warning}{#1} \sl}
	\newcommand{\ewar}{\end{warning}}

\newenvironment{assumption}[1]%
{ \begin{Assumption} \label{A:#1}  \ifShowLabels \TeXref{A:#1} \fi }%
	{ \end{Assumption} }
\newcommand{\ass}[1]{\begin{assumption}{#1} \sl}
	\newcommand{\eass}{\end{assumption}}	
	
\newenvironment{propos}[1]%
{ \begin{Prop} \label{P:#1}  \ifShowLabels \TeXref{P:#1} \fi }%
	{ \end{Prop} }
\newcommand{\prop}[1]{\begin{propos}{#1}\sl }
	\newcommand{\eprop}{\end{propos}}

\newenvironment{corol}[1]%
{ \begin{Cor} \label{C:#1}  \ifShowLabels \TeXref{C:#1} \fi }%
	{ \end{Cor} }
\newcommand{\cor}[1]{\begin{corol}{#1} \sl }
	\newcommand{\ecor}{\end{corol}}

\newenvironment{defeni}[1]%
{ \begin{Def} \label{D:#1}  \ifShowLabels \TeXref{D:#1} \fi }%
	{ \end{Def} }
\newcommand{\defe}[1]{\begin{defeni}{#1} \sl }
	\newcommand{\edefe}{\end{defeni}}

\newenvironment{remark}[1]%
{ \begin{Rem} \label{R:#1}  \ifShowLabels \TeXref{R:#1} \fi }%
	{ \end{Rem} }
\newcommand{\rem}[1]{\begin{remark}{#1}}
	\newcommand{\erem}{\end{remark}}

\newenvironment{conjec}[1]%
{ \begin{Conj} \label{Co:#1}  \ifShowLabels \TeXref{Co:#1} \fi }%
	{ \end{Conj} }
\renewcommand{\conj}[1]{\begin{conjec}{#1} \sl }
	\newcommand{\econj}{\end{conjec}}

\newenvironment{First proof}[1]%
{ \begin{First proof} \label{Co:#1}  \ifShowLabels \TeXref{Co:#1} \fi }%
	{ \end{First proof} }

\newenvironment{Second proof}[1]%
{ \begin{Second proof} \label{Co:#1}  \ifShowLabels \TeXref{Co:#1} \fi }%
	{ \end{Second proof} }

\newcommand{\eq}[1]%
{ \ifShowLabels \TeXref{E:#1} \fi
	\begin{equation} \label{E:#1} }
\newcommand{\eeq}{ \end{equation} }

\newcommand{\prf}{ \begin{proof} }
	\newcommand{\epr}{ \end{proof} }

\newcommand{\prft}{ \begin{proof} }
	\newcommand{\eprt}{ \end{proof} }

\newcommand\nc{\newcommand}

\nc{\HC}{{\mathcal{HC}}}
\nc{\on}{\operatorname}
\nc{\BA}{{\mathbb{A}}}
\nc{\BC}{{\mathbb{C}}}
\nc{\BF}{{\mathbb{F}}}
\nc{\BG}{{\mathbb{G}}}
\nc{\BM}{{\mathbb{M}}}
\nc{\BN}{{\mathbb{N}}}
\nc{\BO}{{\mathbb{O}}}
\nc{\BQ}{{\mathbb{Q}}}
\nc{\BP}{{\mathbb{P}}}
\nc{\BR}{{\mathbb{R}}}
\nc{\BZ}{{\mathbb{Z}}}
\nc{\BS}{{\mathbb{S}}}

\nc{\CA}{{\mathcal{A}}}
\nc{\CB}{{\mathcal{B}}}
\nc{\CalC}{{\mathcal C}}
\nc{\CalD}{{\mathcal D}}
\nc{\CE}{{\mathcal{E}}}
\nc{\CF}{{\mathcal{F}}}
\nc{\CG}{{\mathcal{G}}}
\nc{\CH}{{\mathcal{H}}}
\nc{\CK}{{\mathcal{K}}}
\nc{\CL}{{\mathcal{L}}}
\nc{\CM}{{\mathcal{M}}}
\nc{\CMM}{{\mathcal{M}^{\operatorname{gen}}_\hbar(-\rho)}}
\nc{\CN}{{\mathcal{N}}}
\nc{\CO}{{\mathcal{O}}}
\nc{\CP}{{\mathcal{P}}}
\nc{\CQ}{{\mathcal{Q}}}
\nc{\CR}{{\mathcal{R}}}
\nc{\CS}{{\mathcal{S}}}
\nc{\CT}{{\mathcal{T}}}
\nc{\CU}{{\mathcal{U}}}
\nc{\CV}{{\mathcal{V}}}
\nc{\CW}{{\mathcal{W}}}
\nc{\CX}{{\mathcal{X}}}
\nc{\CY}{{\mathcal{Y}}}
\nc{\CZ}{{\mathcal{Z}}}

\nc{\gen}{{\operatorname{gen}}}
\nc{\cM}{{\check{\mathcal M}}{}}
\nc{\csM}{{\check{\mathcal A}}{}}
\nc{\obM}{{\overset{\circ}{\mathbf M}}{}}
\nc{\oCA}{{\overset{\circ}{\mathcal A}}{}}
\nc{\obA}{{\overset{\circ}{\mathbf A}}{}}
\nc{\ooM}{{\overset{\circ}{M}}{}}
\nc{\osM}{{\overset{\circ}{\mathsf M}}{}}
\nc{\vM}{{\overset{\bullet}{\mathcal M}}{}}
\nc{\nM}{{\underset{\bullet}{\mathcal M}}{}}
\nc{\obD}{{\overset{\circ}{\mathbf D}}{}}
\nc{\cp}{{\overset{\circ}{\mathbf p}}{}}
\nc{\ofZ}{{\overset{\circ}{\mathfrak Z}}{}}

\nc{\fa}{{\mathfrak{a}}}
\nc{\fb}{{\mathfrak{b}}}
\nc{\fg}{{\mathfrak{g}}}
\nc{\fgl}{{\mathfrak{gl}}}
\nc{\fh}{{\mathfrak{h}}}
\nc{\fj}{{\mathfrak{j}}}
\nc{\fm}{{\mathfrak{m}}}
\nc{\fn}{{\mathfrak{n}}}
\nc{\fu}{{\mathfrak{u}}}
\nc{\fp}{{\mathfrak{p}}}
\nc{\frr}{{\mathfrak{r}}}
\nc{\fs}{{\mathfrak{s}}}
\nc{\ft}{{\mathfrak{t}}}
\nc{\fT}{{\mathfrak{T}}}
\nc{\ofT}{{\overline{\mathfrak T}}}
\nc{\ofS}{{\overline{\mathfrak S}}}
\nc{\fsl}{{\mathfrak{sl}}}
\nc{\hsl}{{\widehat{\mathfrak{sl}}}}
\nc{\hgl}{{\widehat{\mathfrak{gl}}}}
\nc{\hg}{{\widehat{\mathfrak{g}}}}
\nc{\chg}{{\widehat{\mathfrak{g}}}{}^\vee}
\nc{\hn}{{\widehat{\mathfrak{n}}}}
\nc{\chn}{{\widehat{\mathfrak{n}}}{}^\vee}

\nc{\fA}{{\mathfrak{A}}}
\nc{\fB}{{\mathfrak{B}}}
\nc{\fD}{{\mathfrak{D}}}
\nc{\fE}{{\mathfrak{E}}}
\nc{\fF}{{\mathfrak{F}}}
\nc{\fG}{{\mathfrak{G}}}
\nc{\fI}{{\mathfrak{I}}}
\nc{\fJ}{{\mathfrak{J}}}
\nc{\fK}{{\mathfrak{K}}}
\nc{\fL}{{\mathfrak{L}}}
\nc{\fM}{{\mathfrak{M}}}
\nc{\fN}{{\mathfrak{N}}}
\nc{\frP}{{\mathfrak{P}}}
\nc{\fQ}{{\mathfrak Q}}
\nc{\fS}{{\mathfrak S}}
\nc{\fU}{{\mathfrak{U}}}
\nc{\fZ}{{\mathfrak{Z}}}

\nc{\ba}{{\mathbf{a}}}
\nc{\bb}{{\mathbf{b}}}
\nc{\bc}{{\mathbf{c}}}
\nc{\bd}{{\mathbf{d}}}
\nc{\be}{{\mathbf{e}}}
\nc{\bi}{{\mathbf{i}}}
\nc{\bj}{{\mathbf{j}}}
\nc{\bn}{{\mathbf{n}}}
\nc{\bp}{{\mathbf{p}}}
\nc{\br}{{\mathbf{r}}}
\nc{\bq}{{\mathbf{q}}}
\nc{\bu}{{\mathbf{u}}}
\nc{\bv}{{\mathbf{v}}}
\nc{\bx}{{\mathbf{x}}}
\nc{\by}{{\mathbf{y}}}
\nc{\bw}{{\mathbf{w}}}
\nc{\bA}{{\mathbf{A}}}
\nc{\bB}{{\mathbf{B}}}
\nc{\bC}{{\mathbf{C}}}
\nc{\bD}{{\mathbf{D}}}
\nc{\bE}{{\mathbf{E}}}
\nc{\bK}{{\mathbf{K}}}
\nc{\bH}{{\mathbf{H}}}
\nc{\bM}{{\mathbf{M}}}
\nc{\bN}{{\mathbf{N}}}
\nc{\bO}{{\mathbf{O}}}
\nc{\bQ}{{\mathbf Q}}
\nc{\bS}{{\mathbf{S}}}
\nc{\bT}{{\mathbf{T}}}
\nc{\bV}{{\mathbf{V}}}
\nc{\bW}{{\mathbf{W}}}
\nc{\bX}{{\mathbf{X}}}
\nc{\bP}{{\mathbf{P}}}
\nc{\bZ}{{\mathbf{Z}}}

\nc{\sA}{{\mathsf{A}}}
\nc{\sB}{{\mathsf{B}}}
\nc{\sC}{{\mathsf{C}}}
\nc{\sD}{{\mathsf{D}}}
\nc{\sF}{{\mathsf{F}}}
\nc{\sK}{{\mathsf{K}}}
\nc{\sM}{{\mathsf{M}}}
\nc{\sO}{{\mathsf{O}}}
\nc{\sQ}{{\mathsf{Q}}}
\nc{\sP}{{\mathsf{P}}}
\nc{\sV}{{\mathsf{V}}}
\nc{\sW}{{\mathsf{W}}}
\nc{\sZ}{{\mathsf{Z}}}
\nc{\sfp}{{\mathsf{p}}}
\nc{\sr}{{\mathsf{r}}}
\nc{\st}{{\mathsf{t}}}
\nc{\sfb}{{\mathsf{b}}}
\nc{\sfc}{{\mathsf{c}}}
\nc{\sd}{{\mathsf{d}}}
\nc{\sg}{{\mathsf{g}}}
\nc{\sk}{{\mathsf{k}}}
\nc{\sfl}{{\mathsf{l}}}

\nc{\BK}{{\bar{K}}}

\nc{\tA}{{\widetilde{\mathbf{A}}}}
\nc{\tB}{{\widetilde{\mathcal{B}}}}
\nc{\tg}{{\widetilde{\mathfrak{g}}}}
\nc{\tG}{{\widetilde{G}}}
\nc{\TM}{{\widetilde{\mathbb{M}}}{}}
\nc{\tN}{{\widetilde{\mathcal{N}}}{}}
\nc{\tO}{{\widetilde{\mathsf{O}}}{}}
\nc{\tU}{{\widetilde{\mathfrak{U}}}{}}
\nc{\TZ}{{\tilde{Z}}}
\nc{\tZ}{\widetilde{Z}{}}
\nc{\tx}{{\tilde{x}}}
\nc{\tbv}{{\tilde{\bv}}}
\nc{\tfP}{{\widetilde{\mathfrak{P}}}{}}
\nc{\tz}{{\tilde{\zeta}}}
\nc{\tmu}{{\tilde{\mu}}}

\nc{\td}{\ddot{\underline{d}}{}}
\nc{\tzeta}{\widetilde{\zeta}{}}
\nc{\hd}{{\widehat{\underline{d}}}}
\nc{\hG}{{\widehat{G}}}
\nc{\hBP}{\widehat{\mathbb P}{}}
\nc{\hQ}{{\widehat{Q}}}
\nc{\hsM}{\widehat{\mathsf M}{}}
\nc{\hfM}{\widehat{\mathfrak M}{}}
\nc{\hCP}{\widehat{\mathcal P}{}}
\nc{\hCR}{\widehat{\mathcal R}{}}
\nc{\hCS}{{\widehat{\mathcal S}}}
\nc{\hfZ}{\widehat{\mathfrak Z}{}}
\nc{\hZ}{\widehat{Z}{}}

\nc{\urho}{\underline{\rho}}
\nc{\uB}{\underline{B}}
\nc{\uC}{{\underline{\mathbb{C}}}}
\nc{\ui}{\underline{i}}
\nc{\ofP}{{\overline{\mathfrak{P}}}}

\nc{\hrho}{{\hat{\rho}}}

\nc{\unl}{\underline}
\nc{\ol}{\overline}
\nc{\one}{{\mathbf{1}}}
\nc{\two}{{\mathbf{t}}}

\nc{\Sym}{{\mathop{\operatorname{Sym}}}}
\nc{\Tot}{{\mathop{\operatorname{\normalshape Tot}}}}
\nc{\Hilb}{{\mathop{\operatorname{\normalshape Hilb}}}}
\nc{\Hom}{{\mathop{\operatorname{Hom}}}}
\nc{\CHom}{{\mathop{\operatorname{{\mathcal{H}}\it om}}}}
\nc{\defi}{{\mathop{\operatorname{\normalshape def}}}}
\nc{\length}{{\mathop{\operatorname{\normalshape length}}}}

\nc{\Cliff}{{\mathsf{Cliff}}}
\nc{\Fl}{{\mathcal{F}\ell}}
\nc{\Fib}{{\mathsf{Fib}}}
\nc{\Coh}{{\mathsf{Coh}}}
\nc{\FCoh}{{\mathsf{FCoh}}}

\nc{\reg}{{\text{\normalshape reg}}}
\nc{\res}{{\operatorname{res}}}
\nc{\cplus}{{\mathbf{C}_+}}
\nc{\cminus}{{\mathbf{C}_-}}
\nc{\cthree}{{\mathbf{C}_*}}
\nc{\Qbar}{{\bar{Q}}}

\nc{\bh}{{\bar{h}}}
\nc{\bOmega}{{\overline{\Omega}}}
\nc\tGr{\widetilde{\Gr}}

\nc{\seq}[1]{\stackrel{#1}{\sim}}
\nc\ogu{\overline{G/U}}
\nc\chlam{\check{\lam}}

\nc\St{\operatorname{St}}
\nc{\oZ}{{\overset{\circ}{Z}}}
\nc{\tF}{\widetilde{\mathcal F}}

\nc\uS{\underline{S}}
\nc\QM{\mathcal{QM}}

\nc{\chmu}{\check{\mu}}
\newcommand\iso{\,\vphantom{j^{X^2}}\smash{\overset{\sim}{\vphantom{\rule{0pt}{0.20em}}\smash{\longrightarrow}}}\,}
\nc{\ul}{\underline}
\nc{\Mvd}{\mathfrak{M}(\underline{v},\underline{d})}
\nc{\MVD}{\mathfrak{M}(V,D)}
\nc{\mt}{\mapsto}
\nc{\sm}{\setminus}
\nc{\ra}{\rightarrow}
\nc{\lar}{\leftarrow}
\nc{\hr}{\hookrightarrow}
\nc{\La}{\Lambda}
\nc{\Lap}{\Lambda^{+}}
\nc{\oZal}{\overset{\circ}{Z^{\alpha}}}
\nc{\sig}{\sigma}
\nc{\al}{\alpha}
\nc{\la}{\lambda}
\nc{\is}{\simeq}
\nc{\ip}{\iota^{+}_{\la, \mu}}
\nc{\im}{\iota^{-}_{\la, \mu}}
\nc{\jp}{j^{+}_{\la, \mu}}
\nc{\jm}{j^{-}_{\la, \mu}}
\nc{\pip}{\pi^{+}_{\la, \mu}}
\nc{\pim}{\pi^{-}_{\la, \mu}}
\nc{\s}{\star}
\nc{\fpt}{[A^{\la},B^{\la},\gamma^{\la},\delta^{\la}]}
\nc{\ulfpt}{[\ul{A}^{\la},\ul{B}^{\la},\ul{\gamma}^{\la},\ul{\delta}^{\la}]}
\nc{\lvee}{\!\scriptscriptstyle\vee}
\nc{\Gr}{{\operatorname{Gr}}}
\nc{\rra}{\twoheadrightarrow}

\nc{\End}{\on{End}}

\nc{\RHom}{R\mathcal{H}om}

\nc{\yg}{Y(\fg)}
\nc{\yvg}{Y_V(\fg)}
\nc{\CAg}{\CA_{\fg}}
\nc{\Ag}{A_{\fg}}
\nc{\Sgt}{S^{\bullet}(\fg[t])}
\nc{\Sg}{S^{\bullet}(\fg)}
\nc{\Ugt}{U(\fg[t])}

\nc{\Spec}{\operatorname{Spec}}

\author{Vasily Krylov}
\address{V. K.: 
Department of Mathematics
Harvard University and CMSA
\newline
1 Oxford Street,
Cambridge, MA 02138,
USA}
\email{vkrylov@math.harvard.edu, krylovasya@gmail.com}

\author{Inna Mashanova-Golikova}
\address{I. M.-G.} 
\email{inna.mashanova@gmail.com}

\author{Leonid Rybnikov}
\address{L. R.: Department of Mathematics and Statistics, University of Montreal, \newline Montreal QC, Canada}
\email{leonid.rybnikov@umontreal.ca}

\begin{document}
	\begin{abstract}
Let $\mathfrak{g}$ be a simple finite-dimensional Lie algebra and $G$ its adjoint group. For each $C\in G$, we consider the Bethe subalgebra $B(C)\subset Y(\mathfrak{g})$, a commutative subalgebra encoding the integrals of the generalized $XXX$ spin chain. Adapting the construction of~\cite{hkrw} of $\mathfrak{g}$-crystals on
spectra of inhomogeneous Gaudin subalgebras in $U(\mathfrak{g})$, we construct a natural $\hat{\mathfrak{g}}$-crystal structure on the spectra of $B(C)$ in Kirillov--Reshetikhin $Y(\mathfrak{g})$-modules in type $A$. We conjecture that such a construction exists for arbitrary $\mathfrak{g}$ and recovers Kirillov--Reshetikhin crystals. 
The main technical ingredient is a degeneration of Bethe subalgebras $B(C)$ to  commutative subalgebras $\mathcal{A}_\chi^{\mathrm{u}} \subset U(\mathfrak{g}[t])$, depending on
$\chi\in\mathfrak{g}$. We call these subalgebras universal inhomogeneous Gaudin subalgebras and show that they arise from the Feigin--Frenkel center at the critical level. This allows us to identify the affine crystals 
above with Kirillov--Reshetikhin crystals.
We then apply these results to prove the
monodromy conjecture of Ilin and the second and third authors
\cite[Conjecture~7.1]{imr}, see also~\cite[Conjecture~1.6]{ikr}, for the
spectra of the algebras $B(C)$ and for the spectra of quantum cohomology
rings of type $A$ quiver varieties.
	\end{abstract}

	\title[]{Bethe subalgebras in Yangians and Kirillov--Reshetikhin crystals}
	
	\maketitle

\tableofcontents






\section{Introduction}

\ssec{}{Crystals}\label{Crystals_and_Gaudin_HKRW}

Let $\mathfrak{g}$ be a simple finite-dimensional Lie algebra over $\BC$. 
Kashiwara $\mathfrak{g}$-crystals are combinatorial objects which model 
crystal bases of representations of $\mathfrak{g}$. Namely, a 
$\mathfrak{g}$-crystal $B$ is a set equipped with operators
\[
{\on{e}}_i,{\on{f}}_i\colon B \to B \cup \{0\},
\]
attached to the simple roots $\al_i$ of $\mathfrak{g}$, together with a 
weight function $\on{wt}$, satisfying certain conditions; see 
Section~\ref{cryst_main} for details. In particular, to each irreducible 
finite-dimensional representation $\mathfrak{g} \curvearrowright V_\la$ associated with the dominant weight $\lambda$ one can attach 
a connected  $\mathfrak{g}$-crystal $B_\la$. Given two $\mathfrak{g}$-crystals $B$ and $B'$, 
one can form their tensor product $B \otimes B'$.

One can replace $\mathfrak{g}$ by the corresponding {\emph{affine}} Kac--Moody Lie algebra $\hat{\mathfrak{g}}$. An important class of finite $\hat{\mathfrak{g}}$-crystals are (tensor products of) so-called Kirillov--Reshetikhin crystals modeling bases of Kirillov--Reshetikhin modules over the Yangian $Y(\mathfrak{g})$. Kirillov--Reshetikhin modules are labeled by dominant $\mathfrak{g}$-weights $\lambda$ such that $\lambda$ is a multiple of a fundamental weight. We denote the corresponding $\hat{\mathfrak{g}}$-crystal by ${\bf{B}}_\lambda$.

\subsection{Gaudin subalgebras} 

Let $\mathfrak{h} \subset \mathfrak{g}$ be a Cartan subalgebra, and let 
$\mathfrak{h}_{\mathbb{R}}$ denote its split real form. To every 
$\chi \in \mathfrak{h}$ and every collection of distinct points 
$z_1,\ldots,z_k \in \BC$, one can associate a commutative subalgebra
\[
\CA_\chi(\ul{z})=\CA_\chi(z_1,\ldots,z_k)
\subset U(\mathfrak{g})^{\otimes k},
\]
called the inhomogeneous Gaudin subalgebra of $U(\mathfrak{g})^{\otimes k}$; 
see, for example, \cite[Section~9]{hkrw} and the references therein. 

For $\chi \in \mathfrak{h}^{\mathrm{reg}}$, the algebra 
$\CA_\chi(z_1,\ldots,z_k)$ is a maximal commutative subalgebra of 
$U(\mathfrak{g})^{\otimes k}$. It contains $\mathfrak{h}$. These 
commutative subalgebras describe the quantum integrable spin chain known as 
the Gaudin magnet; see \cite{er,ffr,fft,ffry,r0}.

\subsection{$\mathfrak{g}$-crystals from Gaudin subalgebras} 
The crystals $B_\la$ and, more generally, tensor products
$B_{\la_1} \otimes \cdots \otimes B_{\la_k}$ can be constructed using 
inhomogeneous Gaudin subalgebras $\mathcal{A}_\chi(\underline{z})$; 
see~\cite{hkrw}. We briefly recall the construction.

Let $\la_1,\ldots,\la_k$ be dominant weights of $\mathfrak{g}$. We have a natural action
\[
\CA_\chi(z_1,\ldots,z_k) \curvearrowright
V_{\la_1} \otimes \cdots \otimes V_{\la_k}.
\]

The following proposition holds.

\prop{}{\em{(\cite{ffry})}}
For any dominant weights $\la_1,\ldots,\la_k$, any
$\chi \in \mathfrak{h}^{\mathrm{reg}}_{\BR}$, and any distinct
$z_1,\ldots,z_k \in \BR$, the action of $\CA_\chi(\ul{z})$ on
$
V_{\la_1} \otimes \cdots \otimes V_{\la_k}
$
has simple spectrum.
\eprop

Let $\CE_\chi(\ul{\la})$ denote the set of eigenlines of $\CA_\chi(\ul{z})$
acting on
\[
\ul{V}:=V_{\la_1} \otimes \cdots \otimes V_{\la_k}.
\]
We suppress the dependence on $\ul{z}$ from the notation. 
We next recall, following~\cite{hkrw}, how $\CE_\chi(\ul{\lambda})$ is endowed with a natural
$\mathfrak{g}$-crystal structure.

Let $H_\alpha\subset \mathfrak{h}_{\BR}$ be the wall corresponding to a
positive root $\alpha$. These walls divide $\mathfrak{h}_{\BR}$ into Weyl
chambers. We denote by $O$ the dominant chamber and assume 
that $\chi\in O^{\mathrm{reg}}$. Since
$\mathfrak{h}\subset \CA_\chi(\ul{z})$, every eigenline of
$\CA_\chi(\ul{z})$ is an $\mathfrak{h}$-eigenline. Thus, for
$L\in \CE_\chi(\ul{\lambda})$, we can define $\on{wt}(L)=\mu$ if $L$ has
$\mathfrak{h}$-weight $\mu$.

The crystal operators are defined by degenerating the Gaudin subalgebras to
the walls of the dominant chamber. More precisely, for a simple root
$\alpha_i$, choose a generic point $\chi_0$ in the wall
$H_{\alpha_i}\cap O$. Consider the limiting subalgebra
\begin{equation*}
\CA_{(\chi_0,\chi)}(\ul{z})
:=
\underset{\epsilon\to 0}{\lim}\,\CA_{\chi_0+\epsilon\chi}(\ul{z});
\end{equation*}
see Section~\ref{rees_limits} for the notion of such limits. This algebra is
generated by $\CA_{\chi_0}(\ul{z})$ together with
$\Delta^{(k)}(h_{\alpha_i})\in U(\mathfrak{g})^{\otimes k}$, where
$h_{\alpha_i}$ is the coroot corresponding to $\alpha_i$, and $\Delta^{(k)}$ is the $k$ times iterated comultiplication map.

The eigenspaces of $\CA_{\chi_0}(\ul{z})$ decompose into one-dimensional
eigenlines of $\Delta^{(k)}(h_{\alpha_i})$. These eigenlines form finite
strings, and the crystal operators $\on{e}_i$ and $\on{f}_i$ move one step up
and down along these strings, respectively; if there is no such step, the
operator is defined to be zero. The set of eigenlines of the limiting algebra
$\CA_{(\chi_0,\chi)}(\ul{z})$ canonically identifies with
$\CE_\chi(\ul{\lambda})$, so this construction defines a
$\mathfrak{g}$-crystal structure on $\CE_\chi(\ul{\lambda})$.

The following theorem is
\cite[Theorem~12.3]{hkrw}; it provides a geometric realization for the $\mathfrak{g}$-crystal $B_\lambda$.

\th{}\label{irr_cryst_via_Gaudin}
For every dominant $\la$ and regular $\chi \in O$, there is an isomorphism of $\mathfrak{g}$-crystals $\CE_\chi(\la) \simeq B_\la$.
\eth

One can generalize Theorem~\ref{irr_cryst_via_Gaudin} to arbitrary tuples of
dominant weights $\lambda_1,\ldots,\lambda_k$. The following theorem follows
from the results of~\cite{hkrw}.
\th{}\label{tensor_g_via_Gaudin!!!}
For  $z_1 < \ldots < z_k$  and regular $\chi \in O$, there is a canonical  isomorphism of $\mathfrak{g}$-crystals
\begin{equation*}
\CE_\chi(\ul{\la}) \simeq \CE_{\chi}(\la_1) \otimes \ldots \otimes \CE_\chi(\la_k).
\end{equation*}
Consequently,
\begin{equation*}
\CE_\chi(\ul{\la}) \simeq B_{\la_1} \otimes \ldots \otimes B_{\la_k}.
\end{equation*}
\eth

\rem{}
The relation between the Bethe ansatz for the Gaudin spin chain and crystal
bases was first observed by Varchenko~\cite{va} for
$\mathfrak{g}=\mathfrak{sl}_2$. The general idea of indexing solutions of
Bethe ansatz equations by combinatorial data related to Young tableaux goes
back to Kerov, Kirillov, and Reshetikhin~\cite{kkr,kr}.
\erem

\subsection{Bethe subalgebras in the Yangian} 
Gaudin magnet is known to be a degenerate version of the
$XXX$ Heisenberg spin chain. The
integrals of the latter are encoded by Bethe subalgebras in the
Yangian $Y(\mathfrak{g})$, a Hopf algebra deformation of
$U(\mathfrak{g}[t])$; see Section~\ref{yangian} for the definition of
$Y(\mathfrak{g})$. These are commutative subalgebras
$$
B(C)\subset Y(\mathfrak{g}),
\qquad C\in G.
$$
See Section~\ref{bethe} for details.

\ssec{}{Kirillov--Reshetikhin crystals from Bethe subalgebras}\label{KR_via_Bethe_type_A}

Motivated by the classical results of Kirillov and Reshetikhin on the combinatorial description of  asymptotics of solutions of Bethe ansatz equations for the $XXX$ spin chain \cite{kr}, we prove in type $A$ an analogue of Theorem \ref{tensor_g_via_Gaudin!!!} in which Gaudin subalgebras are replaced by Bethe subalgebras in the Yangian and $\mathfrak{g}$-crystals by Kirillov--Reshetikhin $\hat{\mathfrak{g}}$-crystals.

Namely, assume that $G=\on{PGL}_n$ 
so $\mathfrak{g}=\mathfrak{sl}_n$. 
We have the evaluation homomorphism 
\begin{equation*}
{\bf{ev}}_{\ul{z}}\colon Y(\mathfrak{g}) \ra U(\mathfrak{g})^{\otimes k}
\end{equation*} 
that depends on the collection of points $z_1,\ldots,z_k \in \BC$.  

Let $\la_1,\ldots,\la_k$ be a collection of dominant weights of $\mathfrak{g}$. The algebra $B(C)$ acts on $V_{\la_1} \otimes \ldots \otimes V_{\la_k}$ via the homomorphism ${\bf{ev}}_{\ul{z}}$. We denote the resulting $B(C)$-module by
\begin{equation*}
V_{\la_1}(z_1) \otimes \ldots \otimes V_{\la_k}(z_k).
\end{equation*}

Assume that $C$ is a regular element of the maximal compact subtorus $S \subset \on{PGL}_n$, consisting of (classes of) unitary diagonal matrices. By the results of~\cite{imr} (see Section~\ref{simpl_spec_kr}), for generic $z_1,\ldots,z_k$ in the appropriate shifts of $i\BR$ and $\la_j=a_j\varpi_{b_j}$ being multiples of fundamental weights, the action of $B(C)$ on $V_{a_1\varpi_{b_1}}(z_1) \otimes \ldots \otimes V_{a_k\varpi_{b_k}}(z_k)$ has simple spectrum. We denote by $\CE_C(\ul{\la})$ the set of eigenlines.

Consider the covering map 
\begin{equation*}
\on{exp}\colon \mathfrak{h}_{\BR} \ra S,\, \chi \mapsto \on{exp}(2\pi i \chi).    
\end{equation*}
As for the Gaudin case, $S$ is separated  by walls (that are root subtori of $S$). Preimages of walls in $S$ under the covering map $\on{exp}$ induce the decomposition of $\mathfrak{h}_{\BR}$ into  alcoves parametrized by the affine Weyl group,  corresponding to $\mathfrak{g}$ (see Section~\ref{alcoves} for details). Passing to limits $\underset{\longrightarrow}{\on{lim}}\,B(C)$ as $C \ra C_0$ being generic element of the appropriate (affine) wall (see Section~\ref{limits_to_the_wall}), we can endow $\CE_C(\ul{\la})$ with a $\hat{\mathfrak{g}}$-crystal structure (see Section~\ref{cryst_str_E_la} for details).

We then prove the following theorem.
\begin{thmx}\label{KR_viaBethe!!!}{\em{(Theorem~\ref{KR_tensor_via_geometry})}}
(a) For every dominant $\la$ that is a multiple of a fundamental weight there is an isomorphism of $\hat{\mathfrak{g}}$-crystals $\CE_{C}(\la) \simeq {\bf{B}}_\la$, where ${\bf{B}}_\la$ is the Kirillov--Reshetikhin crystal, corresponding to $\la$. 

(b) For a collection $\la_1,\ldots,\la_k$ of dominant weights such that every $\la_i$ is a multiple of a fundamental weight and $z_i$ are as above with $\on{Im}z_1 \gg \on{Im}z_2 \gg \ldots \gg \on{Im}z_k$, we have a canonical isomorphism 
\begin{equation*}
\CE_C(\ul{\la}) \simeq  \CE_C(\la_1) \otimes \CE_C(\la_2) \otimes \ldots \otimes \CE_C(\la_k),  
\end{equation*}
so we have the isomorphism of $\hat{\mathfrak{g}}$-crystals 
\begin{equation*}
\CE_C(\ul{\la}) \simeq {\bf{B}}_{\la_1} \otimes \ldots \otimes {\bf{B}}_{\la_k}.
\end{equation*}
\end{thmx}

\subsection{Monodromy}
Let us now explain one important
application of Theorem \ref{KR_viaBethe!!!}. 

The Deligne--Mumford space $\overline{M}_{0,n+2}$ is a compactification of
$M_{0,n+2}\simeq T^{\mathrm{reg}}$. By the results of \cite{ir0}, the space
$\overline{M}_{0,n+2}$ parametrizes all possible limits of Bethe subalgebras
$B(C)\subset Y(\mathfrak{sl}_n)$. We will be interested in the real form
$\overline{S^{\mathrm{reg}}}\subset \overline{M}_{0,n+2}$ which compactifies
$S^{\mathrm{reg}}$.

It follows from \cite[Theorem 8.3]{ce} and \cite[Theorem 11.12]{ikpr} that the
$S_n$-equivariant fundamental group of $\overline{S^{\mathrm{reg}}}$ is naturally
identified with the {\emph{extended affine cactus group}} $\widetilde{AC}_n$. So, $\widetilde{AC}_n$ acts naturally on the $S_n$-equivariant covering $\mathcal{E}_C(\ul{\lambda})$ above.
Theorem \ref{KR_viaBethe!!!} allows us to uniquely determine the action of the generators $s_{ij}$ of the group $\widetilde{AC}_n$, this result should be considered as an analog of \cite[Proposition 13.7]{hkrw} and confirms \cite[Conjecture~7.1]{imr} and \cite[Conjecture 1.6]{ikr}.

\begin{thmx}\label{th in intro:monodromy}{\emph{(Theorem~\ref{th:monodromy})}}
The generator $s_{ij}$ of the monodromy group of $\overline{S^{\mathrm{reg}}}$ acts on $\mathcal{E}_C(\ul{\lambda})$ as a partial Sch\"utzenberger involution with respect to the corresponding subdiagram in the  Dynkin diagram of $\widetilde{A_{n-1}}$. 
\end{thmx}

Theorem \ref{th in intro:monodromy} applies to quantum cohomology of type $A$ quiver
varieties. Let $\mathfrak{M}=\mathfrak{M}({\mathsf{w}})$ be the disjoint union of type $A$ Nakajima quiver varieties
$\mathfrak{M}(\mathsf{v},\mathsf{w})$ with fixed framing vector ${\mathsf{w}}$. Let ${\mathsf{A}} \times \mathbb{C}^\times$ be the natural torus acting on $\mathfrak{M}$. By the results of Varagnolo \cite{varagnolo_yangians_on_quivers} and
Maulik--Okounkov \cite{mo}, the equivariant cohomology $H^*_{\mathsf{A}\times \mathbb C^\times}(\mathfrak{M})$ 
is identified with a tensor product $V_{\la_1}(z_1) \otimes \ldots \otimes V_{\la_k}(z_k)$ of {\emph{fundamental}} evaluation modules for the
Yangian $Y(\mathfrak{g})$, and the action of the quantum cohomology ring
$QH^*_{\mathsf{A}\times \mathbb C^\times}(\mathfrak{M})$ is identified with the action of the family of
Bethe subalgebras $B(C)$. 
So, the following theorem is a corollary of our results. 
\begin{thmx}\label{intro_thm_monodromy_qh}\emph{(Theorem~\ref{thm_monodromy_qh})}
The eigenvalues of multiplication by elements of the quantum cohomology ring $QH^*_{\mathsf{A} \times \BC^\times}(\mathfrak{M})$ are in bijection with ${\bf{B}}_{\lambda_1} \otimes \ldots \otimes {\bf{B}}_{\lambda_k}$ in such a way that the  monodromy group of these eigenvalues  
is generated by partial Sch\"utzenberger involutions. 
\end{thmx}

Finally, combining Theorem \ref{intro_thm_monodromy_qh} with the results of \cite{rep_th_gies_II} on cross-walling bijections for
type $A$ quiver varieties, we see that, after the appropriate identifications, the
monodromy of $QH^*_{T\times\mathbb C^\times}(\mathfrak{M})$ coincides with the cross-walling
bijection for the symplectically dual quiver variety. This supports the general expectation formulated in the last paragraph of \cite[Section 1.13]{hkrw}.

\ssec{}{Universal Gaudin subalgebras as limits of Bethe subalgebras}\label{gaudin_as_limit_of_bethe!}
One of the main technical ingredients of our proof of Theorem \ref{KR_viaBethe!!!} is the construction of the explicit degeneration from the ``Yangian Bethe'' setting to the ``Gaudin''
setting. 

Namely, for arbitrary simple $\mathfrak{g}$, we define
universal inhomogeneous Gaudin subalgebras
$$
\CA^{\mathrm{u}}_\chi\subset U(\mathfrak{g}[t]),
\qquad \chi\in\mathfrak{g},
$$
and prove that they arise as {\emph{limits}} of Bethe subalgebras; see
Theorem~\ref{th:intro_thm_limit} below. This allows us to explicitly relate our setting with the one considered in \cite{hkrw} and, in particular, show that the crystal structures are compatible.

All results in this subsection hold for arbitrary simple $\mathfrak{g}$ and we believe are of independent interest. We now explain the construction in more detail.

Recall that $\CA_\chi(\ul{z})$ is a  subalgebra of $U(\mathfrak{g})^{\otimes k}$. It turns out that $\CA_{-\chi}(\ul{z})$ can be realized as $\on{ev}_{\ul{z}}(\CA^{\mathrm{u}}_\chi)$, where $\on{ev}_{\ul{z}}\colon U(\mathfrak{g}[t]) \ra U(\mathfrak{g})^{\otimes k}$ is the evaluation homomorphism (at the points $z_1,\ldots,z_k$) and $\CA^{\mathrm{u}}_\chi \subset U(\mathfrak{g}[t])$ is a certain subalgebra of $U(\mathfrak{g}[t])$, which we call the {\em{inhomogeneous universal Gaudin subalgebra}}. The algebra $\CA^{\mathrm{u}}_\chi$ is defined as the image of the Feigin--Frenkel center $\CZ$ in the quantum Hamiltonian reduction $U(\hat{\mathfrak{g}})_{-1/2} /\!\!/\!\!/_{\chi} t^{-1}\mathfrak{g}[[t^{-1}]]$ (see Section~\ref{shift_univ_gaudin_def} for details).

It turns out that the algebra $\CA^{\mathrm{u}}_\chi$ is the limit of Bethe subalgebras in the following sense. Consider the family $B(\on{exp}(\epsilon \chi)) \subset Y(\mathfrak{g})$, $\epsilon \in \BC^\times$. Yangian $Y(\mathfrak{g})$ is equipped with a filtration $F_2$ such that $\on{gr}_{2}Y(\mathfrak{g}) \simeq U(\mathfrak{g}[t])$ (see Section~\ref{two_filtr_yang} for the  discussion of filtrations on $Y(\mathfrak{g})$). We can pass to the Rees algebra $Y_{\hbar}(\mathfrak{g})$, corresponding to the filtration above ($Y_{\hbar}(\mathfrak{g})$ is nothing else but the well-known homogeneous version of the Yangian $Y(\mathfrak{g})$). Then, for every $\epsilon$ as above, we obtain the algebra $Y_{\epsilon}(\mathfrak{g}):=Y_{\hbar}(\mathfrak{g})/(\hbar-\epsilon)$. For $\epsilon \in \BC^\times$, we have the natural isomorphism $Y(\mathfrak{g}) \iso Y_{\epsilon}(\mathfrak{g})$ and denote by $B_{\epsilon}(\on{exp}(\epsilon \chi)) \subset Y_{\epsilon}(\mathfrak{g})$ the image of $B(\on{exp}(\epsilon \chi))$. The following 
Theorem gives a precise relation between Bethe subalgebras in $Y(\mathfrak{g})$ and universal inhomogeneous Gaudin subalgebras of $U(\mathfrak{g}[t])$.

\begin{thmx}\label{th:intro_thm_limit}{\em{(Theorem~\ref{lim_Bethe})}}
We have 
\begin{equation*}
\underset{\epsilon \ra 0}{\on{lim}}\,B_\epsilon(\on{exp}(\epsilon \chi))=\CA^{\mathrm{u}}_\chi.    
\end{equation*}
\end{thmx}

\ssec{}{Structure of the paper} 

The paper is organized as follows. 

In Section~\ref{rees_limits} we discuss several notions of limits of families
of algebras, or more generally of vector spaces, and recall the Rees
construction.

In Section~\ref{univ_inhomogeneous_def_section} we define the
universal inhomogeneous Gaudin subalgebra $\CA^{\mathrm{u}}_\chi \subset U(\mathfrak{g}[t])$ and its classical analogue $\ol{\CA}{}^{\mathrm{u}}_\chi \subset S^\bullet(\mathfrak{g}[t])$. 
In Section \ref{inhomogeneous_gaudin_via_conformal} we study the image of the algebra $\CA^{\mathrm{u}}_\chi$ under the evaluation homomorphism $\on{ev}_{z_1,\ldots,z_k} \colon U(\mathfrak{g}[t]) \to U(\mathfrak{g})^{\otimes k}$ and identify it with the inhomogeneous Gaudin subalgebra of $U(\mathfrak{g})^{\otimes k}$ (see Proposition \ref{realiz_our_via_conf}).

In Section \ref{size_univ_section} we  compute the Poincar\'e series of the algebras $\CA^{\mathrm{u}}_{\chi}$ and $\ol{\CA}{}^{\mathrm{u}}_\chi$ (see Propositions~\ref{poinc_gaud_U},~\ref{poinc_univ_shift}).
In Section  \ref{univ_inhom_as_central_section} we realize  $\CA^{\mathrm{u}}_\chi$ and $\ol{\CA}{}^{\mathrm{u}}_\chi$ as centralizers of certain quadratic elements $\tilde{\Omega}_\chi \in \CA^{\mathrm{u}}_\chi$ and $\Omega_\chi \in \ol{\CA}{}^{\mathrm{u}}_\chi $ (see Proposition~\ref{ca_as_centr}).

Section~\ref{yangian} recalls basic facts about the Yangian $Y(\mathfrak{g})$, its representation theory as well as various 
filtrations on it that are used in the paper.

In Section~\ref{bethe} we recall the definition of  Bethe subalgebras $B(C) \subset Y(\mathfrak{g})$ and then prove Theorem~\ref{lim_Bethe} realizing $\mathcal{A}_\chi^{\mathrm{u}}$ as a limit of Bethe subalgebras.

In Section~\ref{cryst_main} we recall the notion of $\mathfrak{sl}_n$, $\hat{\mathfrak{sl}}_n$- and Kirillov--Reshetikhin crystals. The main result of this section is Proposition~\ref{restr_cryst_class}.

Section~\ref{yangian_type_A} discusses Yangians, Bethe and Gaudin subalgebras in type $A$. In Section~\ref{gen_A_u_chi_exp_resid} we study generators of the universal inhomogeneous Gaudin subalgebras $\CA^{\mathrm{u}}_{\chi} \subset U(\mathfrak{gl}_n[t])$ (see Proposition~\ref{expl_gen_A_u_chi}).  As a corollary, we obtain generators of the inhomogeneous Gaudin subalgebras $\CA_\chi(z_1,\ldots,z_k) \subset U(\mathfrak{gl}_n)^{\otimes k}$ (see Corollary~\ref{expl_gen_A_chi_z} and Proposition~\ref{gen_A_chi_z_via_res}).

 In Section~\ref{alcoves} we recall affine and extended affine Weyl groups of type $A$, alcoves and (affine) walls. In Section~\ref{limits_to_the_wall} we study limits of Bethe and Gaudin subalgebras to generic points  of a wall.

 In Section~\ref{simpl_spec_kr} we formulate the results of~\cite{imr} that give a criterion for the action of Bethe subalgebras on the tensor product $V_{\la_1}(z_1) \otimes \ldots \otimes V_{\la_k}(z_k)$ to have a simple spectrum. In Section \ref{further_prop_of_family} we study the image in $\on{End}(V_{\la_1} \otimes \ldots \otimes V_{\la_k})$ of a certain two-parameter family of algebras.

 In Section~\ref{cryst_str_E_la} we define the structure of $\hat{\mathfrak{sl}}_n$-crystal on the set $\CE_C(\ul{\la})$ of eigenlines of $B(C) \curvearrowright V_{\la_1}(z_1) \otimes \ldots \otimes V_{\la_k}(z_k)$. In Section~\ref{E_la_via_tens_prod} we prove that the crystal $\CE_C(\ul{\la})$ is isomorphic to the tensor product ${\bf{B}}_{\la_1} \otimes \ldots \otimes {\bf{B}}_{\la_k}$ of Kirillov--Reshetikhin crystals.  In Section \ref{monodromy} we apply our results to compute the monodromy of the spectra of Bethe subalgebras and of the quantum cohomology of type $A$ quiver varieties.

\ssec{}{Notation} A complete list of notation may be found in Appendix A.
  
\ssec{}{Acknowledgements}
 We are  grateful to Roman Bezrukavnikov, Aleksei Ilin, Joel Kamnitzer, Vitaly Tarasov and Curtis Wendlandt for useful discussions. L.R. and V. K. were partially supported by the Foundation for the Advancement of Theoretical Physics and Mathematics
 ``BASIS". V. K. was supported by the Simons Foundation Award 888988 as part of the Simons Collaboration on Global Categorical Symmetries.  L.R. was supported by the Fondation Courtois. A significant part of the work was accomplished during L.R.'s stay at the Institut des Hautes Études Scientifiques (IHES). L.R. would like to thank IHES and especially Maxim Kontsevich for the hospitality.

\section{Rees construction and various limits}\label{rees_limits}
\ssec{}{Limits of families of subspaces of a fixed vector space}\label{lim_of_fam_in_vect}
In this section we define limits of families of subspaces of a filtered vector space and discuss some properties of this construction that will be useful later.
Let $U$ be a vector space over $\BC$ equipped with an increasing $\BZ_{\geqslant 0}$-filtration $Q^\bullet U$ by finite-dimensional vector subspaces.  
Let $Z$ be either a formal disc $D=\on{Spec}\BC[[t]]$ (more generally $\on{Spec}R$, where $R$ is a discrete valuation ring) or an affine line $\BA^1$ and let $\mathring{Z}$ be either $\mathring{D}=\on{Spec}\BC((t))$ or $\BA^1 \setminus \{0\}$.
By a point $\epsilon$ of $Z$ we will mean a ${\mathbf{k}}$-point of $Z$, where ${\mathbf{k}}=\BC$ for $Z=\BA^1$ and ${\mathbf{k}}$ is either $\BC$ or $\BC((t))$ for $Z=D$ ($\BC$ or $\on{Frac}R$ for $Z=\on{Spec}R$).
An algebraic family of subspaces $H_\epsilon \subset U \otimes {\bf{k}}$, $\epsilon \in Z$ is a collection of compatible morphisms 
$\mathring{f}_i\colon \mathring{Z} \ra \on{Gr}(d(i),Q^iU)$, $i \in \BZ_{\geqslant 0}$ i.e. a collection of morphisms $\mathring{f}_i$ as above such that for every $\epsilon \in Z$ and $i \leqslant j$ we have $Q^jH_\epsilon \cap (Q^iU \otimes {\bf{k}}) = Q^iH_\epsilon$, where $Q^jH_\epsilon=\mathring{f}_j(\epsilon)$, $Q^iH_\epsilon=\mathring{f}_i(\epsilon)$.

\begin{Rem}
{\em{
In other words, we are given an algebraic family of subspaces $H_\epsilon \subset U \otimes {\bf{k}},\, \epsilon \neq 0$ such that $d(i)=\on{dim}_{\bf{k}}(H_\epsilon \cap Q^iU)$ does not depend on $\epsilon$.}} 
\end{Rem}

Since $\on{Gr}(d(i),Q^iU)$ is a proper variety, we  can uniquely extend each $\mathring{f}_i$ to the map $f_i\colon Z \ra \on{Gr}(d(i),Q^iU)$ and define $\underset{\epsilon \ra 0}{\on{lim}}\,H_\epsilon$ as $\bigcup_i f_i(0) \subset U$. We also set $\underset{\epsilon \ra 0}{\on{lim}}\,Q^i H_\epsilon:=f_i(0)$.

\lem{}\label{lim_in_lim}
Pick $a \in  \underset{\epsilon \ra 0}{\on{lim}}\,Q^i H_\epsilon$.  There exists an algebraic morphism $\tilde{a}\colon Z \ra Q^i U$ such that $\tilde{a}(\epsilon) \in H_\epsilon$ for $\epsilon \neq 0$ and $\tilde{a}(0)=a$.
\elem
\prf
Recall that we have a morphism $f_i\colon Z \ra \on{Gr}(d(i),Q^i U)$, which sends $\epsilon \in Z$ to 
$Q^i H_\epsilon$
and sends $0$ to $\underset{\epsilon \ra 0}{\on{lim}}\,Q^i H_\epsilon$. Let  $\mathcal{V}$ be the tautological vector bundle on $\on{Gr}(d(i),Q^iU)$. Consider the pull-back $f^{*}_i\mathcal{V}$, it is a vector bundle on $Z$. Note that every vector bundle on $Z$ is trivial. Note also that $a$ can be considered as an element of the fiber of $f^{*}_i\mathcal{V}$ over a point $0 \in Z$. Vector bundle $f^{*}_i\mathcal{V}$ is trivial so there exists a section $\tilde{a}\colon Z \ra f^{*}_i\mathcal{V}$ such that $\tilde{a}(0)=a$. Note that we have a natural embedding of $f_{i}^*\mathcal{V}$ into the trivial vector bundle $Z\times Q^i U \ra Z$. In other words, one can consider a section $\tilde{a}(\epsilon)$ as a map $\tilde{a}\colon Z \ra Q^i U$ such that $\tilde{a}(\epsilon) \in Q^i H_\epsilon$ for every $\epsilon \neq 0$ and $\tilde{a}(0)=a$.
\epr

\lem{}\label{lim_can_extend}
Let $\tilde{a}\colon Z \ra Q^iU$ be an algebraic morphism such that $\tilde{a}(\epsilon) \in H_\epsilon$ for $\epsilon \neq 0$ then $\tilde{a}(0) \in \underset{\epsilon \ra 0}{\on{lim}}\,Q^i H_\epsilon$.
\elem
\prf
Since $0 \in \underset{\epsilon \ra 0}{\on{lim}}\,Q^i H_\epsilon$ we can assume that $\tilde{a}(0) \neq 0$. Let $e_1,\ldots,e_{d(i)-1}$ be any subset of $\underset{\epsilon \ra 0}{\on{lim}}\,Q^i H_\epsilon$ such that $\{\tilde{a}(0),e_1,\ldots,e_{d(i)-1}\}$ are linearly independent. 
By Lemma~\ref{lim_in_lim},
we can find morphisms $\tilde{e}_i \colon Z \ra Q^iU$ such that $\tilde{e}_i(\epsilon) \in H_\epsilon$ for $\epsilon \neq 0$ and $\tilde{e}_i(0)=e_i$.
Note that there exists a Zariski open neighbourhood of zero $0 \in W \subset Z$ such that 
the morphism $\tilde{a} \wedge \tilde{e}_1 \wedge \ldots \wedge \tilde{e}_{d(i)-1}\colon W \ra \La^{d(i)}Q^i U$ maps to $\La^{d(i)}Q^i U \setminus \{0\}$ (i.e. for $\epsilon \in W$ elements of the set $\{\tilde{a}(\epsilon),\tilde{e}_1(\epsilon),\ldots,\tilde{e}_{d(i)-1}(\epsilon)\}$ are linearly independent).
Since $\on{dim} Q^i H_\epsilon=d(i)$ and $\{\tilde{a}(\epsilon),\tilde{e}_1(\epsilon),\ldots,\tilde{e}_{d(i)-1}(\epsilon)\} \subset Q^i H_\epsilon$ it follows that the elements $\{\tilde{a}(\epsilon),\tilde{e}_1(\epsilon),\ldots,\tilde{e}_{d(i)-1}(\epsilon)\}$ form a basis of $Q^i H_\epsilon$. 
Recall that we have a closed embedding $\on{Gr}(d(i),Q^iU) \subset \BP(\La^{d(i)}Q^i U)$ and we can consider $Q^i H_\epsilon \in \on{Gr}(d(i),Q^iU)$ as $[\tilde{a}(\epsilon) \wedge \tilde{e}_1(\epsilon) \wedge \ldots \wedge \tilde{e}_{d(i)-1}(\epsilon)] \in \mathbb{P}(\La^{d(i)}Q^iU)$. 
Consider the algebraic morphism 
$[\tilde{a} \wedge \tilde{e}_1 \wedge \ldots \wedge \tilde{e}_{d(i)-1}]\colon W \ra \BP(\La^{d(i)}Q^i U)$, its value at $\epsilon \neq 0$ is $H_\epsilon$ and the value at $\epsilon=0$ is $[\tilde{a}(0) \wedge e_1 \wedge \ldots \wedge e_{d(i)-1}]$.
We conclude that $\underset{\epsilon \ra 0}{\on{lim}}\,Q^i H_\epsilon=[\tilde{a}(0) \wedge e_1 \wedge \ldots \wedge e_{d(i)-1}]$
so $\tilde{a}(0) \in \underset{\epsilon \ra 0}{\on{lim}}\,Q^i H_\epsilon$.
\epr

\cor{}\label{descr_of_lim_elem_level}
The limit $\underset{\epsilon \ra 0}{\on{lim}}\,Q^i H_\epsilon$ can be described as follows: it consists of elements $a \in Q^i U$ such that there exists a morphism $\tilde{a}\colon Z \ra Q^i U$ such that $\tilde{a}(\epsilon) \in H_\epsilon$ for $\epsilon \neq 0$ and $\tilde{a}(0)=a$. 
\ecor
\prf
Follows from Lemmas~\ref{lim_in_lim},~\ref{lim_can_extend}.
\epr

To every ${\bf{k}}$-vector subspace $W \subset U \otimes {\bf{k}}$ we can associate its dimension with respect to the filtration $Q^\bullet$:
\begin{equation*}
\on{dim}_Q W:=\sum_{i \geqslant 0} \on{dim}_{\bf{k}}\left(\frac{W \cap (Q^iU \otimes {\bf{k}})}{W \cap (Q^{i-1}U \otimes {\bf{k}})}\right) q^i \in \BZ[q].    
\end{equation*}
For two series 
$
a(q)=\sum_{i \geqslant 0}a_i q^i,\, b(q)=\sum_{i \geqslant 0}b_i q^i   \in \BZ[q]  
$
we say that $a(q) \geqslant b(q)$ if $\sum_{i=0}^na_i \geqslant \sum_{i=0}^n b_i$ for every $n \in \BZ_{\geqslant 0}$.

\lem{}\label{lim_properties}
$(1)$ We have $\on{dim}_Q H_\epsilon \leqslant \on{dim}_Q (\underset{\epsilon \ra 0}{\on{lim}}\,H_\epsilon)$.

(2) If the filtration $Q^iU$ was induced by some grading with respect to which $H_\epsilon$ are graded then we have $\on{dim}_Q H_\epsilon=\on{dim}_Q (\underset{\epsilon \ra 0}{\on{lim}}\,H_\epsilon)$.
\elem
\prf
Follows from the definitions.
\epr

\rem{}
{\em{Note that it is not true in general that $\on{dim}_Q H_\epsilon=\on{dim}_Q (\underset{\epsilon \ra 0}{\on{lim}}\,H_\epsilon)$. Indeed, take, for example, $U=\BC[x]$ with the filtration by the degree of the polynomial and let $H_\epsilon \subset \BC[x]$ be the subalgebra generated by $\epsilon x^2+x$. Then $\underset{\epsilon \ra 0}{\on{lim}}\,H_\epsilon=\BC[x]$, its dimension is strictly greater than the dimension of $H_1=\BC[x+x^2]$.}}
\erem

\lem{}
$(1)$ If $U$ has an algebra structure such that $Q^iU \cdot Q^jU \subset Q^{i+j}U$ and $H_\epsilon \subset U \otimes {\bf{k}}$ are subalgebras then  $\underset{\epsilon \ra 0}{\on{lim}}\,H_\epsilon$ is a subalgebra of $U$. Moreover, if  $H_\epsilon$ are commutative then $\underset{\epsilon \ra 0}{\on{lim}}\,H_\epsilon$ is commutative. 

$(2)$ If $U$ is itself commutative and equipped with a Poisson bracket $\{\,,\,\}$ 
such that
$H_\epsilon \subset U$ are Poisson commutative subalgebras then $\underset{\epsilon \ra 0}{\on{lim}}\,H_\epsilon$ is Poisson commutative.
\elem
\prf
Let us prove part $(1)$. Pick two elements $a_1\in \underset{\epsilon \ra 0}{\on{lim}}\,Q^iH_\epsilon,\, a_2 \in \underset{\epsilon \ra 0}{\on{lim}}\,Q^jH_\epsilon$. By Lemma~\ref{lim_in_lim} we can find morphisms $\tilde{a}_1\colon Z \ra Q^{i}U,\, \tilde{a}_2 \colon Z \ra Q^{j} U$ such that $\tilde{a}_i(\epsilon) \in H_\epsilon$ for $\epsilon \neq 0$ and $\tilde{a}_i(0)=a_i$. Consider the morphism $\tilde{a}_1\tilde{a}_2\colon Z \ra Q^{i+j}U$ and note that $\tilde{a}_1(\epsilon)\tilde{a}_2(\epsilon) \in H_{\epsilon}$ for $\epsilon \neq 0$ (use that $H_\epsilon \subset U \otimes {\bf{k}}$ is a subalgebra). It follows from Lemma~\ref{lim_can_extend} that $a_1a_2 \in \underset{\epsilon \ra 0}{\on{lim}}\,Q^{i+j}H_\epsilon$. 
Assume now that $H_\epsilon$ are commutative. The composition $[\tilde{a}_1,\tilde{a}_2]\colon Z \ra Q^iU \times Q^jU \ra U$ is clearly continuous.
Note  that $[\tilde{a}_1(\epsilon),\tilde{a}_2(\epsilon)]=0$ for $\epsilon \neq 0$ so we must have $[\tilde{a}_1(0),\tilde{a}_2(0)]=0$.

To prove part $(2)$ consider two elements $a_1 \in \underset{\epsilon \ra 0}{\on{lim}}\,Q^iH_\epsilon,\, a_2 \in \underset{\epsilon \ra 0}{\on{lim}}\,Q^jH_\epsilon$. By Lemma~\ref{lim_in_lim} we can find $\tilde{a}_1\colon Z \ra Q^{i}U,\, \tilde{a}_2 \colon Z \ra Q^{j} U$ such that $\tilde{a}_i(\epsilon) \in H_\epsilon$ for $\epsilon \neq 0$ and $\tilde{a}_i(0)=a_i$.
The composition $Z \ra Q^iU \times Q^jU \ra U$ is clearly continuous.
Note now that $\{\tilde{a}_1(\epsilon),\tilde{a}_2(\epsilon)\}=0$ for $\epsilon \neq 0$ so we must have $\{\tilde{a}_1(0),\tilde{a}_2(0)\}=0$.
\epr

\sssec{}{``Continuous" version of Corollary \ref{descr_of_lim_elem_level}}
In this section we formulate a lemma that should be considered as a ``continuous" analog of Corollary \ref{descr_of_lim_elem_level} above. 
The results of this section will be used in Sections \ref{limits_to_the_wall}, \ref{further_prop_of_family} of the text.

Let $W$ be a finite-dimensional vector space over $\BC$ and $d \in \BZ_{\geqslant 1}$.
Let $(P_n)_{n \in \BZ_{\geqslant 1}}$ be a sequence of points of $\on{Gr}(d,W)$ (considered as a smooth manifold)  and assume that the limit $\underset{n \ra \infty}{\on{lim}}\,P_n$ exists and is equal to some vector space $P \in \on{Gr}(d,W)$.
\lem{}\label{cont_ver_cor}
 Vector space $P$ can be described as follows: it consists of elements $a \in W$ such that there exists a sequence $a_n \in P_n$ with  $\underset{n \ra \infty}{\on{lim}}\,a_n=a$.
\elem
\prf
Let us show that if $a_n \in P_n$ and  $\underset{n \ra \infty}{\on{lim}}\,a_n=a$ then $a \in P$. Consider the natural embedding $\on{Gr}(d,W) \hookrightarrow \BP(\La^dW)$. Consider the tautological line bundle $\CO_{\BP(\La^dW)}(-1)$. It can be trivialized in some neighbourhood of $P$ so we can lift $P_n$, $P$ to some elements $\alpha_n, \alpha \in \La^dW$ such that $\underset{n \ra \infty}{\on{lim}}\,\alpha_n=\alpha$. It follows that $\underset{n \ra \infty}{\on{lim}}\,\alpha_n \wedge a_n = \alpha \wedge a$. Note now that $\alpha_n \wedge a_n=0$ so $\alpha \wedge a=0$, hence, $a \in P$.

Let us show that if $a \in P$ then one can find a sequence $a_n \in P_n$ such that $\underset{n \ra \infty}{\on{lim}}\,a_n=a$. Consider liftings $\al_n$, $\al$ as above. Note that $a \in P$ so there exists a functional $\xi \in (W^{\otimes (d-1)})^*$ such that $\partial_{\xi}(\al)=a$. Set $a_n:=\partial_{\xi}(\al_n)$. It follows from the definitions that $a_n \in P_n$ and $\underset{n \ra \infty}{\on{lim}}\,a_n=a$.
\epr

\ssec{}{Rees construction and limits}\label{rees_constr}
In this section we recall the classical Rees construction 
and discuss its compatibility with taking limits. 
Let $A$ be an algebra equipped with an increasing $\BZ_{\geqslant 0}$-filtration by $\BC$-vector spaces
\begin{equation*}
0=F^{-1}A \subset F^{0}A \subset F^1A \subset \ldots    
\end{equation*}
such that  $F^iA \cdot F^jA \subset F^{i+j}A$ for $i,j \in \BZ_{\geqslant 0}$ and $1 \in F^0A$. Consider the following $\BC[\hbar]$ subalgebra of $A[\hbar]$ 
\begin{equation*}
Rees(A):=\bigoplus_{i \geqslant 0}\hbar^iF^iA. 
\end{equation*}
For $\epsilon \in \BC$ we set $A_{\epsilon}:=Rees(A)/(\hbar-\epsilon)$. More generally, for every ${\bf{k}}$-point $\epsilon \in  \BA^1$ we define $A_\epsilon$ as $A \otimes_{\BC[\hbar]} {\bf{k}}$.
We have a canonical isomorphism of $\BC[\hbar^{\pm 1}]$-algebras:
\begin{equation*}
Rees(A) \otimes_{\BC[\hbar]} \BC[\hbar^{\pm 1}] \iso A[\hbar^{\pm 1}],\, (\hbar^i a) \otimes \hbar^l \mapsto \hbar^{l+i} a,
\end{equation*}
where $i \in \BZ_{\geqslant 0},\, l \in \BZ$ and $a \in F^iA$.

For $\epsilon \neq 0$ we obtain an isomorphism 
\begin{equation*}
A_\epsilon \iso A_1 \iso A,\, [\hbar^ia] \mapsto [\epsilon^i\hbar^ia] \mapsto \epsilon^ia. \end{equation*}

Assume now that the algebra $A$ is equipped with some filtration by $\BC$-vector spaces
\begin{equation*}
0=L^{-1}A \subset L^0 A \subset L^1 A  \subset \ldots \end{equation*}
such that $\on{dim}L^jA<\infty$ for every $j \in \BZ_{\geqslant 0}$. 

We can define a filtration $L^\bullet$ on $Rees(A)$ in the following way:
\begin{equation*}
L^j Rees(A)=\bigoplus_{i \geqslant 0}\hbar^i (F^iA \cap L^jA).    
\end{equation*}
Note that $L^j Rees(A)$ is a $\BC[\hbar]$-submodule of $Rees(A)$ so the filtration $L^j Rees(A)$ induces a filtration on every $A_\epsilon$. Note also that  $L^jRees(A)$ is a finitely generated $\BC[\hbar]$ module: indeed, 
$
L^jRees(A)
$
is a submodule of the $\BC[\hbar]$-module $\bigoplus_{i \geqslant 0}\hbar^i L^jA$ that is free over $\BC[\hbar]$ with any basis of $L^jA$ as the set of generators. Using that $\BC[\hbar]$ is Noetherian, we conclude that 
$L^j Rees(A)$
is finitely generated. 

It follows that $\on{dim}L^iA_\epsilon < \infty$ for every $i \in \BZ_{\geqslant 0}$.
Recall that for every ${\bf{k}}$-point $\epsilon$ of $\BA^1$ and every vector subspace $B_\epsilon \subset A_\epsilon$, we can define the dimension $\on{dim}_L B_\epsilon \in \BZ[q]$ as 
\begin{equation*}
\on{dim}_L B_\epsilon:=\sum_{i \geqslant 0} \on{dim}_{\bf{k}}\left(\frac{L^iA_\epsilon \cap B_\epsilon}{L^{i-1}A_\epsilon \cap B_\epsilon}\right)q^i \in \BZ[q].   
\end{equation*}

Recall that $L^i Rees(A)$ is a finitely generated $\BC[\hbar]$-module i.e. a coherent sheaf on $\BA^1$. Note also that $L^i Rees(A)$ is torsion free finitely generated, hence, is free.  We denote by $(L^i Rees(A))^{\vee}$ the $\BC[\hbar]$-module $\on{Hom}_{\BC[\hbar]}(L^i Rees(A),\BC[\hbar])$. Consider the  spectrum of the relative symmetric power $\on{Spec}(S^\bullet_{\BC[\hbar]}((L^i Rees(A)))^\vee)$ and denote it by $L^iA_{\BA^1}$. Note that we have a natural map $L^iA_{\BA^1} \ra \BA^1$ that is a vector bundle with fiber over $\epsilon \in \BA^1$ being equal to $A_\epsilon$.


Recall that we have a vector bundle $L^iA_{\BA^1} \ra \BA^1$ with fiber over $\epsilon \in \BA^1$ being $L^iA_{\epsilon}$. Consider the relative Grassmannian $\on{Gr}_{\BA^1}(d(i),L^iA_{\BA^1}) \ra \BA^1$, its fiber over $\epsilon \in \BA^1$ is $\on{Gr}(d(i),L^i A_\epsilon)$.
Let $Z$ be a scheme as in Section \ref{lim_of_fam_in_vect} and let us fix a morphism $Z \ra \BA^1$. We denote by $L^i A_{Z} \ra Z$, $\on{Gr}_{Z}(d(i),L^iA_{Z}) \ra Z$ pull backs to $Z$ of the corresponding $\BA^1$-schemes.

Let us fix an algebraic family $B_\epsilon \subset A_\epsilon$ for $\epsilon \neq 0$ such that $d(i)=\on{dim}(B_\epsilon \cap L^iA_\epsilon)$ does not depend on $\epsilon$. In other words, consider a collection of compatible morphisms $\mathring{f}_i\colon \mathring{Z} \ra \on{Gr}_{Z}(d(i),L^iA_{Z})$, then $L^iB_\epsilon=\mathring{f}(\epsilon)$ and $B_\epsilon=\cup_{i \geqslant 0} L^iB_\epsilon$.
By the valuative criterion of properness (applied to $\on{Gr}(d(i),L^i A_{Z}$), $\mathring{Z}$ considered as schemes over $Z$) this section extends uniquely to the section $f\colon Z \ra \on{Gr}_{Z}(d(i),L^iA_{Z})$.
We define $\underset{\epsilon \rightarrow 0}{\on{lim}}\, L^iB_\epsilon$ as $f(0) \subset A_0$. We set $\underset{\epsilon \rightarrow 0}{\on{lim}}\,B_\epsilon := \bigcup_{i \geqslant 0}\underset{\epsilon \rightarrow 0}{\on{lim}}\,L^iB_\epsilon$.

Vector bundle $L^iA_{\BA^1} \ra \BA^1$ is trivial so the following lemma can be proved in the same way as Lemmas~\ref{lim_in_lim},~\ref{lim_can_extend}.
\lem{}\label{lim_descr_via_elements_rees}
The limit $\underset{\epsilon \ra 0}{\on{lim}}\,L^i B_\epsilon$ can be described as follows. Element $a \in L^i A_0$ lies in the limit $\underset{\epsilon \ra 0}{\on{lim}}\,L^i B_\epsilon$ if and only if there exists a section $\tilde{a}\colon Z \ra L^iA_{Z}$ of the vector bundle $L^iA_{Z} \ra Z$ such that $\tilde{a}(\epsilon) \in L^iA_\epsilon$ for $\epsilon \neq 0$ and $\tilde{a}(0)=a$. 
\elem

\lem{}\label{dim_lim_Rees}
We have 
$\on{dim}_L B_1 \leqslant \on{dim}_L(\underset{\epsilon \rightarrow 0}{\on{lim}} B_\epsilon)$.
\elem
\prf
Follows from the definitions. 
\epr

\section{Universal inhomogeneous Gaudin subalgebras: definitions}\label{univ_inhomogeneous_def_section}
\ssec{}{Three filtrations on $S^\bullet(\mathfrak{g}[t]),\, U(\mathfrak{g}[t])$}\label{three_main_filtrations} 

Let \(\fg\) be a finite-dimensional 
simple
Lie algebra. Let 
$(\,,\,)$ be the Killing form on $\mathfrak{g}$.  We fix an orthonormal basis $\{x_a\}|_{a=1,\ldots,\on{dim}\mathfrak{g}}$ of $\mathfrak{g}$. Set $\mathfrak{g}((t^{-1})):=\mathfrak{g} \otimes \BC((t^{-1}))$. For an element $x \in \mathfrak{g}$ and $n \in \BZ$, we set $x[n]:=x \otimes t^n \in \mathfrak{g}((t^{-1})).$ We define a Lie algebra structure on $\mathfrak{g}((t^{-1}))$ as follows: $[x[n],y[m]]:=[x,y][n+m]$, $x,y \in \mathfrak{g}$, $n,m \in \BZ$. Set $\mathfrak{g}[t]:=\mathfrak{g} \otimes \BC[t]$, $t^{-1}\mathfrak{g}[[t^{-1}]]:=\mathfrak{g} \otimes t^{-1}\BC[[t^{-1}]]$. The natural decomposition $\mathfrak{g}((t^{-1}))=\mathfrak{g}[t] \oplus t^{-1}\mathfrak{g}((t^{-1}))$ is the decomposition of Lie subalgebras. 

Let us discuss filtrations on $S^\bullet(\mathfrak{g}[t]),\, U(\mathfrak{g}[t])$ that we will use. 
We have the $PBW$-filtration on $U(\mathfrak{g}[t])$ defined by putting 
\begin{equation*}
\on{deg}_{PBW}x[n]=1.    
\end{equation*}Note that the associated graded $\on{gr}_{PBW}U(\mathfrak{g}[t])$ is isomorphic to the $\BZ_{\geqslant 0}$-graded algebra $S^\bullet(\mathfrak{g}[t])=\bigoplus_{p \geqslant 0} S^p(\mathfrak{g}[t])$. Grading on $S^\bullet(\mathfrak{g}[t])$ induces the filtration that we will also call the $PBW$-filtration on $S^\bullet(\mathfrak{g}[t])$. We will denote by $\on{gr}_{PBW}$ the associated graded with respect to the $PBW$-filtrations on $U(\mathfrak{g}[t])$, $S^\bullet(\mathfrak{g}[t])$.  We also have filtrations $F_1$ on $U(\mathfrak{g}[t]),\, S^\bullet(\mathfrak{g}[t])$ defined by putting $$\on{deg}_1x[n]=n+1.$$ We will denote by $\on{gr}_{1}$ the associated graded with respect to these filtrations.
Finally, we have filtrations $F_2$ on $U(\mathfrak{g}[t]),\, S^\bullet(\mathfrak{g}[t])$ defined by putting
 $$\on{deg}_2x[n]=n.$$
We will denote by $\on{gr}_{2}$ the associated graded with respect to these filtrations.

Note that we have 
\begin{equation*}
\on{dim}F_1^i S^\bullet(\mathfrak{g}[t]),\, 
\on{dim}F_1^i U(\mathfrak{g}[t]) < \infty
\end{equation*}
for every $i \in \BZ_{\geqslant 0}$. So, as in Section~\ref{rees_limits}, for any vector subspace $W$ of $S^\bullet(\mathfrak{g}[t])$ or of $U(\mathfrak{g}[t])$ we can define 
\begin{equation}\label{dim_F_1}
\on{dim}_{F_1}W=\sum_{i \geqslant 0}(\on{dim}F_1^iW) q^i \in \BZ[q],  
\end{equation}
where $F^\bullet_1W$ is the induced filtration on $W$.

\ssec{}{Feigin--Frenkel center and its classical version} 
In this section we recall the Feigin--Frenkel center $\mathcal{Z}$ that is the center of a certain completion of the universal enveloping algebra of the affine Lie algebra $\hat{\mathfrak{g}}$ at the critical level. We also recall  the description of the associated graded $\on{gr}_{PBW}\CZ$. The results of this sections follow from~\cite{bd,ffr,fr,er}.

Recall that 
\begin{equation}\label{def_g_hat}
\hat{\mathfrak{g}}=\mathfrak{g}((t^{-1})) \oplus \BC \mathsf{K},\, [x[n]+b\mathsf{K},y[m]+c\mathsf{K}]=[x,y][n+m]+n(x,y)\delta_{n+m,0}\mathsf{K}
\end{equation}
is the affine Kac-Moody algebra, the central extension of the loop Lie algebra \(\fg((t^{-1}))\).
Define the completion $\widetilde{U}(\hat{\mathfrak{g}})$ of $U(\hat{\mathfrak{g}})$ as the inverse limit of $U(\hat{\mathfrak{g}})/U(\hat{\mathfrak{g}})(t^{-n}\mathfrak{g}[[t^{-1}]])$, $n>1$ and set
\begin{equation*}
\widetilde{U}(\hat{\mathfrak{g}})_{-1/2}:=\widetilde{U}(\hat{\mathfrak{g}})/(\mathsf{K}+1/2).
\end{equation*}

Algebra $\widetilde{U}(\hat{\mathfrak{g}})_{-1/2}$ is equipped with the $PBW$ filtration. The associated graded $\on{gr}_{PBW}\widetilde{U}(\hat{\mathfrak{g}})_{-1/2}$ is isomorphic to the completion: 
\begin{equation*}
\widetilde{S}^\bullet(\mathfrak{g}((t^{-1}))):=\underset{\longleftarrow}{\on{lim}}\,S^{\bullet}(\mathfrak{g}((t^{-1})))/S^{\bullet}(\mathfrak{g}((t^{-1})))(t^{-n}\mathfrak{g}[[t^{-1}]]).
\end{equation*}
Let \(\CZ\) be the center of \(\widetilde{U}(\hat{\fg})_{-1/2}\). Note that $\mathcal{Z}$ contains the following quadratic elements:
\begin{equation*}
S_1^{(r)}=\sum_{a,\,p+q=r} x_a[p]x_a[q]    \in \mathcal{Z} \subset \widetilde{U}(\hat{\mathfrak{g}})_{-1/2},\, r \in \BZ.
\end{equation*}
It is known that there are elements $S^{(r)}_1,\ldots,S^{(r)}_{\on{rk}\mathfrak{g}} \in \CZ,\, r \in \BZ$, such that the image of $\CZ$ in every $U(\hat{\mathfrak{g}})_{-1/2}/U(\hat{\mathfrak{g}})_{-1/2}(t^{-n}\mathfrak{g}[[t^{-1}]])$ is generated by the images of these elements.

There are no simple explicit formulas for the elements $S_i^{(r)}$ for $i>1$, 
but one can describe $\on{gr}_{PBW}S_i^{(r)} \in \widetilde{S}^\bullet(\mathfrak{g}((t^{-1})))$ explicitly.
To do this, let us first describe the associated graded $Z:=\on{gr}_{PBW}\mathcal{Z}$ that is a Poisson-commutative subalgebra in the completion $\widetilde{S}^\bullet(\mathfrak{g}((t^{-1})))$. 

To an element $x \in \mathfrak{g}$ we can associate the following infinite sum
\begin{equation*}
x(z):=\sum_{r \in \BZ}x[r]z^{-r}   
\end{equation*}
that we consider as an element of $\prod_{r \in \BZ}  \widetilde{S}(\mathfrak{g}((t^{-1})))z^r$. We can uniquely extend the map $x \mapsto x(z)$ to the homomorphism of algebras 
\begin{equation*}
S^{\bullet}(\mathfrak{g}) \ra   \prod_{r \in \BZ}  \widetilde{S}(\mathfrak{g}((t^{-1}))) z^r,\, f \mapsto f(z).
\end{equation*}
We can then decompose 
$
f(z)=\sum_{r \in \BZ}f^{(r)}z^{-r}.    
$

Let $\{\Phi_i \in S^\bullet(\mathfrak{g})^{\mathfrak{g}}\}_{i=1,\ldots,\on{rk}\mathfrak{g}}$ be free homogeneous  generators of the commutative algebra $S^\bullet(\mathfrak{g})^{\mathfrak{g}}$ (i.e. $S^\bullet(\mathfrak{g})^{\mathfrak{g}}=\BC[\Phi_i\,|\, i=1,\ldots,\on{rk}\mathfrak{g}]$). 

\prop{}
For $n \in \BZ_{\geqslant 0}$ let $\pi_n$ be the natural surjection from $\widetilde{S}^\bullet(\mathfrak{g}((t^{-1})))$ to $S^{\bullet}(\mathfrak{g}((t^{-1})))/S^{\bullet}(\mathfrak{g}((t^{-1})))(t^{-n}\mathfrak{g}[[t^{-1}]])$.
The algebra $\pi_n(Z)$ is generated by the set 
\begin{equation*}
\{\pi_n(\Phi^{(r)}_i)\,|\, r \in \BZ,1\leqslant i\leqslant \on{rk} \fg, \, \Phi_i \in S^\bullet(\mathfrak{g})^{\mathfrak{g}}\}
\end{equation*}
so the subalgebra  of $Z$ generated by $\Phi^{(r)}_i$ is dense.
We have $\on{gr}_{PBW}S_i^{(r)}=\Phi_i^{(r)}$.
\eprop

Recall  that we have a Casimir element
$
\sum_{a}x_a^2 \in S^{2}(\mathfrak{g})^{\mathfrak{g}} \subset S^\bullet(\mathfrak{g})^{\mathfrak{g}} 
$
that is nothing else but $\Phi_1$.
We conclude that the algebra $Z$ contains the elements 
\begin{equation*}
\Phi_1^{(r)}=\sum_{a,\,p+q=r}x_a[p]x_a[q],\, r \in \BZ.   
\end{equation*}

\ssec{}{Inhomogeneous universal Gaudin subalgebra \(\CA_{\chi}^\mathrm{u}\)}\label{shift_univ_gaudin_def}
Using the center $\CZ$, one can define certain commutative subalgebras of the algebra $U(\mathfrak{g}[t])$.
We regard the Lie algebra \(\fg[t]\) as a ``half'' of the corresponding affine Kac-Moody algebra \(\hat{\fg}\).

Pick $\chi \in \mathfrak{g}$. 
Note that $\chi$ defines a character of the Lie algebra $t^{-1}\fg[[t^{-1}]]=:\hat{\fg}_-$ that sends $x[n]$ to $(\chi,x)\delta_{-1,n}$. We will denote this character by the same letter $\chi$ and will sometimes denote the pairing $(\chi,x)$ by $\chi(x)$.

The image of the natural homomorphism from $\mathcal{Z}$ to the quantum Hamiltonian reduction 
\begin{equation}\label{quant_ham_red}
U(\hat{\mathfrak{g}})_{-1/2} /\!\!/\!\!/_{\chi} t^{-1}\mathfrak{g}[[t^{-1}]]:=\left(U(\hat{\mathfrak{g}})_{-1/2}/U(\hat{\mathfrak{g}})_{-1/2}\{u -\chi(u)\,|\, u \in t^{-1}\mathfrak{g}[[t^{-1}]]\}\right)^{t^{-1}\mathfrak{g}[[t^{-1}]]}   
\end{equation}
 is a commutative subalgebra there.

 Using the decomposition $\hat{\mathfrak{g}}=\mathfrak{g}[t] \oplus t^{-1}\mathfrak{g}[[t^{-1}]] \oplus \BC \mathsf{K}$ and the $PBW$ decomposition, we obtain the embedding of the quantum Hamiltonian reduction algebra~(\ref{quant_ham_red}) into the algebra \(U(\fg[t])\). The image of $\mathcal{Z}$ can be regarded as a commutative subalgebra \(\CA_{\chi}^{\mathrm{u}}\subset U(\fg[t])\), which we call the {\em{inhomogeneous universal Gaudin subalgebra}} of \(U(\fg[t])\).

 Recall the quadratic elements $S_1^{(r)}=\sum_{a,\,p+q=r} x_a[p]x_a[q]$ that lie in $\mathcal{Z}$. We conclude that the algebra $\mathcal{A}^{\mathrm{u}}_\chi$ contains elements 
 \begin{equation*}
\tilde{\omega}_\chi=\frac{1}{2}\sum_a x_a^2+\chi[1],\, \tilde{\Omega}_\chi=\sum_a x_a x_a[1]+\chi[2].     
 \end{equation*}

\rem{}\label{A_chi_cont_omega}
{\em{To see that $\CA_\chi^{\mathrm{u}}$ contains $\tilde{\omega}_\chi,\,\tilde{\Omega}_\chi$ we just note that the image of $S_1^{(0)}$ is 
\begin{equation*}
\sum_a x_a^2+\sum_a 2x_a[1]x_a[-1]=\sum_a x_a^2+2\chi(x_a)x_a=2\tilde{\omega}_\chi.    
\end{equation*}
 The image of $S^{(1)}_1$ is 
 \begin{equation*}
\sum_a 2x_ax_a[1]+\sum_a 2x_a[-1]x_a[2]=\sum_a 2x_ax_a[1]+2\chi(x_a)x_a[2]=2\tilde{\Omega}_\chi.    
\end{equation*}
}}
\erem

\lem{}\label{A_chi_in_inv}
We have $\CA^{\mathrm{u}}_\chi \subset U(\mathfrak{g}[t])^{\mathfrak{z}_{\mathfrak{g}}(\chi)}$.
\elem
\prf
Pick $z \in \CZ$ and $x \in \mathfrak{z}_{\mathfrak{g}}(\chi)$. Note that $[z,x]=0$ considered as elements of the completion $\widetilde{U}(\hat{\mathfrak{g}})_{-1/2}$. We claim that 
\begin{equation}\label{x_inv_quot}
[x,t^{-1}\mathfrak{g}[[t^{-1}]]] \subset \widetilde{U}(\hat{\mathfrak{g}})_{-1/2}\{u - \chi(u)\,|\, u \in t^{-1}\mathfrak{g}[[t^{-1}]]\}.    
\end{equation}
Indeed, if we pick an element $y=y_1[-1]+y_2[-2]+\ldots \in t^{-1}\mathfrak{g}[[t^{-1}]]$ then we have $[x,y]=[x,y_1][-1]+[x,y_2][-2]+\ldots$. It remains to note that $(\chi,[x,y_1])=([\chi,x],y_1)=0$ so~(\ref{x_inv_quot}) indeed holds. It follows that the class of $x$ in $U(\hat{\mathfrak{g}})_{-1/2}/U(\hat{\mathfrak{g}})_{-1/2}\{u-\chi(u)\,|\, u \in t^{-1}\mathfrak{g}[[t^{-1}]]\}$ defines an element of the quantum Hamiltonian reduction $U(\hat{\mathfrak{g}})_{-1/2} /\!\!/\!\!/_{\chi} t^{-1}\mathfrak{g}[[t^{-1}]]$. Now the claim follows from the equality $[z,x]=0$ above.
\epr

\ssec{}{Classical version of the inhomogeneous Gaudin subalgebra} Let us  introduce the classical version $\ol{\CA}{}^{\mathrm{u}}_{\chi}$ of the algebra $\CA^{\mathrm{u}}_\chi$ that will be a Poisson commutative subalgebra of $S^{\bullet}(\mathfrak{g}[t])$.  
We relate $\CA^{\mathrm{u}}_\chi,\, \ol{\CA}{}^{\mathrm{u}}_{\chi}$ in Proposition~\ref{quant_sh_via_class}. 

Recall the Poisson-commutative subalgebra $Z \subset \widetilde{S}^\bullet(\mathfrak{g}((t^{-1})))$. We can consider the Hamiltonian reduction 
\begin{equation}\label{ham_red_symm_classic!}
\left(S^{\bullet}(\mathfrak{g}((t^{-1})))/S^{\bullet}(\mathfrak{g}((t^{-1})))\{u - \chi(u)\,|\, u \in t^{-1}\mathfrak{g}[[t^{-1}]]\}\right)^{t^{-1}\mathfrak{g}[[t^{-1}]]}.
\end{equation}
Again, using the decomposition $\mathfrak{g}((t^{-1}))=\mathfrak{g}[t] \oplus t^{-1}\mathfrak{g}[[t^{-1}]]$, we can embed the Hamiltonian reduction (\ref{ham_red_symm_classic!})  in $S^\bullet(\mathfrak{g}[t])$. We denote by $\ol{\CA}{}^{\mathrm{u}}_{\chi}$ the image of $Z$. In the same way as in Lemma~\ref{A_chi_in_inv} we see that $\ol{\CA}{}^{\mathrm{u}}_{\chi} \subset S^\bullet(\mathfrak{g}[t])^{\mathfrak{z}_{\mathfrak{g}}(\chi)}$.

One can describe explicitly the algebra $\ol{\CA}{}^{\mathrm{u}}_0$, the classical Gaudin subalgebra. It is freely generated by all Fourier components of \(\BC[[t^{-1}]]\)-valued functions \(\Phi_l(x(t))\) on \(t^{-1}\fg[[t^{-1}]]=\on{Spec}S^\bullet(\fg[t])\) for all free homogeneous generators \(\Phi_l\) of the algebra of adjoint invariants \(S^\bullet(\fg)^{\fg}\).
Consider the derivation \(D\) of $S^\bullet(\mathfrak{g}[t])$ given by $D(x[n-1])=nx[n]$.
Recall that \(\Phi_l,l=1,\ldots,\on{rk}\fg\), are free generators of \(S(\fg[0])^{\fg}\subset\Sgt\).

\prop{}\label{generators_of_class_univ_gaud} (\cite[Proposition~4.6]{ir2})
Classical universal Gaudin subalgebra $\ol{\CA}{}^{\mathrm{u}}_0 \subset S^\bullet(\mathfrak{g}[t])^{\mathfrak{g}}$ is the subalgebra freely generated by all \(D^k\Phi_l,k\geqslant 0,l=1,\ldots,\on{rk}\fg\).
\eprop

Recall that the element $\Phi_1^{(r)}=\sum_{a,\, p+q=r}x_a[p]x_a[q]$ lies in $Z$. We see that for $r<-2$, the image of $\Phi_1^{(r)}$ in $S^\bullet(\mathfrak{g}[t])$ is equal to zero. 
The image of $\Phi_1^{(-2)}$ is $\sum_a \chi(x_a)^2=(\chi,\chi)$.
The image of $\Phi_1^{(-1)}$ is $\sum_a \chi(x_a)x_a[0]=\chi[0]$.
The image of $\Phi_1^{(0)}$ is 
\begin{equation*}
\sum_a x_a[0]^2 +\sum_a 2x_a[-1]x_a[1]=\sum_a x_a[0]^2+\sum_a 2\chi(x_a)x_a[1]=
\sum_{a}x_a[0]^2+2\chi[1].
\end{equation*}

The image of $\Phi_1^{(1)}$ is
\begin{multline*}
\sum_a 2x_a[0]x_a[1] +\sum_a 2x_a[-1]x_a[2]=\sum_a 2x_a[0]x_a[1] +\sum_a 2\chi(x_a)x_a[2]=\\
=
\sum_{a}2x_a[0]x_a[1]+2\chi[2].
\end{multline*}
We then set 
\begin{equation*}
\omega_{\chi}:=\frac{1}{2}\sum_{a}x_a[0]^2+\chi[1] \in \ol{\CA}{}^{\mathrm{u}}_{\chi}, 
\end{equation*}
\begin{equation*}
\Omega_\chi:=\sum_a x_a[0]x_a[1]+\chi[2] \in \ol{\CA}{}^{\mathrm{u}}_{\chi}.    
\end{equation*}

\section{Inhomogeneous Gaudin subalgebras and conformal blocks}\label{inhomogeneous_gaudin_via_conformal}
Recall that we have constructed a commutative subalgebra $\CA^{\mathrm{u}}_\chi \subset U(\mathfrak{g}[t])$. 
Note also that for every $k$-tuple of (distinct) points $z_1,\ldots,z_k \in \BC$ we can consider the evaluation homomorphism $\on{ev}_{z_1,\ldots,z_k}\colon U(\mathfrak{g}[t]) \ra U(\mathfrak{g})^{\otimes k}$ and define the commutative algebra $\CA_\chi(z_1,\ldots,z_k):=\on{ev}_{z_1,\ldots,z_k}(\CA^{\mathrm{u}}_{-\chi})$ that is a subalgebra of $U(\mathfrak{g})^{\otimes k}$. 

\begin{warning}{}
Note that $\CA_\chi(z_1,\ldots,z_k)$ is the image of $\CA_{-\chi}^{\mathrm{u}}$, not of $\CA_{\chi}^{\mathrm{u}}$. 
\end{warning}

The goal of this section is to show that the algebra $\CA_\chi(z_1,\ldots,z_k)
$ is the so-called inhomogeneous Gaudin subalgebra 
defined in~\cite{r0},~\cite{fft}, see also~\cite[Section~9]{hkrw}. Let us start from the definition of the inhomogeneous Gaudin subalgebra.

Consider the affine Lie algebra $\hat{\mathfrak{g}}':=\mathfrak{g}((t)) \oplus \BC \mathsf{K}$ and consider the quantum Hamiltonian reduction 
\begin{equation}\label{ham_red_plus}
U({\hat{\mathfrak{g}}'})_{-1/2}/\!\!/\!\!/_0\,\mathfrak{g}[[t]]:=(U(\hat{\mathfrak{g}}')_{-1/2}/U(\hat{\mathfrak{g}}')_{-1/2} \mathfrak{g}[[t]])^{\mathfrak{g}[[t]]}.	
\end{equation}
Using the decomposition $\mathfrak{g}((t))=\mathfrak{g}[[t]] \oplus t^{-1}\mathfrak{g}[t^{-1}]$, we can embed the Hamiltonian reduction (\ref{ham_red_plus}) in the universal enveloping algebra $U(t^{-1}\mathfrak{g}[t^{-1}])$.
Let $\CA' \subset U(t^{-1}\mathfrak{g}[t^{-1}])$ be the image of this embedding, $\CA'$ is a commutative subalgebra of $U(t^{-1}\mathfrak{g}[t^{-1}])$. 

\rem{}
{\em{It follows from~\cite{ff},~\cite{fr} that the center $\mathcal{Z}$ maps onto the Hamiltonian reduction~(\ref{ham_red_plus}) so there is no need to mention $\CZ$ in the  definition of the algebra $\CA'$.
}}
\erem
Note now that for nonzero  $z_i$ we can consider the evaluation homomorphism 
$U(t^{-1}\mathfrak{g}[t^{-1}]) \xrightarrow{\on{ev}'_{z_1,\ldots,z_k}} U(\mathfrak{g})^{\otimes k}$. We also have the  ``evaluation at the infinity''  homomorphism $\on{ev}'_{\infty}\colon U(t^{-1}\mathfrak{g}[t^{-1}]) \ra S^\bullet(\mathfrak{g})$ that is induced by the homomorphism $t^{-1}\mathfrak{g}[t] \ra \mathfrak{g}$ that extracts the coefficient in front of $t^{-1}$. We obtain the homomorphism
\begin{equation*}
\on{ev}'_{z_1,\ldots,z_k,\infty}:=\on{ev}'_{z_1,\ldots,z_k} \otimes \on{ev}'_{\infty}\colon U(t^{-1}\mathfrak{g}[t^{-1}]) \ra U(\mathfrak{g})^{\otimes k}\otimes S^\bullet(\mathfrak{g}).	
\end{equation*}
Recall that we are given $\chi \in \mathfrak{g}$ that defines the evaluation at $\chi$ homomorphism  $S^\bullet(\mathfrak{g}) \ra \BC$ (after the identification $\mathfrak{g} \simeq \mathfrak{g}^*$ via the $\mathfrak{g}$-invariant bilinear form $(\,,\,)$). We obtain the composite map to be denoted 
\begin{equation*}
\on{ev}'_{z_1,\ldots,z_k,\chi}\colon U(t^{-1}\mathfrak{g}[t^{-1}])	\ra U(\mathfrak{g})^{\otimes k} \otimes S^\bullet(\mathfrak{g}) \ra U(\mathfrak{g})^{\otimes k}. 
\end{equation*}
Let $\CA'_{\chi}(z_1,\ldots,z_k)$ be the image of $\CA'$ under the homomorphism $\on{ev}'_{w-z_1,\ldots,w-z_k,\chi}(\CA')$, $w \in \BC \setminus \{z_1,\ldots,z_k\}$. 

\rem{}
{\em{Note that $\CA'_{\chi}(z_1,\ldots,z_k)=\on{ev}_{z_1,\ldots,z_k,-\chi}(\CA')$.}}
\erem

Our goal is to show that 
\begin{equation}\label{identific_two_gaudin}
\CA_\chi(z_1,\ldots,z_k)=\CA'_\chi(z_1,\ldots,z_k).
\end{equation}

Assume that we are given some $\mathfrak{g}$-modules $M_{1},\ldots,M_{k}$. Then the algebras $\CA_\chi(z_1,\ldots,z_k),\, \CA_\chi'(z_1,\ldots,z_k) \subset U(\mathfrak{g})^{\otimes k}$ act naturally on the tensor product $M_{1} \otimes \ldots \otimes M_{k}$ and their images in $\on{End}(M_{1} \otimes \ldots \otimes M_{k})$ are commutative subalgebras that we denote by $\CA_\chi(M_{1},\ldots,M_{k}),\, \CA'_\chi(M_{1},\ldots,M_{k})$ respectively.
To see that $\CA_\chi(z_1,\ldots,z_k)=\CA'_\chi(z_1,\ldots,z_k)$ it is enough to show that
\begin{equation}\label{eq_a-a'_goal_prove}
\CA_\chi(M_{1},\ldots,M_{k})=\CA'_\chi(M_{1},\ldots,M_{k})\end{equation}
for every $M_{1},\ldots,M_{k}$. 
To prove (\ref{eq_a-a'_goal_prove}) we will recall the general approach that is used to construct commutative subalgebras of such sort (approach via conformal blocks) and will identify both $\CA_\chi(M_{1},\ldots,M_{k})$ and $\CA'_\chi(M_{1},\ldots,M_{k})$ with a well-known commutative subalgebra of $\on{End}(M_{1} \otimes \ldots \otimes M_{k})$ (see Proposition~\ref{realiz_our_via_conf}).


We start with some general recollections about conformal blocks (following~\cite{ffr},~\cite{fft} and~\cite{ff}).
Recall that we fix distinct points $z_{1},\ldots,z_{k} \in \BP^1$. In the neighbourhood of each point we have a local coordinate $t-z_i$. Set $\tilde{\mathfrak{g}}(z_i):=\mathfrak{g}((t-z_i))$ and let $\hat{\mathfrak{g}}(z_i)=\tilde{\mathfrak{g}}(z_i) \oplus \BC {\mathsf{K}}_i$ be the one-dimensional central extension of $\tilde{\mathfrak{g}}(z_i)$ (see (\ref{def_g_hat})).
We set $\tilde{\mathfrak{g}}(\ul{z}):=\bigoplus_{i=1}^k \tilde{\mathfrak{g}}(z_i)$ and denote by $\hat{\mathfrak{g}}(\ul{z})$ the one-dimensional central extension of $\tilde{\mathfrak{g}}(\ul{z})$   that is $\Big(\bigoplus_{i=1}^k \hat{\mathfrak{g}}(z_i)\Big)/({\mathsf{K}}_i-{\mathsf{K}}_j\,|\, i \neq j)$.  Let ${\mathsf{K}} \in \hat{\mathfrak{g}}(\ul{z})$ be the central element (class of any ${\mathsf{K}}_i$).

Let $M_{1},\ldots,M_{k}$ be any $\mathfrak{g}$-modules.  
We consider $M_{i}$ as a module over $\mathfrak{g}[[t-z_i]] \oplus \BC \mathsf{K}$, where $(t-z_i)\mathfrak{g}[[t-u_i]]$ acts trivially, $\mathsf{K}$ acts via the multiplication by $-1/2$ and $\mathfrak{g}$ acts via its given action on $M_{i}$. 
Consider also the induced modules $\mathbb{M}_{z_i}:=\on{Ind}_{\mathfrak{g}[[t-z_i]] \oplus \BC \mathsf{K}}^{\hat{\mathfrak{g}}(z_i)}M_{i}$. For $z \in \BP^1$ we set
$\mathbb{V}_{0,z}:=\on{Ind}_{\mathfrak{g}[[t-z]] \oplus \BC \mathsf{K}}^{\hat{\mathfrak{g}}(z)} \BC$. 

\rem{} 
{\em{Note that $\mathbb{V}_{0,z}$ is nothing else than $\mathbb{M}_{z}$ for $M=\BC$, the trivial $\mathfrak{g}$-module.
}}
\erem

Consider  the tensor product $\mathbb{M}_{z_1} \otimes \ldots \otimes \mathbb{M}_{z_k}$. Let $\mathfrak{g}_{\ul{z}}$ be the Lie algebra of regular functions on $\mathbb{P}^1 \setminus \{z_1,\ldots,z_k\}$. 
We have an embedding $\mathfrak{g}_{\ul{z}} \subset \hat{\mathfrak{g}}(\ul{z})$ and denote by $H(\mathbb{M}_{z_1},\ldots,\mathbb{M}_{z_k})$ the space of coinvariants 
\begin{equation*}
H(\mathbb{M}_{z_1},\ldots,\mathbb{M}_{z_k}):=(\mathbb{M}_{z_1} \otimes \ldots \otimes \mathbb{M}_{z_k})/{\mathfrak{g}_{\ul{z}}}.
\end{equation*}

The following proposition is standard.
\prop{}\label{coinv_loop_via_coinv}
The embedding $M_{1} \otimes \ldots \otimes M_{k} \subset \mathbb{M}_{z_1} \otimes \ldots \otimes \mathbb{M}_{z_k}$ induces the isomorphism $(M_{1} \otimes \ldots \otimes M_{k})/{\mathfrak{g}} \iso H(\mathbb{M}_{z_1},\ldots,\mathbb{M}_{z_k})$. 
\eprop

\cor{}
Pick $z \in \BP^1 \setminus \{z_1,\ldots,z_k\}$.
The embedding $\mathbb{M}_{z_1} \otimes \ldots \otimes \mathbb{M}_{z_k} \subset \mathbb{M}_{z_1} \otimes \ldots \otimes \mathbb{M}_{z_k} \otimes \mathbb{V}_{0,z}$ induces the isomorphism $H(\mathbb{M}_{z_1},\ldots,\mathbb{M}_{z_k}) \iso H(\mathbb{M}_{z_1},\ldots,\mathbb{M}_{z_k},\mathbb{V}_{0,z})$.
\ecor

We now define the commutative subalgebra 
\begin{equation*}
\CA(\mathbb{M}_{z_1},\ldots,\mathbb{M}_{z_k}) \subset \on{End}((M_{1} \otimes \ldots \otimes M_{k})/{\mathfrak{g}})
\end{equation*}
as follows. We consider the natural action of $\CZ^{\otimes k}$ on $H({\mathbb{M}}_{z_1},\ldots,{\mathbb{M}}_{z_k}) \simeq (M_1 \otimes \ldots \otimes M_k)/\mathfrak{g}$ and denote by  $\CA(\mathbb{M}_{z_1},\ldots,\mathbb{M}_{z_k})$ the image of $\CZ^{\otimes k}$ in 
$\on{End}((M_{1} \otimes \ldots \otimes M_{k})/{\mathfrak{g}})$.

\prop{}\label{different_realis_coinv_alg} 
Consider any nonempty subset $\{i_1,\ldots,i_l\} \subset \{1,2,\ldots,k\}$, $l \in \{1,2,\ldots,k\}$ and consider the image of the natural homomorphism 
\begin{equation*}
1 \otimes \ldots \otimes 1 \otimes \underset{i_1}{\mathcal{Z}} \otimes 1 \ldots \otimes 1 \otimes \underset{i_2}{\mathcal{Z}} \otimes 1 \ldots \otimes 1 \otimes \underset{i_l}{\mathcal{Z}} \otimes 1 \ldots \otimes 1 \ra \on{End}(H(\mathbb{M}_{z_1},\ldots,\mathbb{M}_{z_k})).
\end{equation*}
This image coincides with $\CA(\mathbb{M}_{z_1},\ldots,\mathbb{M}_{z_k})$. In particular, the algebra $\CA(\mathbb{M}_{z_1},\ldots,\mathbb{M}_{z_k})$ coincides with the image of the natural homomorphism 
\begin{equation*}
\underbrace{1 \otimes 1 \otimes \ldots \otimes 1}_{k} \otimes \CZ \ra \on{End}(H(\mathbb{M}_{z_1},\ldots,\mathbb{M}_{z_k},\mathbb{V}_{0,z}))	
\end{equation*}
and does not depend on the point $z \in \BP^1 \setminus \{z_1,\ldots,z_k\}$.
\eprop
    \prf 
    Set $\mathring{D}:=\on{Spec}\BC((t))$, for $z \in \BP^1$ we set $\mathring{D}_z:=\on{Spec}\BC((t-z))$ i.e. $\mathring{D}_z$ is the punctured formal neighbourhood of $z \in \BP^1$. Let $X$ be a smooth curve or $\mathring{D}_z$. Let $\on{Op}(X)$ be the moduli space of $G^\vee$-opers (see, for example, \cite[Section 4.2]{fr_loop} for the definition), here $G^\vee$ is the Langlands dual group to $G$.
    Recall now that by \cite{ff}, \cite{fr2} (see also \cite{fr_loop}) we have the natural identification $\CZ \simeq \CO(\on{Op}(\mathring{D}))$ which induces the identification $\CZ^{\otimes k} \simeq \CO(\on{Op}(\mathring{D}_{z_1}) \times \ldots \times \on{Op}(\mathring{D}_{z_k}))$. We have the homomorphism 
    \begin{equation}\label{ext_to_glob}
    \CO(\on{Op}(\mathring{D}_{z_1}) \times \ldots \times \on{Op}(\mathring{D}_{z_k})) \ra  \CO(\on{Op}(\BP^1 \setminus \{z_1,\ldots,z_k\}))
    \end{equation}
    induced by the natural (restriction) morphisms $\on{Op}(\BP^1 \setminus \{z_1,\ldots,z_k\}) \ra  \on{Op}(\mathring{D}_{z_i})$, $i=1,\ldots,k$. It follows from the definitions (see the proof of \cite[Theorem 5.7]{fft} for details) that the action of $\CZ^{\otimes k}$ on $H(\mathbb{M}_{z_1},\ldots,\mathbb{M}_{z_k})$ factors through (\ref{ext_to_glob}). It remains to note that for every $i=1,\ldots,k$, the homomorphism $\CO(\on{Op}(\mathring{D}_{z_i})) \ra \CO(\on{Op}(\BP^1 \setminus \{z_1,\ldots,z_k\}))$ is surjective (since the corresponding morphism $\on{Op}(\BP^1 \setminus \{z_1,\ldots,z_k\}) \ra \on{Op}(\mathring{D}_{z_i})$ is a closed embedding).

\epr

Let $\BC_\chi$ be the one-dimensional module over $t^{-1}\mathfrak{g}[[t^{-1}]] \oplus \BC\mathsf{K}$, where $t^{-1}\mathfrak{g}[[t^{-1}]]$ acts via the character $\chi$ and $\mathsf{K}$ acts via $-1/2$.
Consider the following module: $\mathbb{I}_{\chi,\infty}=\on{Ind}_{t^{-1}\mathfrak{g}[[t^{-1}]] \oplus \BC \mathsf{K}}^{\hat{\mathfrak{g}}} \BC_\chi=\on{Ind}_{\mathfrak{g}[[t^{-1}]] \oplus \BC \mathsf{K}}^{\hat{\mathfrak{g}}}I_\chi$, where $I_{\chi}=\on{Ind}_{t^{-1}\mathfrak{g}[[t^{-1}]] \oplus \BC \mathsf{K}}^{\mathfrak{g}[[t^{-1}]] \oplus \BC \mathsf{K}} \BC_\chi$. 

\rem{}
{\em{Note that we can realize $\mathbb{I}_{\chi,\infty}$ as $\mathbb{M}_z$ for an appropriate choice of a $\mathfrak{g}$-module $M$ and a point $z \in \BP^1$. Indeed, taking $M=I_\chi$ and $z=\infty \in \BP^1$, we see  that $\mathbb{I}_{\chi,\infty}=\mathbb{M}_\infty$.}}
\erem

\lem{}
Assume that $\{z_1,\ldots,z_k\} \subset \BA^1=\BP^1 \setminus \{\infty\}$.
The natural embedding $M_{1} \otimes \ldots \otimes M_{k} \subset \mathbb{M}_{z_1} \otimes \ldots \otimes \mathbb{M}_{z_k} \otimes \mathbb{I}_{\chi,\infty}$ induces the isomorphism $M_{1} \otimes \ldots \otimes M_{k} \iso H(\mathbb{M}_{z_1},\ldots,\mathbb{M}_{z_k},\mathbb{I}_{\chi,\infty})$.
\elem
\prf
Follows from Proposition~\ref{coinv_loop_via_coinv} using that $I_\chi$ is isomorphic to $U(\mathfrak{g})$ as a $\mathfrak{g}$-module.
\epr

\prop{}\label{realiz_our_via_conf}
We 
have 
\begin{equation*}
\CA_\chi(z_1,\ldots,z_k) = \CA(\mathbb{M}_{z_1},\ldots,\mathbb{M}_{z_k}, \mathbb{I}_{-\chi,\infty}) = \CA'_\chi(z_1,\ldots,z_k).
\end{equation*}
In particular, 
\begin{equation*}
\CA_\chi(z_1,\ldots,z_k)= \CA'_\chi(z_1,\ldots,z_k).\end{equation*}
\eprop
\prf
Recall that the algebra $\CA^{\mathrm{u}}_\chi$ is the image of $\CZ$ in the quantum Hamiltonian reduction~(\ref{quant_ham_red}) and the latter identifies naturally with 
$\on{End}_{\hat{\mathfrak{g}}}(\mathbb{I}_{-\chi,\infty}) \subset \on{End}_{\mathfrak{g}[t]}(\mathbb{I}_{-\chi,\infty})=U(\mathfrak{g}[t])^{\mathrm{opp}}$. Pick $X \in U(\mathfrak{g}[t])$ and let us denote by $X_{z_i} \in U(\mathfrak{g}[[t-z_i]])$ the corresponding element of $U(\mathfrak{g}[[t-z_i]])$.
Note that for every $v_i \in \mathbb{M}_{z_i},\, v_\infty \in \mathbb{I}_{-\chi,\infty}$ the following equality holds in $H(\mathbb{M}_{z_1},\ldots,\mathbb{M}_{z_k},\mathbb{I}_{-\chi,\infty})$:
\begin{equation}\label{eq_in_H}
\Big[\Big(\bigotimes_{i=1}^k v_i \Big) \otimes X(v_\infty)+\sum_{i=1}^k v_1 \otimes \ldots \otimes v_{i-1} \otimes X_{z_i}(v_i) \otimes v_{i+1} \otimes \ldots \otimes v_k \otimes v_{\infty}	 \Big]=0.
\end{equation}
Now, taking $v_i \in M_{i} \subset \mathbb{M}_{z_i}$, we see that $X_{z_i}(v_i)=X(z_i)v_i$ so we conclude that
\begin{equation}\label{ev_via_act}
\sum_{i=1}^k v_1 \otimes \ldots \otimes v_{i-1} \otimes X_{z_i}(v_i) \otimes v_{i+1} \otimes \ldots \otimes v_k=\on{ev}_{z_1,\ldots,z_{k}}(X)(v_1 \otimes \ldots \otimes v_k).
\end{equation}
Equations~(\ref{eq_in_H}),~(\ref{ev_via_act}) imply that in $H(\mathbb{M}_{z_1},\ldots,\mathbb{M}_{z_k},\mathbb{I}_{-\chi,\infty})$ we have 
\begin{equation*}
\Big[\on{ev}_{z_1,\ldots,z_{k}}(X)(v_1 \otimes \ldots \otimes v_k) \otimes v_{\infty}\Big]=- \Big[\Big(\bigotimes_{i=1}^k v_i \Big) \otimes X(v_\infty)\Big]   
\end{equation*}
and the equality $\CA_\chi(z_1,\ldots,z_k)=\CA(\mathbb{M}_{z_1},\ldots,\mathbb{M}_{z_k},\mathbb{I}_{-\chi,\infty})$ then follows.

Let us now prove the equality $\CA'_\chi(z_1,\ldots,z_k)=\CA(\mathbb{M}_{z_1},\ldots,\mathbb{M}_{z_k},\mathbb{I}_{-\chi,\infty})$. Recall that by Proposition~\ref{different_realis_coinv_alg} the algebra $\CA(\mathbb{M}_{z_1},\ldots,\mathbb{M}_{z_k},\mathbb{I}_{-\chi,\infty})$ coincides with the image of $1 \otimes \ldots \otimes 1 \otimes \CZ$ in $\on{End}(H(\mathbb{M}_{z_1},\ldots,\mathbb{M}_{z_k},\mathbb{I}_{-\chi,\infty},\mathbb{V}_{0,z}))$. Directly from the definitions (cf.~\cite[Sections 2.7,~2.8]{fft}) it then follows that the image of $\underbrace{1 \otimes \ldots \otimes 1}_{k+1} \otimes \CZ$ in $\on{End}(H(\mathbb{M}_{z_1},\ldots,\mathbb{M}_{z_k},\mathbb{I}_{-\chi,\infty},\mathbb{V}_{0,z}))$ is exactly $\CA'_\chi(z_1,\ldots,z_k)$.
\epr

\rem{}
{\em{In type $A$ the equality $\CA_\chi(z_1,\ldots,z_k)=\CA'_\chi(z_1,\ldots,z_k)$ can be deduced from the explicit description of the generating functions for the generators of the algebras $\CA_\chi(z_1,\ldots,z_k),\, \CA'_\chi(z_1,\ldots,z_k)$ (see Section~\ref{gen_A_u_chi_exp_resid}, Corollary~\ref{expl_gen_A_chi_z} and~\cite[Theorem~3.1]{chemo}).}}
\erem

\section{The size of $\ol{\CA}{}^{\mathrm{u}}_{\chi},\, \CA^{\mathrm{u}}_\chi$}\label{size_univ_section}
The goal of this section is to compute dimensions of $\ol{\CA}{}^{\mathrm{u}}_{\chi},\, \CA^{\mathrm{u}}_\chi$ with respect to the filtration $F_1$ (see Section~\ref{three_main_filtrations}). To do this, we first describe the associated graded subalgebras $\on{gr}_{PBW}\ol{\CA}{}^{\mathrm{u}}_{\chi},\, \on{gr}_{PBW}\CA^{\mathrm{u}}_\chi \subset S^\bullet(\mathfrak{g}[t])$. It turns out that both of them are isomorphic to the tensor product $\ol{\CA}{}^{\mathrm{u}}_{0} \otimes_{S^\bullet(\mathfrak{g})^{\mathfrak{g}}} \ol{\CA}_\chi$ (see Proposition~\ref{gr_to_univ}), where $\ol{\CA}_\chi \subset S^\bullet(\mathfrak{g})$ is the subalgebra of $S^\bullet(\mathfrak{g})$ generated by $\frac{\partial^k \Phi_l}{\partial^k \chi}$, $l=1,\ldots,\on{rk}\mathfrak{g}$, $k \geqslant 0$.

\begin{Rem}
{\em{Algebra $\ol{\CA}_\chi$ is called the shift of argument subalgebra (or Mishchenko-Fomenko subalgebra). It can be considered as a classical version of the algebra $\CA_\chi(z)$ defined in Section \ref{inhomogeneous_gaudin_via_conformal}.}}
\end{Rem}

We then recall the dimensions (see (\ref{dim_F_1})) of $\ol{\CA}{}^{\mathrm{u}}_{0},\, \ol{\CA}_\chi$  and obtain the desired formula for the dimension of $\ol{\CA}{}^{\mathrm{u}}_{\chi},\, \CA^{\mathrm{u}}_\chi$. Moreover, we conclude  that $\on{dim}_{F_1}\ol{\CA}{}^{\mathrm{u}}_{\chi}= \on{dim}_{F_1}\ol{\CA}{}^{\mathrm{u}}_0(\mathfrak{z}_{\mathfrak{g}}(\chi))$, where $\ol{\CA}{}^{\mathrm{u}}_0(\mathfrak{z}_{\mathfrak{g}}(\chi))$ is the universal (classical) Gaudin subalgebra of $S^\bullet(\mathfrak{z}_{\mathfrak{g}}(\chi))$ (see \cite[Remark in Section 4.7]{ir2}). This equality will be important in the next section.

Set $\mathfrak{z}_{\mathfrak{g}}(\chi)^{\mathrm{der}}:=[\mathfrak{z}_{\mathfrak{g}}(\chi),\mathfrak{z}_{\mathfrak{g}}(\chi)]$. Let $\chi_1$ be a regular Cartan element of $\mathfrak{z}_{\mathfrak{g}}(\chi)^{\mathrm{der}}$.
Let $\ol{\CA}_{(\chi,\chi_1)}$ be the subalgebra of $S^\bullet(\mathfrak{g})$ generated by $\ol{\CA}_{\chi},\, \ol{\CA}_{\chi_1}(\mathfrak{z}_{\mathfrak{g}}(\chi))$.
According to ~\cite{sh,t}, the following holds.
\prop{}\label{ext_A_chi} 
(a) The subalgebra $\ol{\CA}_{(\chi,\chi_1)} \subset S^\bullet(\mathfrak{g})$
is a free polynomial algebra with
$\frac{1}{2}(\on{dim}\mathfrak{g}+\on{rk}\mathfrak{g})$ generators. The set of generators is the union of standard
generators of $\ol{\CA}_{\chi_1}(\mathfrak{z}_{\mathfrak{g}}(\chi))$ and $\partial^k_\chi(\Phi_i)$ with $k = 1,\ldots, d_i'$ for some $d_i' \leqslant d_i$.

$(b)$ The subalgebra $\ol{\CA}_{(\chi,\chi_1)} \subset S^\bullet(\mathfrak{g})$ is maximal Poisson-commutative and is equal to the limit $\underset{\epsilon \rightarrow 0}{\on{lim}}\, \ol{\CA}_{\chi+\epsilon\chi_1}$, where the limit is taken with respect to the filtration on $S^\bullet(\mathfrak{g})$ by the degree.
\eprop

\prop{}\label{A_chi_max}
For  $\chi \in \mathfrak{g}$, the Poisson subalgebra $\ol{\CA}_{\chi} \subset S^{\bullet}(\mathfrak{g})$ is maximal  Poisson-commutative in $S^{\bullet}(\mathfrak{g})^{\mathfrak{z}_{\mathfrak{g}}(\chi)}$.
\eprop
\prf
This follows from~\cite[Corollary~9.9]{hkrw}. We give a sketch of the argument. Let $\chi_1$ be a regular element of $\mathfrak{z}_{\mathfrak{g}}(\chi)$.
Recall the subalgebra $\ol{\CA}_{(\chi,\chi_1)} \subset S^\bullet(\mathfrak{g})$ generated by $\ol{\CA}_\chi$ and $\ol{\CA}_{\chi_1}(\mathfrak{z}_{\mathfrak{g}}(\chi)) \subset S^\bullet(\mathfrak{z}_{\mathfrak{g}}(\chi))$.

Pick an element $x \in S^{\bullet}(\mathfrak{g})^{\mathfrak{z}_{\mathfrak{g}}(\chi)}$ that commutes with $\ol{\CA}_{\chi}$. Since $x$ is $\mathfrak{z}_{\mathfrak{g}}(\chi)$-invariant it follows that $x$ commutes with $\ol{\CA}_{\chi_1}(\mathfrak{z}_{\mathfrak{g}}(\chi))$ so it commutes with $\ol{\CA}_{(\chi,\chi_1)}$, hence, lies in $\ol{\CA}_{(\chi,\chi_1)}$ (use part $(b)$ of Proposition~\ref{ext_A_chi}). Recall now that by the result of Knop from~\cite{kn} we have
\begin{equation*}
S^\bullet(\mathfrak{g})^{\mathfrak{z}_\mathfrak{g}(\chi)} \cdot S^\bullet(\mathfrak{z}_\mathfrak{g}(\chi))=S^\bullet(\mathfrak{g})^{\mathfrak{z}_\mathfrak{g}(\chi)} \otimes_{S^\bullet(\mathfrak{z}_\mathfrak{g}(\chi))^{\mathfrak{z}_\mathfrak{g}(\chi)}} S^\bullet(\mathfrak{z}_\mathfrak{g}(\chi))
\end{equation*}
so 
\begin{equation*}
\ol{\CA}_{(\chi,\chi_1)}=\ol{\CA}_\chi \cdot \ol{\CA}_{\chi_1}(\mathfrak{z}_{\mathfrak{g}}(\chi))=\ol{\CA}_\chi \otimes_{S^\bullet(\mathfrak{z}_\mathfrak{g}(\chi))^{\mathfrak{z}_\mathfrak{g}(\chi)}}  \ol{\CA}_{\chi_1}(\mathfrak{z}_{\mathfrak{g}}(\chi)).
\end{equation*}
It now follows from the fact that $\ol{\CA}_{\chi,\chi_1}$ is freely generated by the standard generators of $\ol{\CA}_{\chi_1}(\mathfrak{z}_{\mathfrak{g}}(\chi))$ together with some $\partial^k_{\chi} \Phi_l \in \ol{\CA}_{\chi}$ (see Proposition~\ref{ext_A_chi}) that $S^\bullet(\mathfrak{g})^{\mathfrak{z}_{\mathfrak{g}}(\chi)} \cap \ol{\CA}_{(\chi,\chi_1)}=\ol{\CA}_\chi$. We conclude that $x \in S^\bullet(\mathfrak{g})^{\mathfrak{z}_{\mathfrak{g}}(\chi)} \cap \ol{\CA}_{(\chi,\chi_1)}=\ol{\CA}_\chi$ as desired.
\epr

\prop{}\label{prod_inv_funct}
Consider the subalgebras $S^\bullet(\mathfrak{g}[t])^{\mathfrak{g}},\, S^\bullet(\mathfrak{g}) \subset S^\bullet(\mathfrak{g}[t])$. We have 
\begin{equation*}
S^{\bullet}(\mathfrak{g}[t])^{\mathfrak{g}}\cdot S^\bullet(\mathfrak{g})=S^\bullet(\mathfrak{g}[t])^{\mathfrak{g}} \otimes_{S^\bullet(\mathfrak{g})^{\mathfrak{g}}} S^\bullet(\mathfrak{g})=Z_{S^{\bullet}(\mathfrak{g}[t])}(S^\bullet(\mathfrak{g})^\mathfrak{g}).	
\end{equation*}
\eprop
\prf
Follows from~\cite[Lemma~6.9]{ir2}.
\epr

\prop{}\label{gr_to_univ}
The algebras $\on{gr}_{PBW}\CA^{\mathrm{u}}_\chi,\, \on{gr}_{PBW}\ol{\CA}{}^{\mathrm{u}}_{\chi} \subset S^\bullet(\mathfrak{g}[t])$ are both equal to the subalgebra of $S^{\bullet}(\mathfrak{g}[t])^{\mathfrak{z}_{\mathfrak{g}}(\chi)}$ generated by $\ol{\CA}{}^{\mathrm{u}}_{0} \subset S^{\bullet}(\mathfrak{g}[t])^{\mathfrak{g}}$ and $\ol{\CA}_\chi \subset S^{\bullet}(\mathfrak{g})^{\mathfrak{z}_{\mathfrak{g}}(\chi)}$. Moreover, this algebra coincides with the tensor product $\ol{\CA}{}^{\mathrm{u}}_{0} \otimes_{S^\bullet(\mathfrak{g})^{\mathfrak{g}}} \ol{\CA}_\chi$ and is a maximal Poisson-commutative subalgebra of $S^\bullet(\mathfrak{g}[t])^{\mathfrak{z}_{\mathfrak{g}}(\chi)}$. 
\eprop
\prf
Let us show that $S^\bullet(\mathfrak{g})^{\mathfrak{g}} \subset  \on{gr}_{PBW}\CA^{\mathrm{u}}_\chi \cap \on{gr}_{PBW}\ol{\CA}{}^{\mathrm{u}}_{\chi}$. Note that for every $X=x_1\ldots x_l \in S^\bullet(\mathfrak{g})$ we can write 
\begin{equation*}
X^{(0)}=x_1[0]\ldots x_l[0]+\text{other terms}.  
\end{equation*}
The image of $X^{(0)}$ in the Hamiltonian reduction (\ref{ham_red_symm_classic!}) inside $S^\bullet(\mathfrak{g}[t])$ is equal to
\begin{equation*}
X+\text{lower degree terms}
\end{equation*}
so the image of $X^{(0)}$ in $\on{gr}_{PBW}S^\bullet(\mathfrak{g}[t])$ is $X$. We conclude that $S^\bullet(\mathfrak{g})^{\mathfrak{g}} \subset \on{gr}_{PBW}\ol{\CA}{}^{\mathrm{u}}_{\chi}$. In the same way we show that $S^\bullet(\mathfrak{g})^{\mathfrak{g}} \subset \on{gr}_{PBW}\CA^{\mathrm{u}}_\chi$.

Let us now show that the algebra generated by $\ol{\CA}{}^{\mathrm{u}}_{0}$ and $\ol{\CA}_\chi$ is contained in $\on{gr}_{PBW}\CA^{\mathrm{u}}_\chi \cap \on{gr}_{PBW}\ol{\CA}{}^{\mathrm{u}}_{\chi}$. To see that $\ol{\CA}{}^{\mathrm{u}}_{0}$ lies in both $\on{gr}_{PBW}\CA^{\mathrm{u}}_{\chi},\, \on{gr}_{PBW}\ol{\CA}{}^{\mathrm{u}}_{\chi}$ recall (see Proposition~\ref{generators_of_class_univ_gaud}) that $\ol{\CA}{}^{\mathrm{u}}_{0}$ is generated by the elements $D^r\Phi,\, \Phi \in S^\bullet(\mathfrak{g})^{\mathfrak{g}},\, r\in \BZ_{\geqslant 0}$ and we already know that $S^\bullet(\mathfrak{g})^{\mathfrak{g}} \subset \on{gr}_{PBW}\CA^{\mathrm{u}}_{\chi} \cap \on{gr}_{PBW}\ol{\CA}{}^{\mathrm{u}}_{\chi}$.  Note now that if $X=x_1\ldots x_l \in S^\bullet(\mathfrak{g})$ then we can write 
\begin{equation*}
X^{(r)}=\sum_{j_1,\ldots,j_l\geqslant 0,\, j_1+\ldots+j_l=r}x_1[j_1]\ldots x_l[j_l]+\text{other terms}.    
\end{equation*}
The image of $X^{(r)}$ in $S^\bullet(\mathfrak{g}[t])$ is equal to 
\begin{equation*}
\frac{1}{r!}D^rX+\text{lower degree terms}.
\end{equation*}
We conclude that the image of $r!X$ in $\on{gr}_{PBW}S^\bullet(\mathfrak{g}[t])$ is exactly $D^rX$. It follows that $\ol{\CA}{}^{\mathrm{u}}_0 \subset \on{gr}_{PBW}\ol{\CA}{}^{\mathrm{u}}_\chi$. In the same way we see that $\ol{\CA}{}^{\mathrm{u}}_0 \subset \on{gr}_{PBW}\CA_\chi^{\mathrm{u}}$.

To see that $\ol{\CA}_\chi$ lies in $\on{gr}_{PBW}\CA^{\mathrm{u}}_{\chi} \cap \on{gr}_{PBW}\ol{\CA}{}^{\mathrm{u}}_\chi$ recall that $\ol{\CA}_\chi$ is generated by $\partial_\chi^r \Phi,\, \Phi \in S^{\bullet}(\mathfrak{g})^{\mathfrak{g}}$. Note that if $X:=x_1\ldots x_k \in S^{\bullet}(\mathfrak{g})$ then we can write
\begin{multline*}
X^{(-r)}=\sum_{1 \leqslant i_1<\ldots<i_r \leqslant l}x_1[0]\ldots x_{i_1-1}[0]x_{i_1}[-1]x_{i_1+1}[0]\ldots x_{i_r-1}[0]x_{i_r}[-1]x_{i_r+1}[0]\ldots x_{l}[0] +\\
+\text{other terms}.
\end{multline*}
The image of $X^{(-r)}$ in $S^{\bullet}(\mathfrak{g}[t])$ is equal to 
\begin{equation*}
\frac{1}{r!}\partial^r_\chi X+\text{lower degree terms}. \end{equation*}
We conclude that the class of the image of $r! X^{(-r)}$ in  $\on{gr}_{PBW}S^{\bullet}(\mathfrak{g}[t])$ is exactly $\partial^r_\chi X$. This observation finishes the proof of the fact that $\ol{\CA}_\chi \subset \on{gr}_{PBW}\ol{\CA}{}^{\mathrm{u}}_{\chi}$. In the same way we show that $\ol{\CA}_\chi \subset  \on{gr}_{PBW}\CA^{\mathrm{u}}_\chi$

Let us now describe the subalgebra of $S^{\bullet}(\mathfrak{g}[t])^{\mathfrak{z}_{\mathfrak{g}}(\chi)}$ generated by $\ol{\CA}{}^{\mathrm{u}}_0,\, \ol{\CA}_\chi$. 
Note that $\ol{\CA}{}^{\mathrm{u}}_0 \subset S^\bullet(\mathfrak{g}[t])^{\mathfrak{g}},\, \ol{\CA}_\chi \subset S^\bullet(\mathfrak{g})$ and by Proposition~\ref{prod_inv_funct} we have 
\begin{equation*}
S^{\bullet}(\mathfrak{g}[t])^{\mathfrak{g}}\cdot S^\bullet(\mathfrak{g})=S^\bullet(\mathfrak{g}[t])^{\mathfrak{g}} \otimes_{S^\bullet(\mathfrak{g})^{\mathfrak{g}}} S^\bullet(\mathfrak{g})=Z_{S^\bullet(\mathfrak{g}[t])}(S^\bullet(\mathfrak{g})^{\mathfrak{g}}).
\end{equation*}
We conclude that 
\begin{multline*}
\ol{\CA}{}^{\mathrm{u}}_0 \cdot \ol{\CA}_\chi=\ol{\CA}{}^{\mathrm{u}}_0 \otimes_{S^\bullet(\mathfrak{g})^{\mathfrak{g}}} \ol{\CA}_\chi
\subset S^\bullet(\mathfrak{g}[t])^{\mathfrak{g}}\otimes_{S^\bullet(\mathfrak{g})^{\mathfrak{g}}} S^\bullet(\mathfrak{g})^{\mathfrak{z}_{\mathfrak{g}}(\chi)}=\\
=S^\bullet(\mathfrak{g}[t])^{\mathfrak{g}} \cdot S^\bullet(\mathfrak{g})^{\mathfrak{z}_{\mathfrak{g}}(\chi)}=Z_{S^\bullet(\mathfrak{g}[t])^{\mathfrak{z}_{\mathfrak{g}}(\chi)}}(S^\bullet(\mathfrak{g})^{\mathfrak{g}})
\subset 
S^\bullet(\mathfrak{g}[t])^{\mathfrak{z}_{\mathfrak{g}}(\chi)}.
\end{multline*}
It then follows from \cite[Corollary~4.10]{ir2} and Proposition~\ref{A_chi_max} that $\ol{\CA}{}^{\mathrm{u}}_0 \cdot \ol{\CA}_\chi$ is the maximal Poisson-commutative subalgebra of $S^\bullet(\mathfrak{g}[t])^{\mathfrak{g}} \cdot S^\bullet(\mathfrak{g})^{\mathfrak{z}_{\mathfrak{g}}(\chi)}$.

Note also that if $x \in S^\bullet(\mathfrak{g}[t])^{\mathfrak{z}_{\mathfrak{g}}(\chi)}$ commutes with $\ol{\CA}{}^{\mathrm{u}}_0 \cdot \ol{\CA}_\chi$ then (since $S^\bullet(\mathfrak{g})^{\mathfrak{g}} \subset \ol{\CA}{}^{\mathrm{u}}_0 \cap \ol{\CA}_\chi$) we must have $x \in \mathfrak{z}_{S^\bullet(\mathfrak{g}[t])^{\mathfrak{z}_{\mathfrak{g}}(\chi)}}(S^\bullet(\mathfrak{g})^{\mathfrak{g}})=S^\bullet(\mathfrak{g}[t])^{\mathfrak{g}}\otimes_{S^\bullet(\mathfrak{g})^{\mathfrak{g}}} S^\bullet(\mathfrak{g})^{\mathfrak{z}_{\mathfrak{g}}(\chi)}$. Recall that $\ol{\CA}{}^{\mathrm{u}}_0 \cdot \ol{\CA}_\chi \subset S^\bullet(\mathfrak{g}[t])^{\mathfrak{g}} \cdot S^\bullet(\mathfrak{g})^{\mathfrak{z}_{\mathfrak{g}}(\chi)}$ is maximal Poisson-commutative so we conclude that $x \in \ol{\CA}{}^{\mathrm{u}}_0 \cdot \ol{\CA}_\chi$. We have shown that $\ol{\CA}^{\mathrm{u}}_0 \cdot \ol{\CA}_\chi \subset S^\bullet(\mathfrak{g}[t])^{\mathfrak{z}_{\mathfrak{g}}(\chi)}$ is maximal Poisson-commutative. 

Recall now that $\ol{\CA}{}^{\mathrm{u}}_0 \cdot \ol{\CA}_\chi$ is contained in both $\on{gr}_{PBW}\ol{\CA}{}^{\mathrm{u}}_\chi,\, \on{gr}_{PBW}\CA_\chi^{\mathrm{u}}$ so from the maximality of  $\ol{\CA}{}^{\mathrm{u}}_0 \cdot \ol{\CA}_\chi$ and commutativity of $\on{gr}_{PBW}\ol{\CA}{}^{\mathrm{u}}_\chi,\, \on{gr}_{PBW}\CA_\chi^{\mathrm{u}}$ we conclude that 
\begin{equation*}
\on{gr}_{PBW}\ol{\CA}{}^{\mathrm{u}}_\chi=\ol{\CA}{}^{\mathrm{u}}_0 \cdot \ol{\CA}_\chi=\on{gr}_{PBW}\CA_\chi^{\mathrm{u}}.
\end{equation*}
\epr

\cor{}
Subalgebras 
\begin{equation*}
\CA^{\mathrm{u}}_\chi \subset U(\mathfrak{g}[t])^{\mathfrak{z}_{\mathfrak{g}}(\chi)},\, \ol{\CA}{}^{\mathrm{u}}_\chi \subset S^{\bullet}(\mathfrak{g}[t])^{\mathfrak{z}_{\mathfrak{g}}(\chi)}
\end{equation*}
are maximal (resp., maximal Poisson) commutative. In particular, $\ol{\CA}{}^{\mathrm{u}}_\chi \subset S^{\bullet}(\mathfrak{g}[t])^{\mathfrak{z}_{\mathfrak{g}}(\chi)}$ is algebraically closed.
\ecor
\prf
Follows from Proposition~\ref{gr_to_univ} which claims that the associated graded $\on{gr}_{PBW}\CA^{\mathrm{u}}_\chi=\on{gr}_{PBW}\ol{\CA}{}^{\mathrm{u}}_{\chi}$ is maximal.
\epr

Let us now recall the dimensions of the algebras $\ol{\CA}{}^{\mathrm{u}}_0,\, \ol{\CA}_\chi$. Let $d_i,\, i=1,\ldots,\on{rk}\mathfrak{g}$ be the degree of $\Phi_i$ and let $p_i,\, i=1,\ldots, \on{rk}\mathfrak{z}_{\mathfrak{g}}(\chi)^{\mathrm{der}}$ be the degrees of the generators of the algebra $S^\bullet(\mathfrak{z}_{\mathfrak{g}}(\chi)^{\mathrm{der}})^{\mathfrak{z}_{\mathfrak{g}}(\chi)^{\mathrm{der}}}$. 
\prop{}\label{poinc_gaud_U}
We have \begin{equation*}
\on{dim}_{F_1}\ol{\CA}{}^{\mathrm{u}}_0=\prod_{i=1}^{\on{rk}\mathfrak{g}}\prod_{l=d_i}^{\infty}\frac{1}{1-q^l},~ \on{dim}_{PBW}\ol{\CA}_\chi
=\frac{\prod_{i=1}^{\on{rk}\mathfrak{z}_{\mathfrak{g}}(\chi)^{\mathrm{der}}}\prod_{l=1}^{p_i-1}(1-q^l)}{\prod_{i=1}^{\on{rk}\mathfrak{g}}\prod_{l=1}^{d_i}(1-q^l)}.
\end{equation*}
\eprop
\prf
The claim about the dimension of the algebra $\ol{\CA}{}^{\mathrm{u}}_0$ follows from the fact that $\ol{\CA}{}^{\mathrm{u}}_0$ is freely generated by $D^r\Phi_i,\, r \geqslant 0$ (see Proposition~\ref{generators_of_class_univ_gaud}).

Let us now compute the dimension of $\ol{\CA}_\chi$.
Recall that by part $(b)$ of Proposition~\ref{ext_A_chi} we have $\ol{\CA}_{(\chi,\chi_1)}=\underset{\epsilon \rightarrow 0}{\on{lim}}\, \ol{\CA}_{\chi+\epsilon\chi_1}$.
Since $\chi+\epsilon\chi_1 \in \mathfrak{g}$ is regular and $\ol{\CA}_{\chi+\epsilon\chi_1} \subset S^\bullet(\mathfrak{g})$ are graded with respect to the grading by the degree (that induces the filtration that we use to define the limit) we conclude from Lemma~\ref{lim_properties} that
\begin{equation*}
\on{dim}_{PBW}\ol{\CA}_{(\chi,\chi_1)}=\prod_{i=1}^{\on{rk}\mathfrak{g}}\prod_{l=1}^{d_i}\frac{1}{1-q^l}.
\end{equation*}

It now follows from part $(a)$ of Proposition~\ref{ext_A_chi} that 
\begin{equation*}
\on{dim}_{PBW}\ol{\CA}_\chi \cdot  \prod_{i=1}^{\on{rk}\mathfrak{z}_{\mathfrak{g}}(\chi)^{\mathrm{der}}}\prod_{l=1}^{p_i-1}\frac{1}{1-q^l}  = \on{dim}_{PBW}\ol{\CA}_{(\chi,\chi_1)}=\prod_{i=1}^{\on{rk}\mathfrak{g}}\prod_{l=1}^{d_i}\frac{1}{1-q^l}. \end{equation*}
We conclude that \begin{equation*}
\on{dim}_{PBW}\ol{\CA}_\chi=\frac{\prod_{i=1}^{\on{rk}\mathfrak{g}}\prod_{l=1}^{d_i}\frac{1}{1-q^l}}{\prod_{i=1}^{\on{rk}\mathfrak{z}_{\mathfrak{g}}(\chi)^{\mathrm{der}}}\prod_{l=1}^{p_i-1}\frac{1}{1-q^l}}.
\end{equation*}
\epr

We are now ready to compute the dimension of the algebras $\ol{\CA}{}^{\mathrm{u}}_\chi,\, \CA^{\mathrm{u}}_\chi$.
\prop{}\label{poinc_univ_shift}
We have 
\begin{equation*}
\on{dim}_{F_1}\ol{\CA}{}^{\mathrm{u}}_\chi=\on{dim}_{F_1}\CA^{\mathrm{u}}_\chi=
\prod_{i=1}^{\on{rk}\mathfrak{z}_{\mathfrak{g}}(\chi)^{\mathrm{der}}}
\prod_{l=p_i}^{\infty}\frac{1}{1-q^l}
\cdot \prod_{l=1}^\infty \frac{1}{(1-q^l)^{\on{rk}\mathfrak{z}_{\mathfrak{g}}(\chi)-\on{rk}\mathfrak{z}_{\mathfrak{g}}(\chi)^{\mathrm{der}}}}
.
\end{equation*}
\eprop
\prf
Since 
\begin{equation*}
\on{gr}_{PBW}\ol{\CA}{}^{\mathrm{u}}_\chi=\on{gr}_{PBW}\CA^{\mathrm{u}}_\chi=\ol{\CA}{}^{\mathrm{u}}_{0} \otimes_{S^\bullet(\mathfrak{g})^{\mathfrak{g}}} \ol{\CA}_\chi
\end{equation*}
and 
\begin{equation*}
\on{dim}_{F_1}\ol{\CA}{}^{\mathrm{u}}_{0}=\prod_{i=1}^{\on{rk}\mathfrak{g}}\prod_{l=d_i}^{\infty}\frac{1}{1-q^l},\, \on{dim}_{PBW}\ol{\CA}_\chi=\frac{\prod_{i=1}^{\on{rk}\mathfrak{g}}\prod_{l=1}^{d_i}\frac{1}{1-q^l}}{\prod_{i=1}^{\on{rk}\mathfrak{z}_{\mathfrak{g}}(\chi)^{\mathrm{der}}}\prod_{l=1}^{p_i-1}\frac{1}{1-q^l}},\, \on{dim}_{PBW}S^\bullet(\mathfrak{g})^{\mathfrak{g}}=\prod_{i=1}^{\on{rk}\mathfrak{g}}\frac{1}{1-q^{d_i}}
\end{equation*}
it follows that 
\begin{multline*}
\on{dim}_{F_1}\ol{\CA}{}^{\mathrm{u}}_\chi=\frac{\prod_{i=1}^{\on{rk}\mathfrak{g}}\prod_{l=d_i+1}^{\infty}\frac{1}{1-q^l}\cdot \prod_{i=1}^{\on{rk}\mathfrak{g}}\prod_{l=1}^{d_i}\frac{1}{1-q^l}}{\prod_{i=1}^{\on{rk}\mathfrak{z}_{\mathfrak{g}}(\chi)^{\mathrm{der}}}\prod_{l=1}^{p_i-1}\frac{1}{1-q^l}}=\\
=\frac{\prod_{l=1}^{\infty}\frac{1}{(1-q^l)^{\on{rk}\mathfrak{g}}}}{\prod_{i=1}^{\on{rk}\mathfrak{z}_{\mathfrak{g}}(\chi)^{\mathrm{der}}}\prod_{l=1}^{p_i-1}\frac{1}{1-q^l}}=
\prod_{i=1}^{\on{rk}\mathfrak{z}_{\mathfrak{g}}(\chi)^{\mathrm{der}}}
\prod_{l=p_i}^{\infty}\frac{1}{1-q^l}
\cdot \prod_{l=1}^\infty \frac{1}{(1-q^l)^{\on{rk}\mathfrak{z}_{\mathfrak{g}}(\chi)-\on{rk}\mathfrak{z}_{\mathfrak{g}}(\chi)^{\mathrm{der}}}}
.
\end{multline*}
\epr

\cor{}\label{a_chi_via_a_centr}
We have 
$
\on{dim}_{F_1}\ol{\CA}{}^{\mathrm{u}}_\chi =\on{dim}_{F_1}\ol{\CA}^{\mathrm{u}}_0(\mathfrak{z}_{\mathfrak{g}}(\chi)).   
$
\ecor
\prf
We have 
\begin{multline*}
\on{dim}_{F_1}\ol{\CA}{}^{\mathrm{u}}_\chi=\frac{\prod_{l=1}^{\infty}\frac{1}{(1-q^l)^{\on{rk}\mathfrak{g}}}}{\prod_{i=1}^{\on{rk}\mathfrak{z}_{\mathfrak{g}}(\chi)^{\mathrm{der}}}\prod_{l=1}^{p_i-1}\frac{1}{1-q^l}}=\\
=
\prod_{i=1}^{\on{rk}\mathfrak{z}_{\mathfrak{g}}(\chi)^{\mathrm{der}}}
\prod_{l=p_i}^{\infty}\frac{1}{1-q^l}
\cdot \prod_{l=1}^\infty \frac{1}{(1-q^l)^{\on{rk}\mathfrak{z}_{\mathfrak{g}}(\chi)-\on{rk}\mathfrak{z}_{\mathfrak{g}}(\chi)^{\mathrm{der}}}}
=\on{dim}_{F_1}\ol{\CA}_0^{\mathrm{u}}(\mathfrak{z}_{\mathfrak{g}}(\chi)).
\end{multline*}
\epr

\section{Universal inhomogeneous Gaudin subalgebras as  centralizers
}\label{univ_inhom_as_central_section}
It follows from~\cite{ir2} that for $\chi=0$ we have 
\begin{equation*}
\ol{\CA}{}^{\mathrm{u}}_0=Z_{S^\bullet(\mathfrak{g}[t])^{\mathfrak{g}}}(\Omega_0),\, \CA^{\mathrm{u}}_0=Z_{U(\mathfrak{g}[t])^{\mathfrak{g}}}(\tilde{\Omega}_0).    
\end{equation*}
The main goal of this section is to prove that the same equalities hold for $\ol{\CA}{}^{\mathrm{u}}_\chi,\, \CA^{\mathrm{u}}_{\chi}$ i.e. that 
\begin{equation*}
\ol{\CA}{}^{\mathrm{u}}_\chi=Z_{S^\bullet(\mathfrak{g}[t])^{\mathfrak{z}_{\mathfrak{g}}(\chi)}}(\Omega_\chi),\, \CA^{\mathrm{u}}_\chi=Z_{U(\mathfrak{g}[t])^{\mathfrak{z}_{\mathfrak{g}}(\chi)}}(\tilde{\Omega}_\chi). \end{equation*}


\lem{}\label{cent_chi}
For every $m>0$ we have $Z_{S^\bullet(\mathfrak{g}[t])}(\chi[m])=S^\bullet(\mathfrak{z}_{\mathfrak{g}}(\chi)[t])$.
\elem
\prf
Same proof as of~\cite[Lemma~4]{r}.
\epr

\prop{}\label{class_gaud_centr}
We have $\ol{\CA}{}^{\mathrm{u}}_\chi=Z_{S^{\bullet}(\mathfrak{g}[t])^{\mathfrak{z}_{\mathfrak{g}}(\chi)}}(\Omega_\chi)$.
\eprop
\prf
Consider the family $Z_{S^\bullet(\mathfrak{g}[t])^{\mathfrak{z}_{\mathfrak{g}}(\chi)}}(\Omega_{\kappa\chi}),\, \kappa \in \BC^\times$.  Note that $$Z_{S^\bullet(\mathfrak{g}[t])^{\mathfrak{z}_{\mathfrak{g}}(\chi)}}(\Omega_{\kappa\chi}) \simeq Z_{S^\bullet(\mathfrak{g}[t])^{\mathfrak{z}_{\mathfrak{g}}(\chi)}}(\Omega_\chi)$$ via the Poisson automorphism of $S^\bullet(\mathfrak{g}[t])$ given by the map $x[n] \mapsto \kappa^{n}x[n]$.  It follows that we can consider the limit 
\begin{equation*}
\underset{\kappa \ra \infty}{\on{lim}}\,Z_{S^\bullet(\mathfrak{g}[t])^{\mathfrak{z}_{\mathfrak{g}}(\chi)}}(\Omega_{\kappa\chi}).
\end{equation*}
Pick an element 
\begin{equation*}
a_0 \in \underset{\kappa \ra \infty}{\on{lim}}\,Z_{S^\bullet(\mathfrak{g}[t])^{\mathfrak{z}_{\mathfrak{g}}(\chi)}}(\Omega_{\kappa\chi}).
\end{equation*}
Note that $a_0 \in F^i_1S^{\bullet}(\mathfrak{g}[t])$ for some $i \geqslant 0$. Recall that we have a map $$f\colon \BA^1 \ra \on{Gr}(d(i),F^i_1 S^{\bullet}(\mathfrak{g}[t]))$$ that sends $\epsilon \in \BA^1$ to $F^i_1 S^{\bullet}(\mathfrak{g}[t])) \cap Z_{S^\bullet(\mathfrak{g}[t])^{\mathfrak{z}_{\mathfrak{g}}(\chi)}}(\epsilon\Omega_0+\chi[2])$ and sends $0$ to $\underset{\kappa \ra \infty}{\on{lim}}\,(F^i_1 S^{\bullet}(\mathfrak{g}[t]) \cap Z_{S^\bullet(\mathfrak{g}[t])^{\mathfrak{z}_{\mathfrak{g}}(\chi)}}(\Omega_{\kappa\chi}))$. 
It follows from Lemma~\ref{lim_in_lim} that there exists a morphism
$a\colon \BA^1 \ra F^i_1 S^{\bullet}(\mathfrak{g}[t]))$ such that $a(\epsilon) \in Z_{S^\bullet(\mathfrak{g}[t])^{\mathfrak{z}_{\mathfrak{g}}(\chi)}}(\epsilon\Omega_0+\chi[2])$ for every $\epsilon \in \BA^1 \setminus \{0\}$ and $a(0)=a_0$. Then we can write
$a(\epsilon)=a_0+\epsilon a_1+\epsilon^2 a_2+\ldots$ with $a_i \in F^i_1S^\bullet(\mathfrak{g}[t])$. 

We have $\{\epsilon\Omega_0+\chi[2],a\}=0$ so \begin{equation*}
\{\chi[2],a_0\}=0,\, \{\Omega_0,a_0\}+\{\chi[2],a_1\}=0.
\end{equation*}
From $\{\chi[2],a_0\}=0$ we conclude by Lemma~\ref{cent_chi} that 
\begin{equation}\label{a_0_in_centr}
a_0 \in S^{\bullet}(\mathfrak{z}_{\mathfrak{g}}(\chi)[t]).
\end{equation}

Pick  a basis of $\mathfrak{g}$, consisting of root vectors and consider the corresponding decomposition $\mathfrak{g}=\mathfrak{z}_{\mathfrak{g}}(\chi) \oplus \mathfrak{m}$. It induces the decomposition $S^\bullet(\mathfrak{g}[t])=S^{\bullet}(\mathfrak{z}_{\mathfrak{g}}(\chi)[t]) \oplus S^\bullet(\mathfrak{g}[t])\mathfrak{m}[t]$. Let $\pi\colon S^\bullet(\mathfrak{g}[t]) \twoheadrightarrow S^{\bullet}(\mathfrak{z}_{\mathfrak{g}}(\chi)[t])$ be the projection. Note that $\{\chi[2],S^\bullet(\mathfrak{z}_{\mathfrak{g}}(\chi)[t])\}=0$, $\{\chi[2],S^\bullet(\mathfrak{g}[t])\mathfrak{m}[t]\} \subset S^\bullet(\mathfrak{g}[t])\mathfrak{m}[t]$ so it is clear that $\pi(\{\chi[2],a_1\})=0$. It is also clear (use (\ref{a_0_in_centr})) that  we have $\pi(\{\Omega_0,a_0\})=\{\Omega_{0}(\mathfrak{z}_\mathfrak{g}(\chi)),a_0\}$, where $\Omega_{0}(\mathfrak{z}_\mathfrak{g}(\chi))$ is the element $\Omega_0$ for $\mathfrak{z}_\mathfrak{g}(\chi)$. 

We conclude that $\{\Omega_{0}(\mathfrak{z}_\mathfrak{g}(\chi)),a_0\}=0$. We also know that $a_0 \in S^{\bullet}(\mathfrak{z}_{\mathfrak{g}}(\chi)[t])^{\mathfrak{z}_{\mathfrak{g}}(\chi)}$. It follows that   
\begin{equation*}
\underset{\kappa \ra \infty}{\on{lim}}\,Z_{S^\bullet(\mathfrak{g}[t])^{\mathfrak{z}_{\mathfrak{g}}(\chi)}}(\Omega_{\kappa\chi}) \subset Z_{S^\bullet(\mathfrak{z}_{\mathfrak{g}}(\chi)[t])^{\mathfrak{z}_{\mathfrak{g}}(\chi)}}(\Omega_0({\mathfrak{z}_{\mathfrak{g}}(\chi)}))=\ol{\CA}{}^{\mathrm{u}}_0(\mathfrak{z}_{\mathfrak{g}}(\chi)),
\end{equation*}
where the last equality follows from~\cite[Proposition~4.9]{ir2}.
It remains to note that $\ol{\CA}{}^{\mathrm{u}}_{\chi} \subset Z_{S^\bullet(\mathfrak{g}[t])^{\mathfrak{z}_{\mathfrak{g}}(\chi)}}(\Omega_{\chi})$, so
\begin{equation}\label{dim_chi_infty}
\on{dim}_{F_1}\ol{\CA}{}^{\mathrm{u}}_{\chi} \leqslant \on{dim}_{F_1}Z_{S^\bullet(\mathfrak{g}[t])^{\mathfrak{z}_{\mathfrak{g}}(\chi)}}(\Omega_{\chi}) \leqslant  \on{dim}_{F_1}\left(\underset{\kappa \ra \infty}{\on{lim}}\,Z_{S^\bullet(\mathfrak{g}[t])^{\mathfrak{z}_{\mathfrak{g}}(\chi)}}(\Omega_{\kappa\chi})\right) \leqslant \on{dim}_{F_1} \ol{\CA}{}^{\mathrm{u}}_0(\mathfrak{z}_{\mathfrak{g}}(\chi)). 
\end{equation}
Recall that  $\on{dim}_{F_1}\ol{\CA}{}^{\mathrm{u}}_{\chi}=\on{dim}_{F_1} \ol{\CA}{}^{\mathrm{u}}_0(\mathfrak{z}_{\mathfrak{g}}(\chi))$ by Corollary~\ref{a_chi_via_a_centr} so the inequalities in~(\ref{dim_chi_infty}) are actually equalities and we must have $\ol{\CA}{}^{\mathrm{u}}_{\chi} = Z_{S^\bullet(\mathfrak{g}[t])^{\mathfrak{z}_{\mathfrak{g}}(\chi)}}(\Omega_{\chi})$ as desired.

\epr

\rem{}
{\em{It follows from the proof of Proposition~\ref{class_gaud_centr} that we have 
\begin{equation*}
\underset{\kappa \ra \infty}{\on{lim}}\,\ol{\CA}{}^{\mathrm{u}}_{\kappa\chi}=\ol{\CA}{}^{\mathrm{u}}_0(\mathfrak{z}_{\mathfrak{g}}(\chi)).
\end{equation*}
}}
\erem

\rem{}\label{cent_poiss_w_chi_reg}
{\em{Note that for regular $\chi$ we also have $\ol{\CA}{}^{\mathrm{u}}_\chi=Z_{S^\bullet(\mathfrak{g}[t])}(\omega_\chi)$. Indeed, we have the filtration $F_2$ on $S^\bullet(\mathfrak{g}[t])$ with $\on{deg}_2x[m]=m$. Note that $\on{gr}_2(\omega_\chi)=\on{gr}_2(\chi[1]])$ so we conclude that $\on{gr}_2(Z_{S^\bullet(\mathfrak{g}[t])}(\omega_\chi)) \subset Z_{S^\bullet(\mathfrak{g}[t])}(\chi[1])=S^\bullet(\mathfrak{h}[t])$ (see Lemma \ref{cent_chi}). It follows that $\on{dim}_{F_1}Z_{S^\bullet(\mathfrak{g}[t])}(\omega_\chi) \leqslant \on{dim}_{F_1}S^\bullet(\mathfrak{h}[t])$. Recall also that $\ol{\CA}{}^{\mathrm{u}}_\chi \subset Z_{S^\bullet(\mathfrak{g}[t])}(\omega_\chi)$ and $\on{dim}_{F_1}\ol{\CA}{}^{\mathrm{u}}_\chi=\on{dim}_{F_1}S^\bullet(\mathfrak{h}[t])$. We conclude that $\ol{\CA}{}^{\mathrm{u}}_\chi = Z_{S^\bullet(\mathfrak{g}[t])}(\omega_\chi)$.}}
\erem


So we have shown that the algebra $\ol{\CA}{}^{\mathrm{u}}_\chi$ coincides with the centralizer $Z_{S^\bullet(\mathfrak{g}[t])^{\mathfrak{z}_{\mathfrak{g}}(\chi)}}(\Omega_\chi)$. Let us now quantize this statement i.e. prove that 
\begin{equation*}
\CA^{\mathrm{u}}_{\chi}=Z_{U(\mathfrak{g}[t])^{\mathfrak{z}_{\mathfrak{g}}(\chi)}}(\tilde{\Omega}_\chi).    
\end{equation*}

We first relate the algebras $\ol{\CA}{}^{\mathrm{u}}_\chi,\,\CA^{\mathrm{u}}_{\chi}$.
Recall that the $PBW$-filtration on $U(\mathfrak{g}[t])$ and  consider the Rees family $U_{\epsilon}(\mathfrak{g}[t]),\, \epsilon \in \BC$ as in Section~\ref{rees_constr}. Note that $U_{\epsilon}(\mathfrak{g}[t]) \simeq U(\mathfrak{g}[t])$ for $\epsilon \neq 0$ and $U_{0}(\mathfrak{g}[t])=S^\bullet(\mathfrak{g}[t])$. 
For $\epsilon \neq 0$ we pick the natural isomorphisms 
\begin{equation*}
 \varphi_\epsilon\colon U(\mathfrak{g}[t]) \iso U_{\epsilon}(\mathfrak{g}[t])
\end{equation*}
that  send an element $a$ of degree $r$ to $[\epsilon^{-r}\hbar^r a]$. Note that these isomorphisms are $\mathfrak{g}$-equivariant.

Recall that the algebra $U(\mathfrak{g}[t])$ is graded with the grading $\on{deg}_2 x[n]=n$. Note that this grading induces the grading on $Rees(U(\mathfrak{g}[t]))$ ($\on{deg}_2 \hbar x[n]=n$) such that $\BC[\hbar] \subset U_{\epsilon}(\mathfrak{g}[t])$ lies in the zero graded component. Thus, we obtain the grading on every $U_\epsilon(\mathfrak{g}[t])$. We can consider the automorphism $d_\epsilon\colon U_\epsilon(\mathfrak{g}[t]) \iso U_\epsilon(\mathfrak{g}[t])$ that sends $a$ such that $\on{deg}_2a=i$ to $\epsilon^{-i}a$. The automorphism $d_\epsilon$ is $\mathfrak{g}$-equivariant.  We obtain the composition 
\begin{equation*}
 d_\epsilon \circ \varphi_\epsilon\colon U(\mathfrak{g}[t])  \iso U_{\epsilon}(\mathfrak{g}[t]).
\end{equation*}

Set $\tilde{\varphi}_\epsilon:=d_\epsilon \circ \varphi_\epsilon$.
For $\epsilon \neq 0$ we can  embed $\CA_{\chi}^{\mathrm{u}} \subset U_{\epsilon}(\mathfrak{g}[t])$ via $\tilde{\varphi}_\epsilon$.

\prop{}\label{quant_sh_via_class}
After the embedding $\tilde{\varphi}_\epsilon\colon \CA^{\mathrm{u}}_\chi\subset U_\epsilon(\mathfrak{g}[t])$, we have  $\underset{\epsilon \ra 0}{\on{lim}}\,\CA_{\chi}^{\mathrm{u}}=\ol{\CA}{}^{\mathrm{u}}_{\chi}$.
\eprop
\prf
Let us check that $\Omega_\chi \in \underset{\epsilon \ra 0}{\on{lim}}\,\CA_{\chi}^{\mathrm{u}}$. Recall that $\tilde{\Omega}_\chi \in \CA^{\mathrm{u}}_\chi$ (see Remark~\ref{A_chi_cont_omega}). After the identification $\tilde{\varphi}_\epsilon$ this element becomes 
\begin{equation*}
\tilde{\Omega}_{\chi,\epsilon}=\sum_a \hbar^2 \epsilon^{-3}x_a x_a[1] + \epsilon^{-3}\hbar\chi[2]
\end{equation*}
and the limit of $\epsilon^3\tilde{\Omega}_{\chi,\epsilon}$ as $\epsilon \ra 0$ is exactly $\Omega_\chi$.

We conclude (use Lemma \ref{lim_descr_via_elements_rees}) that $\underset{\epsilon \ra 0}{\on{lim}}\,\CA_{\chi}^{\mathrm{u}} \subset Z_{S^\bullet(\mathfrak{g}[t])^{\mathfrak{z}_\mathfrak{g}(\chi)}}(\Omega_\chi)=\ol{\CA}{}^{\mathrm{u}}_{\chi}$, where the last equality holds by Proposition~\ref{class_gaud_centr}. Since the dimension of $\underset{\epsilon \ra 0}{\on{lim}}\,\CA_{\chi}^{\mathrm{u}}$ is at least  $\on{dim}_{F_1}\CA^{\mathrm{u}}_\chi$ (Lemma \ref{dim_lim_Rees}) and by Proposition~\ref{gr_to_univ} we have $\on{dim}_{F_1}\CA^{\mathrm{u}}_\chi=\on{dim}_{F_1}\ol{\CA}{}^{\mathrm{u}}_\chi$ we conclude  that $\underset{\epsilon \ra 0}{\on{lim}}\,\CA_{\chi}^{\mathrm{u}}=\ol{\CA}{}^{\mathrm{u}}_\chi$.
\epr

\prop{}\label{ca_as_centr}
We have $\CA^{\mathrm{u}}_\chi=Z_{U(\mathfrak{g}[t])^{\mathfrak{z}_{\mathfrak{g}}(\chi)}}(\tilde{\Omega}_\chi)$.
\eprop
\prf
Since the element $\epsilon^3\tilde{\Omega}_{\chi,\epsilon}=\epsilon^3\tilde{\varphi}_\epsilon(\tilde{\Omega}_\chi)$ goes to $\Omega_\chi$ as $\epsilon$ goes to zero we conclude that 
\begin{equation*}
\underset{\epsilon \ra 0}{\on{lim}}\,Z_{U_{\epsilon}(\mathfrak{g}[t])^{\mathfrak{z}_{\mathfrak{g}}(\chi)}}(\tilde{\Omega}_{\chi,\epsilon}) \subset Z_{S^\bullet(\mathfrak{g}[t])^{\mathfrak{z}_{\mathfrak{g}}(\chi)}}(\Omega_\chi)=\ol{\CA}{}^{\mathrm{u}}_\chi.
\end{equation*}
Note now that since $\CA^{\mathrm{u}}_{\chi} \subset Z_{U(\mathfrak{g}[t])^{\mathfrak{z}_{\mathfrak{g}}(\chi)}}(\tilde{\Omega}_\chi)$ and passing to limit may only increase the dimension then
\begin{equation}\label{dim_est}
\on{dim}_{F_1}\CA^{\mathrm{u}}_\chi \leqslant \on{dim}_{F_1}Z_{U(\mathfrak{g}[t])^{\mathfrak{z}_{\mathfrak{g}}(\chi)}}(\tilde{\Omega}_\chi) \leqslant \on{dim}_{F_1}\underset{\epsilon \ra 0}{\on{lim}}\,Z_{U_\epsilon(\mathfrak{g}[t])^{\mathfrak{z}_{\mathfrak{g}}(\chi)}}(\tilde{\Omega}_{\chi,\epsilon}) \leqslant \on{dim}_{F_1}\ol{\CA}{}^{\mathrm{u}}_\chi. 
\end{equation}
Since by Proposition~\ref{gr_to_univ} we have $\on{dim}_{F_1}\CA^{\mathrm{u}}_\chi=\on{dim}_{F_1}\ol{\CA}{}^{\mathrm{u}}_\chi$ we conclude that inequalities in~(\ref{dim_est}) are equalities so  
$\CA^{\mathrm{u}}_\chi = Z_{U(\mathfrak{g}[t])^{\mathfrak{z}_{\mathfrak{g}}(\chi)}}(\tilde{\Omega}_\chi)$.
\epr

\rem{}\label{cent_til_w_chi_reg_chi}
{\em{Note that for regular $\chi$ we also have $\CA^{\mathrm{u}}_\chi=Z_{U(\mathfrak{g}[t])}(\tilde{\omega}_\chi)$. The proof is the same as the one of Proposition~\ref{ca_as_centr}, the only difference is that we use Remark~\ref{cent_poiss_w_chi_reg} instead of Proposition~\ref{class_gaud_centr}.
}}
\erem

\section{Yangian}\label{yangian}
Recall that $\mathfrak{g}$ is a simple Lie algebra.

\ssec{}{Definition of the Yangian}
\begin{defeni}{}
Yangian $Y(\mathfrak{g})$ is the associative $\BC$-algebra generated by 
\begin{equation*}
\{x,\, J(x)\,|\, x \in \mathfrak{g}\}
\end{equation*}
subject to the following relations:
\begin{equation*}
xy-yx=[x,y],\, J([x,y])=[J(x),y],\, J(cx+dy)=cJ(x)+dJ(y), \end{equation*}
\begin{equation*}
[J(x),[J(y),z]]=[x,[J(y),J(z)]]=\sum_{a_1, a_2, a_3}([x,x_{a_1}],[[y,x_{a_2}],[z,x_{a_3}]])\{x_{a_1},x_{a_2},x_{a_3}\},
\end{equation*}
\begin{multline*}
[[J(x),J(y)],[z,J(w)]]+[[J(z),J(w)],[x,J(y)]]=\\
=\sum_{a_1, a_2, a_3}(([x,x_{a_1}],[[y,x_{a_2}],[[z,w],x_{a_3}]])+([z,x_{a_1}],[[w,x_{a_2}],[[x,y]x_{a_3}]]))\{x_{a_1},x_{a_2},J(x_{a_3})\}
\end{multline*}
for all $x,y,z,w \in \mathfrak{g},\, c,d \in \BC$. Here $\{x_a\}|_{a=1,\ldots,\on{dim}\mathfrak{g}}$ is an orthonormal basis and 
\begin{equation*}
\{x_1,x_2,x_3\}:=\frac{1}{24}\sum_{\sigma \in S_3}x_{\sigma(1)}x_{\sigma(2)}x_{\sigma(3)}
\end{equation*}
for all $x_1, x_2, x_3 \in Y(\mathfrak{g})$.
\end{defeni}

By~\cite[Theorem~2]{d} Yangian $Y(\mathfrak{g})$ is a Hopf algebra. It follows from~\cite[Theorem~3]{d} that there is a unique formal series 
\begin{equation*}
\CR(u) = 1 + \sum_{k=1}^\infty \CR_k u^{-k} \in (Y(\mathfrak{g}) \otimes Y(\mathfrak{g}))[[u^{-1}]]
\end{equation*}
satisfying
\begin{equation*}
(\on{id} \otimes \Delta)\CR(u)=\CR_{12}(u)\CR_{13}(u),    
\end{equation*}
\begin{equation*}
\tau_{0,u}\Delta^{\mathrm{opp}}(x)=\CR(u)^{-1}(\tau_{0,u}\Delta(x))\CR(u)~\text{for all}~x \in Y(\mathfrak{g}), \end{equation*}
where $\tau_{0,u}$ is the automorphism of $Y(\mathfrak{g})[[u]]$ given by $x \mapsto x,\, J(x) \mapsto J(x)+ux$ for all $x \in \mathfrak{g}$.

We have
\begin{equation*}
\CR(u)=1+2\omega' u^{-1}+\Big(\sum_{a} (J(x_a)\otimes x_a -x_a \otimes J(x_a))+2(\omega')^2\Big) u^{-2}+O(u^{-3}),
\end{equation*}
where $\omega'=\frac{1}{2}\sum_a x_a \otimes x_a$.
See \cite{d} and~\cite[Section~3]{w} for  details. We will call $\CR(u)$ the universal $R$-matrix.

\subsection{$RTT$ realization of Yangian}\label{RTT_realiz_Yangian} 
Let $\rho\colon Y(\mathfrak{g}) \ra \on{End}(V)$ be a finite-dimensional representation of $Y(\mathfrak{g})$ that is not a direct sum of trivial representations w.r.t. $\mathfrak{g} \subset Y(\mathfrak{g})$.
We fix a basis $e_1,\ldots,e_{\on{dim}V}$ of $V$.
Let 
$
R(u-v):=(\rho \otimes \rho)(\CR(u-v))
$
be the image of the universal \(R\)-matrix in \(\End(V)^{\otimes 2}\). Using this data, we define the $RTT$-realization \(Y_V(\fg)\) of $Y(\mathfrak{g})$ as follows.



\begin{defeni}{}
The Yangian \(\yvg\) is a unital associative algebra generated by the elements \(t_{ij}^{(r)},1\leqslant i,j\leqslant\dim V;r\geqslant 1\) with the defining relations
\[
R(u-v)T_1(u)T_2(v)=T_2(v)T_1(u)R(u-v)
\in
\End(V)^{\otimes 2}\otimes\yvg[[u^{-1},v^{-1}]],
\]
\[
\eta^2(T(u))=T(u+\frac{1}{2}c_{\fg}),
\]
where \(\eta(T(u))=T(u)^{-1}\) is the antipode map and \(c_{\fg}\) is the value of the Casimir element $\sum_a x_a^2$ of \(\fg\) on the adjoint representation.
Here
\[
T(u)=\left [t_{ij}(u) \right ]_{i,j=1,\ldots,\dim V}\in\End V\otimes\yvg,
\]
\[
t_{ij}(u)=\delta_{ij}+\sum_{r \geqslant 1} t_{ij}^{(r)}u^{-r}
\]
and \(T_1(u)\) (resp. \(T_2(u)\)) is the image of \(T(u)\) in the first (resp. second) copy of \(\End V\). 
\end{defeni}

\th{}(\cite[Theorem~6.2]{w})
The assignment $T(u) \mapsto (\rho \otimes 1)\CR(-u)$ extends to an isomorphism of algebras $\Phi\colon Y_V(\mathfrak{g}) \iso Y(\mathfrak{g})$.
\eth

Recall that we have a  basis $e_1,\ldots,e_{\on{dim}V}$ of $V$. We  denote by  $E_{ij}|_{1 \leqslant i,j \leqslant \on{dim}V} \in \on{End}(V)$ the corresponding matrix units. 
We define $\CF_{ij} \in \on{End}(V)$ by the following identity
\begin{equation*}
\sum_{ij} E_{ij} \otimes \CF_{ij}=-(\rho \otimes 1)\Big(\sum_a x_a \otimes x_a\Big). \end{equation*}
\prop{} (\cite[Corollary~6.3]{w})\label{iso_Y_Y}
The map $\Phi^{-1}$ sends generators $\{\CF_{ij},\, J(\CF_{ij})\}$ of $Y(\mathfrak{g})$ to
\begin{equation*}
\CF_{ij} \mapsto t^{(1)}_{ij},\, J(\CF_{ij}) \mapsto t^{(2)}_{ij}-\frac{1}{2}\sum_{a=1}^{\on{dim}V}t_{ia}^{(1)}t_{aj}^{(1)}+\sum_{p,l=1}^{\on{dim}V}b_{pl}^{(ij)}t^{(1)}_{pl},    
\end{equation*}
where $b^{ij}_{pl}$ are certain scalars.
\eprop

\rem{}
{\em{Note that for $\mathfrak{g}=\mathfrak{sl}_n$ and $V=\BC^n$ (the standard representation of $\mathfrak{g}$), we have $\mathcal{F}_{ij}=-E_{ji}$.}}
\erem

\subsection{Two filtrations on the Yangian}\label{two_filtr_yang} 
The first filtration \(F_1\) on \(\yvg\) is determined by putting \(\deg_1 t_{ij}^{(r)}=r\). More precisely, the \(r\)-th filtered component \(F_1^{(r)}\yvg\) is the linear span of all monomials \(t_{i_1j_1}^{(r_1)}\cdot\ldots\cdot t_{i_mj_m}^{(r_m)}\) with \(r_1+\ldots+r_m\leqslant r\). 
It follows from Proposition~\ref{iso_Y_Y} that after the identification $Y(\mathfrak{g})=Y_V(\mathfrak{g})$ 
we have $\on{deg}_1 x=1,\, \on{deg}_1 J(x)=2$ for every nonzero $x \in \mathfrak{g}$.
By \(\on{gr_{1}}\) we denote the operation of taking associated graded algebra with respect to \(F_1\). For $x \in Y(\mathfrak{g})$ we denote by $\on{gr}_1 x$ the class of $x$ in  in $F_1^{(i)}Y(\mathfrak{g})/F_1^{(i-1)}Y(\mathfrak{g}) \subset \on{gr}_{2}Y(\mathfrak{g})$, where $i$ is the minimal such that $x \in F_1^{(i)}Y(\mathfrak{g})$.

Let $G[[t^{-1}]]$ be the group of $\BC[[t^{-1}]]$-points of $G$ i.e. $g \in G[[t^{-1}]]$ is a morphism $g\colon \on{Spec}\BC[[t^{-1}]] \ra G$.
For $g \in G[[t^{-1}]]$, we denote by 
$ev_g\colon \BC[G] \ra \BC[[t^{-1}]]$
the corresponding homomorphism of algebras.
We have the evaluation at infinity homomorphism $G[[t^{-1}]] \ra G$ and denote by $G_1[[t^{-1}]]$ its kernel. To any $f \in \BC[G]$ one can assign the function 
\begin{equation*}
\tilde{f}\colon G_1[[t^{-1}]] \ra \BC[[t^{-1}]],\, \tilde{f}=\sum_{r \geqslant 0} f^{(r)}t^{-r}
\end{equation*}
given by 
\begin{equation*}
\tilde{f}(g):=ev_g(f),\, g \in G_1[[t^{-1}]].
\end{equation*}
The $\BC^\times$-action on $G_1[[t^{-1}]]$ by dilations  of the variable $t$ determines a grading on $\BC(G_1[[t^{-1}]])$
such that $\on{deg}f^{(r)}=r$ for any $f \in \BC[G]$. 
\prop(\cite[Proposition~2.24]{ir})\label{gr_1_yanfian}
There is an isomorphism of graded algebras
$\on{gr}_{1} Y_V(\mathfrak{g}) \simeq \CO(G_1[[t^{-1}]])$, such that $\on{gr}_{1}t_{ij}^{(r)}=\Delta_{ij}^{(r)}$, where $\Delta_{ij} \in \BC[G]$ are the matrix coefficients of the representation $V$.
\eprop

The second filtration \(F_2\) on \(\yvg\) is determined by putting \(\deg t_{ij}^{(r)}=r-1\) i.e. the \(r\)-th filtered component \(F_2^{(r)}\yvg\) is the linear span of all monomials \(t_{i_1j_1}^{(r_1)}\cdot\ldots\cdot t_{i_mj_m}^{(r_m)}\) with \(r_1+\ldots+r_m\leqslant r+m\). 
It follows from Proposition~\ref{iso_Y_Y} that after the identification $Y(\mathfrak{g})=Y_V(\mathfrak{g})$ the filtration $F_2$ can by described as follows: we have $\on{deg}_2 x=0,\, \on{deg}_2 J(x)=1$ for every nonzero $x \in \mathfrak{g}$.
By $\on{gr}_2$ we denote the operation of taking associated graded algebra with respect to \(F_2\). For $x \in Y(\mathfrak{g})$ we denote by $\on{gr}_{2} x$ the class of $x$ in $F_2^{(i)}Y(\mathfrak{g})/F_2^{(i-1)}Y(\mathfrak{g}) \subset \on{gr}_{2}Y(\mathfrak{g})$, where $i$ is the minimal such that $x \in F_2^{(i)}Y(\mathfrak{g})$.

\prop{}(\cite[Theorem~6.5]{w})
\(\on{gr_{2}}\yvg\simeq U(\fg[t])\), where the grading is given by the $\BC^\times$-action dilating \(t\). Moreover, we have \(t^{r-1}\fg\subset\on{span}\{t_{ij}^{(r)}\}/F_2^{(r-2)}\yvg\). 
\eprop

\ssec{}{Filtrations on $\on{gr}_1 Y(\mathfrak{g})$, $\on{gr}_2 Y(\mathfrak{g})$ and the associated bigraded algebra}
We follow \cite[Sections 2.7, 2.13]{ir2}. The filtration $F_1$ on $Y(\mathfrak{g})$ produces a filtration on $U(\mathfrak{g}[t])=\on{gr}_2 Y(\mathfrak{g})$ which we denote by the same letter $F_1$. The  filtration $F_2$ on $Y(\mathfrak{g})$ produces the filtration on $\BC[G[[t^{-1}]]_1]=\on{gr}_1 Y(\mathfrak{g})$. The goal of this section is to describe these filtrations explicitly and also describe the corresponding associated bigraded algebra $\on{bigr}Y(\mathfrak{g})$ (see (\ref{bigr_defen}) below).

\sssec{}{Algebras with multiple filtrations} We begin with some general facts about algebras with multiple filtrations. For any algebra $A$ endowed with two filtrations $F_1$, $F_2$ one can define the associated bigraded algebra of $A$ as 
\begin{equation}\label{bigr_defen}
\on{bigr}A := \bigoplus_{i,j} (F_1^{(i)}A \cap F_2^{(j)}A)/((F_1^{(i-1)}A \cap F_2^{(j)}A) + (F_1^{(i)}A \cap F_2^{(j-1)}A)).    
\end{equation}

Algebra $\on{bigr}A$ is naturally bigraded. For $x \in F_1^{(i)}A \cap F_2^{(j)}A$ let $\on{bigr}^{(i,j)}(x) \in \on{bigr} A$ be the class of $x$ in $F_1^{(i)}A \cap F_2^{(j)}A/(F_1^{(i-1)}A \cap F_2^{(j)}A + F_1^{(i)}A \cap F_2^{(j-1)}A)$. 
 Note that $\on{bigr}^{(i,j)}(x)$ has bidegree $(i,j)$.
 
\war{} 
Note that $\on{bigr}^{(i,j)}(x)$ may be zero.
\ewar

 We have canonical identifications (see \cite[Lemma 2.8]{ir2})
 \begin{equation*}
\on{gr}_{2}\on{gr}_1 A \simeq \on{bigr}A \simeq \on{gr}_{1}\on{gr}_2A
 \end{equation*}
 so it makes sense to compare elements $\on{gr}_{21}x :=\on{gr}_2 \on{gr}_1 x$, $\on{gr}_{12}x := \on{gr}_1 \on{gr}_2 x$.
 The following Lemma describes necessary and sufficient conditions on $x \in A$ for the equality $\on{gr}_{21}x=\on{gr}_{12}x$ to be true.

\lem{}\label{cond_bigr_eq}
 $(a)$ Pick $x \in A$, $i,j$ such that $x \in F^{(i)}_1A \cap F_2^{(j)}A$.

 The following are equivalent:

$(\on{i})$ We have $x \notin (F_1^{(i-1)}A \cap F_2^{(j)}A)+(F_1^{(i)}A \cap F_2^{(j-1)}A)$, 

$(\on{ii})$ $\on{gr}_{21}x$ has bidegree $(i,j)$,

$(\on{iii})$ $\on{gr}_{21}x=\on{bigr}^{(i,j)}x$,

$(\on{iv})$ $\on{gr}_{12} x$ has bidegree $(i,j)$,

$(\on{v})$ $\on{gr}_{12}x=\on{bigr}^{(i,j)}x$,

$(\on{vi})$ $\on{gr}_{21}x=\on{bigr}^{(i,j)}x=\on{gr}_{12}x.$

$(b)$ Pick $x \in A$ then $\on{gr}_{21}x=\on{gr}_{12}x$ if and only if there exist $i,j$ such that $x \in F_1^{(i)}A \cap F_2^{(j)}A$ and one of six (equivalent) conditions $(i)$-$(vi)$ holds.

\elem
\begin{proof}
Let us prove part $(a)$.
Let us first of all show that 
\begin{equation}\label{equiv_cond_gr_21_is_bigr}
\on{bigr}^{(i,j)}x = \on{gr}_{21}x~\text{
if and only if}~x \notin (F_1^{(i)}A \cap F_2^{(j-1)}A)+F_1^{(i-1)}A.
\end{equation}

Assume that $x \notin (F_1^{(i)}A \cap F_2^{(j-1)}A)+F_1^{(i-1)}A$. Note that by our assumptions we have $x \in F_{1}^{(i)}A \setminus F_1^{(i-1)}A$ so $\on{gr}_1 x$ is equal to the class of $x$ in $F_1^{(i)}A/F_1^{(i-1)}A$. We need to show that 
\begin{equation}\label{prop_1_gr}
\on{gr}_1 x \in (F_1^{(i)}A \cap F_2^{(j)}A)/(F_1^{(i-1)}A \cap F_2^{(j)}A),\end{equation}
\begin{equation}\label{prop_2_gr}
\on{gr}_1 x \notin (F_1^{(i)}A \cap F_2^{(j-1)}A)/(F_1^{(i-1)}A \cap F_2^{(j-1)}A).
\end{equation}
Equation (\ref{prop_1_gr}) holds since $x \in F_1^{(i)}A \cap F_2^{(j)}A$, $(\ref{prop_2_gr})$ holds since $x \notin (F_1^{(i)}A \cap F_{2}^{(j-1)}A)+F_1^{(i-1)}A$.

Assume that $\on{bigr}^{(i,j)}x=\on{gr}_{21}x$. It follows that (\ref{prop_2_gr}) holds so $x \notin (F_1^{(i)}A \cap F_{2}^{(j-1)}A)+F_1^{(i-1)}A$ as desired. 

Same observation shows that 
\begin{equation}\label{equiv_cond_gr_12_is_bigr}
\on{bigr}^{(i,j)}x = \on{gr}_{12}x~\text{
if and only if}~x \notin (F_1^{(i-1)}A \cap F_2^{(j)}A)+F_2^{(j-1)}A.
\end{equation}

Note now that $x \in F_1^{(i)}A \cap F_2^{(j)}A$ implies that 
\begin{equation*}
x \notin (F_1^{(i-1)}A \cap F_2^{(j)}A)+(F_{1}^{(i)}A \cap F_2^{(j-1)}A)~\text{
if and only if}~x \notin (F_1^{(i)}A \cap F_2^{(j-1)}A)+F_1^{(i-1)}A
\end{equation*}
and 
\begin{equation*}
x \notin (F_1^{(i-1)}A \cap F_2^{(j)}A)+(F_{1}^{(i)}A \cap F_2^{(j-1)}A)~\text{
if and only if}~x \notin (F_1^{(i-1)}A \cap F_2^{(j)}A)+F_2^{(j-1)}A.
\end{equation*}
We conclude that $(\on{i}) \Leftrightarrow (\on{iii}) \Leftrightarrow (\on{v}) \Leftrightarrow (\on{vi})$. It is also clear that $(\on{ii}) \Leftrightarrow (\on{iii})$, $(\on{iv}) \Leftrightarrow (\on{v})$.

Let us now prove part $(b)$. We only need to show that if $\on{gr}_{21}x=\on{gr}_{12}x$ has some bidegree $(i,j)$ then $x \in F^{(i)}_1A \cap F^{(j)}_2A$. Indeed, since $\on{gr}_1x$ has degree $i$ (w.r.t. the grading on $\on{gr}_1A$ induced by $F_1$) we must have $x \in F_1^{(i)}A$ and similarly since $\on{gr}_2x$ has degree $j$ we must have $x \in F_2^{(j)}A$.

\end{proof}

\begin{warning}{}
Note that it is not true in general that $\on{gr}_{21}x=\on{gr}_{12}x$ for every $x \in A$. Take, for example, $A=\BC[a,b]$ and define  filtrations as follows:
\begin{equation*}
\on{deg}_{F_1}(a)=1,\, \on{deg}_{F_1}(b)=0,\, \on{deg}_{F_2}(a)=0,\, \on{deg}_{F_2}(b)=1.    
\end{equation*}
Take $x=a+b$ then 
\begin{equation*}
\on{gr}_{21}(x)= a \neq b = \on{gr}_{12}(x).    
\end{equation*}
\end{warning}

\sssec{}{Case of the Yangian} Let us now return to the case $A=Y(\mathfrak{g})$. Recall that $\on{gr}_1Y(\mathfrak{g})=\BC[G[[t^{-1}]]_1]$, $\on{gr}_2 Y(\mathfrak{g})=U(\mathfrak{g}[t])$. Let us describe the induced filtrations on $\BC[G[[t^{-1}]]_1]$, $U(\mathfrak{g}[t])$ and the associated bigraded algebra $\on{gr}_{21}Y(\mathfrak{g}) \simeq \on{bigr} Y(\mathfrak{g}) \simeq \on{gr}_{12}Y(\mathfrak{g})$.

Pick an identification $\on{exp}\colon t^{-1}\mathfrak{g}[[t^{-1}]] \iso G_1[[t^{-1}]]$ and identify $\BC[t^{-1}\mathfrak{g}[[t^{-1}]]] \simeq S^\bullet(\mathfrak{g}[t])$ via the pairing given by 
\begin{equation*}
(x(t),y(t)):=\on{Res}_{t=0}( x(t),y(t)),~x(t) \in \mathfrak{g}[t],\, y(t) \in t^{-1}\mathfrak{g}[[t^{-1}]].    
\end{equation*}

The following proposition holds by \cite[proof of Proposition 2.12, Section 2.13]{ir2}.
\prop{}\label{expl_descr_filtr_bigr_prop}
$(a)$ After an identification $\on{gr}_1 Y(\mathfrak{g}) \simeq \BC[G_1[[t^{-1}]]] \simeq S^\bullet(\mathfrak{g}[t])$ above the grading on $S^\bullet(\mathfrak{g}[t])$ is given by $\on{deg}_1x[n-1]=n$ and the filtration $F_2$ is given by $\on{deg}_2x[n-1]=n$.

$(b)$ The grading on $\on{gr}_2 Y(\mathfrak{g}) \simeq U(\mathfrak{g}[t])$ is given by $\on{deg}_2x[n-1]=n-1$ and the filtration $F_1$ is given by $\on{deg}_1 x[n-1]=n$.

$(c)$ We have a bigraded algebra isomorphism $\on{bigr}Y(\mathfrak{g}) \simeq S^\bullet(\mathfrak{g}[t]) 
$, where the bigrading on $S^\bullet(\mathfrak{g}[t])$ is given by 
$\on{deg}_1 x[n-1]=n$, $\on{deg}_2 x[n-1]=n-1$.
\eprop

\ssec{}{Representation theory of $Y(\mathfrak{g})$}\label{rep_th_yang_facts}
By~\cite[Section~12]{cp} to every dominant $\la$ we can associate an irreducible finite-dimensional representation of $Y(\mathfrak{g})$ to be denoted $V(\la,0)$, $\rho\colon Y(\mathfrak{g}) \ra \on{End}V(\la,0)$. 
As a $\mathfrak{g}$-module we have $V(\la,0)=V_\la \oplus \bigoplus_{\mu<\la}V_\mu^{\oplus l_\mu}$. Moreover, if $\la$ is minuscule then $V(\la,0)=V_\la$ and $J(x)$ acts on $V(\la,0)$ by zero for every $x \in \mathfrak{g}$.

\begin{Rem}
{\em{If $\mathfrak{g}$ is of type $A$ then $V(\la,0)=V_\la$ for every dominant $\la$. The reason for this is the existence of so-called evaluation homomorphism $Y(\mathfrak{g}) \ra U(\mathfrak{g})$ (which exists only in type $A$), see Section \ref{ev_hom_for_y_and_u_section}.}}
\end{Rem}




\section{Bethe subalgebras and their degenerations}\label{bethe}

\subsection{Bethe subalgebras in the Yangian}
Let \(\rho_i\colon\yg\to\End V(\varpi_i,0)\) be the \(i\)-th fundamental representation of \(\yg\).
Let $V$ be the direct sum of $V(\varpi_i,0)$. 
Let 
\(T^i(u)
\) be the submatrix of \(T(u)\)-matrix, corresponding to the \(i\)-th fundamental representation $V(\varpi_i,0)$. Let $\tilde{G}$ be the corresponding to \(\fg\) connected simply-connected group and let $G$ be the adjoint group with Lie algebra $\mathfrak g$. The action of $\mathfrak{g} \subset Y(\mathfrak{g})$ integrates to the action $\tilde{G} \curvearrowright V(\varpi_i,0)$. We denote the corresponding map $\tilde{G}\ra \on{End}V(\varpi_i,0)$ by $\rho_i$.

\begin{defeni}{}
Let \(C\in \tilde{G}\). Bethe subalgebra \(B(C)\subset \yvg\) is the subalgebra generated by all coefficients of the following series with the coefficients in \(\yvg\)
\[
\tau_i(u,C):=\on{tr}_{V(\varpi_i,0)}\rho_i(C)T^i(u), 1\leqslant i\leqslant n.
\]
\end{defeni}

\begin{Rem}
{\em{In fact, \(B(C)\) depends only on the class of \(C\) in $G=\tilde{G}/Z(\tilde{G})$ so from now on we assume that $C \in G$.
}}
\end{Rem}

It is easy to see that $B(C) \subset Y_V(\mathfrak{g})^{\mathfrak{z}_\mathfrak{g}(C)}$. Let $T \subset G$ be a maximal torus.
\begin{Prop}{(\cite{no},\,\cite{i},\,\cite{ir},\,\cite{ir2})}\label{main_Bethe}
\begin{enumerate}
    \item Bethe subalgebra \(B(C)\) is commutative  for any \(C\in G\).
    \item \(B(C)\) is a maximal commutative subalgebra of \(\yvg\) for \(C\) in the regular part of $T$.
    \item  \(B(C)\) is a maximal commutative subalgebra of $Y_V(\mathfrak{g})^{\mathfrak{z}_\mathfrak{g}(C)}$ for $C \in T$.
    \item For $C \in T$  we have $\on{gr}_{2} B(C)=\CA_{0}^{\mathrm{u}}(\mathfrak{z}_{\mathfrak{g}}(C))$ so, in particular, we have $\on{dim}_{F_1}B(C)=\on{dim}_{F_1}\CA_{0}^{\mathrm{u}}(\mathfrak{z}_{\mathfrak{g}}(\chi))$.
\end{enumerate}
\end{Prop}

\ssec{}{Limits of Bethe subalgebras}

Recall now the filtration $F^\bullet_2$ on $Y(\mathfrak{g})$. We can consider the corresponding Rees algebra (see Section~\ref{rees_constr}) $Rees(Y(\mathfrak{g}))=Y_{\hbar}(\mathfrak{g})$. For every $\epsilon \in \BC$ we obtain the algebra $Y_{\epsilon}(\mathfrak{g}):=Y_{\hbar}(\mathfrak{g})/(\hbar-\epsilon)$. Recall also that for $\epsilon \neq 0$ we have the identification 
\begin{equation}\label{y_eps_via_y}
Y_{\epsilon}(\mathfrak{g})  \iso Y(\mathfrak{g}),\, [\hbar^i a]  \mapsto    \epsilon^i a .
\end{equation}
We denote the inverse to~(\ref{y_eps_via_y}) by $\psi_\epsilon$. Pick now an element $\chi \in \mathfrak{h}$ and consider the element $C:=\on{exp}(\epsilon \chi) \in T$. We can consider the corresponding Bethe subalgebra $B(\on{exp}(\epsilon \chi)) \subset Y(\mathfrak{g})$. Let us now denote by $B_{\epsilon}(C)=B_{\epsilon}(\on{exp}(\epsilon \chi))$ the subalgebra $\psi_\epsilon(B(\on{exp}(\epsilon \chi))) \subset Y_{\epsilon}(\mathfrak{g})$. We obtain the (formal) algebraic family $B_\epsilon(C)$ of subalgebras in $Y_\epsilon(\mathfrak{g})$. 

\rem{}
{\em{Note that $B_\epsilon(\on{exp}(\epsilon\chi))$ can be considered as an algebraic family over $\on{Spec}\BC[[\hbar]]$ or as a complex analytic family over $\BC$.}}
\erem

Using the filtration $F_1^\bullet$, we can then define the limit $\underset{\epsilon \rightarrow 0}{\on{lim}}\,B_\epsilon(C)$ (see Section~\ref{rees_constr}) that will be a commutative subalgebra of $Y_0(\mathfrak{g})=\on{gr}_{2} Y(\mathfrak{g})=U(\mathfrak{g}[t])$.

The goal of this section is to prove that $\underset{\epsilon \rightarrow 0}{\on{lim}}\,B_\epsilon(C)=\CA^{\mathrm{u}}_\chi$. The strategy is the following. Recall that by Proposition~\ref{ca_as_centr} we have $\CA^{\mathrm{u}}_\chi=Z_{U(\mathfrak{g}[t])^{\mathfrak{z}_{\mathfrak{g}}(\chi)}}(\tilde{\Omega}_\chi)$. We prove in Proposition~\ref{W_lim} that $\tilde{\Omega}_\chi \in \underset{\epsilon \rightarrow 0}{\on{lim}}\,B_\epsilon(C)$, concluding that 
$\underset{\epsilon \rightarrow 0}{\on{lim}}\,B_\epsilon(C) \subset Z_{U(\mathfrak{g}[t])^{\mathfrak{z}_{\mathfrak{g}}(\chi)}}(\tilde{\Omega}_\chi)=\CA^{\mathrm{u}}_\chi$. Then the comparison of the dimensions finishes the proof.

Let $\al_j,\, j=1,\ldots,\on{rk}\mathfrak{g}$ be simple roots of $\mathfrak{g}$ and recall that $\varpi_j,\, j=1,\ldots,\on{rk}\mathfrak{g}$ are fundamental weights. 
We also denote by $h_j \in \mathfrak{h}$ the element, corresponding to the simple root $\al_j \in \mathfrak{h}^*$ via the identification $\mathfrak{h}^* \iso \mathfrak{h}$ induced by the invariant scalar product $(\,,\,)$. Similarly, $t_{\varpi_j} \in \mathfrak{h}$ is the element, corresponding to $\varpi_j \in \mathfrak{h}^*$. For every positive root $\al \in \Delta_+$ we denote by $x_{\al}^{\pm} \in \mathfrak{g}_{\pm \al}$ elements of $\mathfrak{g}_{\pm \al}$ such that $(x_\al^+,x_{\al}^-)=1$. 
We have 
\begin{equation}\label{expl_form_for_R}
\CR^{(2)}=\sum_{\al \in \Delta_+}(J(x_\al^{\pm}) \otimes x_\al^{\mp} - x_\al^{\pm} \otimes J(x_\al^{\mp}))+\sum_{j=1}^{\on{rk}\mathfrak{g}}\frac{2}{(\al_j,\al_j)}(J(h_j)\otimes t_{\varpi_j}-h_j \otimes J(t_{\varpi_j}))+2\omega'^2,
\end{equation}
where 
\begin{equation*}
\omega'=\sum_{\al \in \Delta^+}\frac{1}{2}x_{\al}^{\pm} \otimes x_{\al}^{\mp}+\sum_{j=1}^{\on{rk}\mathfrak{g}} \frac{1}{(\al_j,\al_j)}h_j \otimes t_{\varpi_j}.
\end{equation*}

\prop{}\label{lim_w}
The element $\tilde{\omega}_\chi$ lies in the limit algebra $\underset{\epsilon \rightarrow 0}{\on{lim}}\,B_\epsilon(\on{exp}(\epsilon \chi))$.
\eprop
\prf
Let $V$ be a finite-dimensional representation of $Y(\mathfrak{g})$ that is not a direct sum of trivial $\mathfrak{g}$ modules. Let $(\,,\,)_V$ be the corresponding invariant form and let $c_V \in \BC^\times$ be such that $(\,,\,)_V=c_V(\,,\,)$. 

Pick an identification $\on{exp}\colon t^{-1}\mathfrak{g}[[t^{-1}]] \iso G_1[[t^{-1}]]$. Then (using Proposition \ref{gr_1_yanfian}) we obtain an isomorphism $\on{gr}_{1}Y(\mathfrak{g}) \simeq \BC[t^{-1}\mathfrak{g}[[t^{-1}]]]$, which sends  $\on{gr}_1(\on{tr}_{V}\rho(C))T$ to the function on $t^{-1}\mathfrak{g}[[t^{-1}]]$ given by $g(t^{-1}) \mapsto \on{tr}_{V}\rho(C)\on{exp}g(t^{-1})$, 
where $g(t^{-1})=a_1t^{-1}+a_2t^{-2}+\ldots$, $a_i \in \mathfrak{g}$.
We begin with two Lemmas.
\lem{}\label{gr_T_2}
We have $\on{tr}_V T^{(2)} \in F_2^{(0)}Y(\mathfrak{g})$ and $\on{gr}_{2}(\on{tr}_V T^{(2)})=c_V\tilde{\omega}_0+b$ for some $b \in \BC$. 
\elem
\prf

Assume for the sake of contradiction that $\on{tr}_V T^{(2)} \notin F^{(0)}_2 Y(\mathfrak{g})$. It follows from the definitions that $\on{tr}_V T^{(2)} \in F_1^{(2)}Y(\mathfrak{g}) \cap F^{(1)}_2 Y(\mathfrak{g})$. We claim that the condition $(\on{i})$ of Lemma \ref{cond_bigr_eq} holds for $x=\on{tr}_V T^{(2)}$ and $(i,j)=(2,1)$. Indeed, note that $F_1^{(1)}Y(\mathfrak{g}) \subset F_2^{(0)}Y(\mathfrak{g})$ so $(F_1^{(1)}Y(\mathfrak{g}) \cap F_2^{(1)}Y(\mathfrak{g}))+(F_1^{(2)}Y(\mathfrak{g}) \cap F_2^{(0)}Y(\mathfrak{g})) \subset F_2^{(0)}Y(\mathfrak{g})$ and $\on{tr}_V T^{(2)} \notin F^{(0)}_2 Y(\mathfrak{g})$ by our assumption. It follows from Lemma \ref{cond_bigr_eq} that
$
\on{gr}_{21}(\on{tr}_V T^{(2)})=\on{bigr}^{(2,1)}(\on{tr}_V T^{(2)})$ has bidegree $(2,1)$.
Note also that $\on{gr}_1(\on{tr}_V T^{(2)})$ is the function 
\begin{equation*}
g \mapsto \frac{1}{2}(a_1,a_1)_V    
\end{equation*}
that has degree $0<1$ w.r.t. the second filtration on $\BC[t^{-1}\mathfrak{g}[[t^{-1}]]] \simeq S^\bullet(\mathfrak{g}[t])$ (use Proposition \ref{expl_descr_filtr_bigr_prop}). Contradiction finishes the proof.

Let us now prove that $\on{gr}_{2}(\on{tr}_V T^{(2)})=c_V\tilde{\omega}_0+b$ for some $b \in \BC$. 
Note that $\on{tr}_V T^{(2)}$ is $\mathfrak{g}$-invariant and lies in $F_1^{(2)}Y(\mathfrak{g}) \cap F_2^{(0)}Y(\mathfrak{g})$ so its class in $\on{gr}_{2}Y(\mathfrak{g})=U(\mathfrak{g}[t])$ must be a $\mathfrak{g}$-invariant element of $U(\mathfrak{g})$ that has degree at most two w.r.t. the $PBW$-filtration on $U(\mathfrak{g})$. It follows that $\on{gr}_{2}(\on{tr}_V T^{(2)})=a\tilde{\omega}_0+b$, $a,b \in \BC$. 

It remains to check that $a = c_V$. Let us first of all prove that 
\begin{equation}\label{bigr_t_2}
\on{gr}_{21}(\on{tr}_V T^{(2)})=\on{bigr}(\on{tr}_V T^{(2)})=\on{gr}_{12}(\on{tr}_V T^{(2)}).    
\end{equation}
To see that, we apply Lemma \ref{cond_bigr_eq} (for $i=2, j=0$): we have $\on{tr}_V T^{(2)} \in F_1^{(2)}Y(\mathfrak{g}) \cap F_2^{(0)}Y(\mathfrak{g})$, note also that
$\on{gr}_1\on{tr}_V T^{(2)}$ is $(g \mapsto \frac{1}{2}(a_1,a_1)_V)$ so (using Proposition \ref{expl_descr_filtr_bigr_prop}) $\on{gr}_{21}\on{tr}_V T^{(2)}=c_V \omega_0$ has bidegree $(2,1)$. 
We conclude that (\ref{bigr_t_2}) holds. It remains to note that 
\begin{equation*}
a\omega_0=\on{gr}_1(a\tilde{\omega}_0+b)=\on{gr}_{12}\on{tr}_V T^{(2)}=\on{gr}_{21}\on{tr}_V T^{(2)}=\on{gr}_{2}(x \mapsto \frac{1}{2}(a_1,a_1)_V)=c_V\omega_0
\end{equation*}
so $a=c_V$ as desired.
\epr

\lem{}\label{gr_T_2_chi}
Assume that $\chi \neq 0$. 

$(a)$ We have $\on{tr}_V \rho(\chi)T^{(2)} \in F_1^{(2)}Y(\mathfrak{g}) \cap F_2^{(1)}Y(\mathfrak{g})$ and $\on{tr}_V \rho(\chi)T^{(2)} \notin F_2^{(0)}Y(\mathfrak{g})$.

$(b)$ We have $\on{gr}_2\Big(\frac{1}{c_V}\on{tr}_V \rho(\chi)T^{(2)}\Big)=\chi[1]$.
\elem
\prf 
It follows from the definitions that $\on{tr}_V \rho(\chi)T^{(2)} \in F_1^{(2)}Y(\mathfrak{g}) \cap F_2^{(1)}Y(\mathfrak{g})$. Let us now check that $\on{tr}_V \rho(\chi)T^{(2)} \notin  F_2^{(0)}Y(\mathfrak{g})$. Assume for the sake of contradiction that $\on{tr}_V \rho(\chi)T^{(2)} \in  F_2^{(0)}Y(\mathfrak{g})$. It follows that the degree of $\on{gr}_1(\on{tr}_V \rho(\chi)T^{(3)})$ w.r.t. the second filtration on $\BC[t^{-1}\mathfrak{g}[[t^{-1}]]]=S^\bullet(\mathfrak{g}[t])$ is at most zero. This contradicts to the fact that $\on{gr}_1(\on{tr}_V \rho(\chi)T^{(2)})$ is
\begin{equation*}
x \mapsto \frac{1}{2}(a_1,a_1)_V+(\chi,a_2)_V,    
\end{equation*}
which has degree one with respect to the second filtration (use that $\chi \neq 0$). 

Let us now decompose $\chi=\sum_{j=1}^{\on{rk}\mathfrak{g}}\frac{2(\chi,\al_j)}{(\al_j,\al_j)}t_{\varpi_j}$ and consider the element 
\begin{equation*}
x:=\sum_{j=1}^{\on{rk}\mathfrak{g}}\frac{2(\chi,\al_j)}{(\al_j,\al_j)}J(t_{\varpi_j}) \in Y(\mathfrak{g}),
\end{equation*}
which lies in $F_1^{(2)}Y(\mathfrak{g}) \cap F_2^{(1)}Y(\mathfrak{g})$. It is clear that  $\on{gr}_2(x)=\chi[1]$. 
Note now that $x \notin F_2^{(0)}$ (since $\on{gr}_2(x)=\chi[1]$). Recall also that $\on{tr}_V\rho(\chi)T^{(2)} \notin F_2^{(0)}$. We conclude from Lemma \ref{cond_bigr_eq} (using  $F_1^{(1)}Y(\mathfrak{g}) \subset F_2^{(0)}Y(\mathfrak{g})$) that 
\begin{equation*}
\on{gr}_{21}(x)=\on{bigr}^{(2,1)}(x)=\on{gr}_{12}(x),
\end{equation*}
\begin{equation*}
\on{gr}_{21}(\on{tr}_V\rho(\chi)T^{(2)})=\on{bigr}^{(2,1)}(\on{tr}_V\rho(\chi)T^{(2)})=\on{gr}_{12}(\on{tr}_V\rho(\chi)T^{(2)}).    
\end{equation*}
Note also that 
\begin{equation*}
\on{bigr}\Big(\frac{1}{c_V}\on{tr}_V\rho(\chi)T^{(2)}\Big)=\on{gr}_{21}\Big(\frac{1}{c_V}\on{tr}_V\rho(\chi)T^{(2)}\Big)=\chi[1]=\on{gr}_{12}(x)=\on{bigr}(x)    
\end{equation*}
so 
\begin{equation*}
\frac{1}{c_V}\on{tr}_V\rho(\chi)T^{(2)} - x \in (F_{1}^{(1)}Y(\mathfrak{g}) \cap F_2^{(1)}Y(\mathfrak{g})) +  (F_1^{(2)}Y(\mathfrak{g}) \cap F_{2}^{(0)}Y(\mathfrak{g}))  \subset F_2^{(0)}Y(\mathfrak{g}), 
\end{equation*}
hence, 
\begin{equation*}
\on{gr}_2\Big(\frac{1}{c_V}\on{tr}_V\rho(\chi)T^{(2)}\Big)=\on{gr}_2(x)=\chi[1].
\end{equation*}
\epr

We are now ready to prove Proposition \ref{lim_w}.



We start from the case $\chi=0$. In this case we need to show that $\tilde{\omega}_0 \in \on{gr}_{2}B(1)$. 
It follows from Lemma \ref{gr_T_2} that $\on{gr}_2(\on{tr}_V T^{(2)}) = c_V\tilde{\omega}_0 + b$ for some $b \in \BC$. It follows that 
\begin{equation*}
\tilde{\omega}_0=\frac{1}{c_V}\Big(\on{gr}_2(\on{tr}_V T^{(2)})-b\Big) \in \on{gr}_{2}B(1).    
\end{equation*}

Let us now deal with arbitrary $\chi \neq 0$.
Consider the image of the  element 
\begin{equation*}
\frac{1}{c_V}\Big(\on{tr}_V\rho_V(\on{exp}(\epsilon \chi))T^{(2)}-b\Big)
\end{equation*}
in $B_{\epsilon}(C)$. 
Let $X(\chi)$ be the coefficient in front of $\epsilon^{0}=1$, which we consider as a function on $\chi$ with values in $Y_0(\mathfrak{g})=U(\mathfrak{g}[t])$. 
It follows from Lemmas \ref{gr_T_2}, \ref{gr_T_2_chi}  that
we have 
\begin{equation*}
\psi_\epsilon\Big(\frac{1}{c_V}\Big(\on{tr}_V\rho_V(\on{exp}(\epsilon \chi))T^{(2)}-b\Big) \Big) = X(\chi) + O(\epsilon).   
\end{equation*}
We claim that 
$X(\chi)$ is exactly $\tilde{\omega}_\chi$. It follows from the above
that $X(0)=\tilde{\omega}_0$.  Recall that terms in $X(\chi)$, depending on $\chi$ (i.e. $X(\chi)-X(0)$), appear only from
\begin{equation*}
\psi_\epsilon(\on{tr}_{V}\rho_V(\epsilon \chi+\ldots)T^{(2)}).
\end{equation*}
Recall now that $X(\chi)$ is the coefficient in front of $1$ so the only term in $X(\chi)$, depending on $\chi$, comes from $\on{tr}_{V}\rho_V(\chi)T^{(2)}$ i.e. $X(\chi)-X(0)= \frac{1}{c_V}\on{gr}_2(\on{tr}_{V}\rho_V(\chi)T^{(2)})$. Recall that by Lemma \ref{gr_T_2_chi} we have $\frac{1}{c_V}\on{gr}_2(\on{tr}_{V}\rho_V(\chi)T^{(2)})=\chi[1]$. We conclude that 
\begin{equation*}
X(\chi)=X(0)+\chi[1]=\tilde{\omega}_{\chi}    
\end{equation*}
so $\tilde{\omega}_{\chi} \in \underset{\epsilon \rightarrow 0}{\on{lim}}\,B_\epsilon(\on{exp}(\epsilon \chi))$.
%
\epr

\begin{Rem}
{\em{Note that for $\mathfrak{g}$, which has a minuscule representation (i.e. for $\mathfrak{g}$ not of type $E_8$, $G_2$) it is easy to see that the element $\tilde{\omega}_\chi$ lies in the limit, using the explicit formula (\ref{expl_form_for_R}) for $\CR^{(2)}$.}}  
\end{Rem}

We also need to find an element 
\begin{equation*}
\tilde{\Omega}_\chi=\sum_{a}x_ax_a[1]+\chi[2]
\end{equation*}
in the limit $\underset{\epsilon \ra 0}{\on{lim}}\,B_{\epsilon}(C)$. 

\prop{}\label{W_lim}
The element $\tilde{\Omega}_\chi$ lies in the limit algebra $\underset{\epsilon \rightarrow 0}{\on{lim}}\,B_\epsilon(\on{exp}(\epsilon \chi))$.
\eprop
\prf 
The proof is similar to the proof of Proposition \ref{lim_w}. Recall that $V$ is a finite-dimensional representation of $Y(\mathfrak{g})$ that is not a direct sum of trivial $\mathfrak{g}$ modules, $(\,,\,)_V$ is the corresponding invariant form and  $c_V \in \BC^\times$ is such that $(\,,\,)_V=c_V(\,,\,)$.  We begin with two Lemmas (cf. Lemmas \ref{gr_T_2}, \ref{gr_T_2_chi} above).
\lem{}\label{gr_T_3}
$(a)$ 
We have $\on{tr}_V T^{(3)} \in F_2^{(1)}Y(\mathfrak{g})$. 

$(b)$ 
We have $\on{gr}_{2}(\on{tr}_V T^{(3)})=c_V\tilde{\Omega}_0$. 
\elem
\prf
Let us prove $(a)$. Assume for the sake of contradiction that $\on{tr}_V T^{(3)} \notin F^{(1)}_2 Y(\mathfrak{g})$. It follows from the definitions that $\on{tr}_V T^{(3)} \in F_1^{(3)}Y(\mathfrak{g}) \cap F^{(2)}_2  Y(\mathfrak{g})$. We claim that the condition $(\on{i})$ of Lemma \ref{cond_bigr_eq} holds for $x=\on{tr}_V T^{(3)}$ and $(i,j)=(3,2)$. Indeed, note that $F_1^{(2)}Y(\mathfrak{g}) \subset F_2^{(1)}Y(\mathfrak{g})$ so $(F_1^{(2)}Y(\mathfrak{g}) \cap F_2^{(2)}Y(\mathfrak{g}))+(F_1^{(3)}Y(\mathfrak{g}) \cap F_2^{(1)}Y(\mathfrak{g})) \subset F_2^{(1)}Y(\mathfrak{g})$ and $\on{tr}_V T^{(3)} \notin F^{(1)}_2 Y(\mathfrak{g})$ by our assumption. It follows from Lemma \ref{cond_bigr_eq} that
\begin{equation*}
\on{gr}_{21}(\on{tr}_V T^{(3)})=\on{bigr}^{(3,2)}(\on{tr}_V T^{(3)})=\on{gr}_{12}(\on{tr}_V T^{(3)}).    
\end{equation*}
Note now that $\on{gr}_1(\on{tr}_V T^{(3)})$ is the function 
\begin{equation*}
g \mapsto (a_1,a_2)_V + f(a_1),   
\end{equation*}
where $f \in (S^3\mathfrak{g})^{\mathfrak{g}}$ is some $\mathfrak{g}$-invariant function of degree $3$. We see that $\on{gr}_1(\on{tr}_V T^{(3)})$ has degree $1<2$ w.r.t. the second filtration on $\BC[t^{-1}\mathfrak{g}[[t^{-1}]]]$. Contradiction finishes the proof.

Let us now prove $(b)$. 
Let us
note that $\on{tr}_V T^{(3)} \in F_1^{(3)}Y(\mathfrak{g}) \cap F_2^{(1)}Y(\mathfrak{g})$ and 
\begin{equation*}
\on{gr}_{21}\on{tr}_V T^{(3)}=\on{gr}_2 (g \mapsto ((a_1,a_2)_V+f(a_1)))=c_V\Omega_0    
\end{equation*}
has bidegree $(3,1)$. It follows from Lemma \ref{cond_bigr_eq} that $c_V\Omega_0  =\on{gr}_{21}\on{tr}_V T^{(3)}=\on{gr}_{12}\on{tr}_V T^{(3)}$ so $\on{gr}_1(\on{gr}_2 \on{tr}_V T^{(3)})=c_V \Omega_0$. 
Note now that $\on{tr}_V T^{(3)}$ is $\mathfrak{g}$-invariant so $\on{gr}_2 \on{tr}_V T^{(3)} \in U(\mathfrak{g}[t])^{\mathfrak{g}}$. Moreover, $\on{gr}_2 \on{tr}_V T^{(3)}$ is homogeneous of degree $1$ w.r.t. the grading on $U(\mathfrak{g}[t])$.
There exists the unique graded lift of $c_V \Omega_0$ to $U(\mathfrak{g}[t])^{\mathfrak{g}}$ so $\on{gr}_2 \on{tr}_V T^{(3)} =c_V\tilde{\Omega}_0$.




\epr


\lem{}\label{gr_T_3_chi}
Assume that $\chi \neq 0$. 

$(a)$ We have $\on{tr}_V \rho(\chi)T^{(3)} \in F_1^{(3)}Y(\mathfrak{g}) \cap F_2^{(2)}Y(\mathfrak{g})$ and $\on{tr}_V \rho(\chi)T^{(3)} \notin F_2^{(1)}Y(\mathfrak{g})$.

$(b)$ We have $\on{gr}_2\Big(\frac{1}{c_V}\on{tr}_V \rho(\chi)T^{(3)}\Big)=\chi[2]$.
\elem
\prf
It follows from the definitions that $\on{tr}_V \rho(\chi)T^{(3)} \in F_1^{(3)}Y(\mathfrak{g}) \cap F_2^{(2)}Y(\mathfrak{g})$. Let us now check that $\on{tr}_V \rho(\chi)T^{(3)} \notin  F_2^{(1)}Y(\mathfrak{g})$. Assume for the sake of contradiction that $\on{tr}_V \rho(\chi)T^{(3)} \in  F_2^{(1)}Y(\mathfrak{g})$. It follows that the degree of $\on{gr}_1(\on{tr}_V \rho(\chi)T^{(3)})$ w.r.t. the second filtration on $\BC[t^{-1}\mathfrak{g}[[t^{-1}]]]$ is at most one. This contradicts to the fact that $\on{gr}_1(\on{tr}_V \rho(\chi)T^{(3)})$ is
\begin{equation*}
g \mapsto (\chi,a_3)_V+(a_1,a_2)_V+f(a_1)    
\end{equation*}
that has degree two with respect to the second filtration (use Proposition \ref{expl_descr_filtr_bigr_prop}). 

Let us now decompose $\chi=\sum_{j=1}^{\on{rk}\mathfrak{g}}\frac{2(\chi,\al_j)}{(\al_j,\al_j)}t_{\varpi_j}$ and consider the element 
\begin{equation*}
x:=\sum_{j=1}^{\on{rk}\mathfrak{g}}\frac{2(\chi,\al_j)}{(\al_j,\al_j)}[J(x^+_{\al_j}),J(x^-_{\al_j})]
\end{equation*}
that clearly lies in $F_1^{(3)}Y(\mathfrak{g}) \cap F_2^{(2)}Y(\mathfrak{g})$ and $\on{gr}_2(x)=\chi[2]$. 
Note now that $x \notin F_2^{(1)}Y(\mathfrak{g})$ (since $\on{gr}_2(x)=\chi[2]$). Recall also that $\on{tr}_V\rho(\chi)T^{(3)} \notin F_2^{(1)}$. We conclude from Lemma \ref{cond_bigr_eq}  (using  $F_1^{(2)}Y(\mathfrak{g}) \subset F_2^{(1)}Y(\mathfrak{g})$) that 
\begin{equation*}
\on{gr}_{21}(x)=\on{bigr}^{(3,2)}(x)=\on{gr}_{12}(x),
\end{equation*}
\begin{equation*}
\on{gr}_{21}(\on{tr}_V\rho(\chi)T^{(3)})=\on{bigr}^{(3,2)}(\on{tr}_V\rho(\chi)T^{(3)})=\on{gr}_{12}(\on{tr}_V\rho(\chi)T^{(3)}).    
\end{equation*}
Note now that 
\begin{equation*}
\on{bigr}\Big(\frac{1}{c_V}\on{tr}_V\rho(\chi)T^{(3)}\Big)=\on{gr}_{21}\Big(\frac{1}{c_V}\on{tr}_V\rho(\chi)T^{(3)}\Big)=\chi[2]=\on{gr}_{12}(x)=\on{bigr}(x)    
\end{equation*}
so 
\begin{equation*}
\frac{1}{c_V}\on{tr}_V\rho(\chi)T^{(3)} - x \in (F_{1}^{(2)}Y(\mathfrak{g}) \cap F_2^{(2)}Y(\mathfrak{g})) +  (F_1^{(3)}Y(\mathfrak{g}) \cap F_{2}^{(1)}Y(\mathfrak{g}))  \subset F_2^{(1)}Y(\mathfrak{g}), 
\end{equation*}
hence, 
\begin{equation*}
\on{gr}_2\Big(\frac{1}{c_V}\on{tr}_V\rho(\chi)T^{(3)}\Big)=\on{gr}_2(x)=\chi[2].
\end{equation*}
\epr



%
We are now ready to prove Proposition \ref{W_lim}. 
We start from the case $\chi=0$. 
Consider the  element $\on{tr}_{V} T^{(3)} \in B(1)$ and recall that by Lemma \ref{gr_T_3} $(b)$ we have 
\begin{equation*}
\on{gr}_2(\on{tr}_{V} T^{(3)})=c_V\tilde{\Omega}_0.
\end{equation*}
We conclude that 
\begin{equation*}
\tilde{\Omega}_0=\frac{1}{c_V}\on{gr}_2(\on{tr}_{V} T^{(3)}) \in \on{gr}_2B(1).
\end{equation*}

Let us now deal with arbitrary $\chi$. Consider the image of the element 
\begin{equation*}
\frac{1}{c_V}\on{tr}_V \rho_V(\on{exp}(\epsilon \chi))T^{(3)}    
\end{equation*}
in $B_\epsilon(C)$. Let $X(\chi)$ be the coefficient in front of $\epsilon^{-1}$, which we consider as a function on $\chi$ with values in $Y_0(\mathfrak{g})=U(\mathfrak{g}[t])$. 
It follows from Lemmas \ref{gr_T_3}, \ref{gr_T_3_chi} that 
\begin{equation*}
\psi_\epsilon\Big(\frac{1}{c_V}\on{tr}_V \rho_V(\on{exp}(\epsilon \chi))T^{(3)}    \Big)= X(\chi) \epsilon^{-1}+ O(1).    
\end{equation*}
We claim that $X(\chi)$ is exactly $\tilde{\Omega}_\chi$. It follows from the above that $X(0)=\tilde{\Omega}_0$. Recall that terms in $X(\chi)$, depending on $\chi$ (i.e. $X(\chi)-X(0)$) appear only from 
\begin{equation*}
\psi_\epsilon (\on{tr}_V \rho_V(\epsilon \chi+\ldots ) T^{(3)}).    
\end{equation*}
Recall now that $X(\chi)$ is the coefficient in front of $\epsilon^{-1}$ so the only term in $X(\chi)$, depending on $\chi$ comes from $\on{tr}_V \rho_V(\chi)T^{(3)}$ i.e. $X(\chi)-X(0)=\frac{1}{c_V} \on{gr}_2(\on{tr}_V \rho_V(\chi)T^{(3)})$. Recall  that by Lemma \ref{gr_T_3_chi}  we have  $\frac{1}{c_V} \on{gr}_2(\on{tr}_V \rho_V(\chi)T^{(3)})=\chi[2]$. We conclude that 
\begin{equation*}
X(\chi)=X(0)+\chi[2]=\tilde{\Omega}_\chi
\end{equation*}
so $\tilde{\Omega}_\chi \in \underset{\epsilon \rightarrow 0}{\on{lim}}\,B_\epsilon(\on{exp}(\epsilon \chi))$.

\epr

\th{}\label{lim_Bethe}
We have $\underset{\epsilon \rightarrow 0}{\on{lim}}\,B_\epsilon(\on{exp}(\epsilon \chi))=\CA^{\mathrm{u}}_\chi$.
\eth
\prf
We set $C=C(\epsilon)=\on{exp}(\epsilon \chi)$.
Recall that  by Proposition~\ref{ca_as_centr} we have
$\CA_\chi^{\mathrm{u}}=Z_{U(\mathfrak{g}[t])}(\tilde{\Omega}_\chi)$.
It follows from Proposition
~\ref{W_lim} that 
$\tilde{\Omega}_\chi \in \underset{\epsilon \rightarrow 0}{\on{lim}}\,B_\epsilon(C)$. We conclude that $\underset{\epsilon \rightarrow 0}{\on{lim}}\,B_\epsilon(C) \subset \CA^{\mathrm{u}}_\chi$. It remains to note that by Lemma~\ref{dim_lim_Rees} 
\begin{equation*}
\on{dim}_{F_1}B(C) \leqslant \on{dim}_{F_1}\underset{\epsilon \rightarrow 0}{\on{lim}}\, B_\epsilon(C) \leqslant \on{dim}_{F_1}\CA^{\mathrm{u}}_\chi.  
\end{equation*}
Using Corollary~\ref{a_chi_via_a_centr} and Proposition~\ref{main_Bethe}, we obtain $\on{dim}_{F_1}B(C)=\on{dim}_{F_1}\CA^{\mathrm{u}}_0(\mathfrak{z}_{\mathfrak{g}}(\chi))=\on{dim}_{F_1}\CA^{\mathrm{u}}_\chi$ and conclude that $\on{dim}_{F_1}\underset{\epsilon \rightarrow 0}{\on{lim}}\, B_\epsilon(C)=\on{dim}_{F_1}\CA^{\mathrm{u}}_\chi$ so from $\underset{\epsilon \rightarrow 0}{\on{lim}}\,B_\epsilon(C) \subset \CA^{\mathrm{u}}_\chi$ it follows that $\underset{\epsilon \rightarrow 0}{\on{lim}}\,B_\epsilon(C) = \CA^{\mathrm{u}}_\chi$.
\epr

\begin{Rem}
{\em{Note that in type $A$ Theorem \ref{lim_Bethe} can be deduced from the results of Section \ref{gen_A_u_chi_exp_resid}.}}
\end{Rem}

\section{Crystals: main properties and examples}\label{cryst_main}
From now on we assume that $\mathfrak{g}=\mathfrak{sl}_n$, $G=\on{PGL}_n$. Slightly abusing notation, we  denote by $\mathfrak{h}$ the subalgebra of diagonal matrices of $\mathfrak{gl}_n$ and by $\mathfrak{h}_0 \subset \mathfrak{h}$ we denote the subalgebra of traceless matrices. Let $T \subset \on{GL}_n$ be the subgroup of diagonal matrices and let $T_0 \subset T$ be the subgroup of matrices with determinant one. We also denote by $\ol{T} \subset \on{PGL}_n$ the subgroup of diagonal matrices.

In this section  we introduce notions of $\mathfrak{sl}_n$, $\hat{\mathfrak{sl}}_n$ - crystals and discuss some of their properties.   

\ssec{}{Definitions of $\mathfrak{sl}_n, \hat{\mathfrak{sl}}_n$-crystals}
Let us  discuss $\mathfrak{sl}_n$-crystals. 
We denote by $P^{\vee}$ the weight lattice of $\mathfrak{sl}_n$. 

\begin{defeni}{}
A $\mathfrak{sl}_n$-crystal is a finite set $B$ together with maps: 
\begin{equation*} 
\on{e}_{i},\, \on{f}_{i}: B \ra B \cup \{ 0\},\, \on{wt}\colon B \ra P^{\vee},\, i=1,2,\ldots n-1,
\end{equation*}
such that for each $i=1,2,\ldots,n-1$ we have:

$(a)$ let $b \in B$. If $\on{e}_{i} \cdot \, b \in B$ for some $i \in \{1,2,\ldots,n-1\}$, then 
$\on{wt}(\on{e}_{i} \cdot \, b) = \on{wt}(b) + \alpha_{i}^{\vee} 
,$
 if $\on{f}_{i} \cdot\,b \in B$ for some $i \in \{1,2,\ldots,n-1\}$, then
$\on{wt}(\on{f}_{i} \cdot \, b) = \on{wt}(b) - \alpha_{i}^{\vee} 
,$

$(b)$ for all $b, b' \in B$, \hspace{0,1cm}  $\on{e}_{i} \cdot \, b = b'$ if and only if $b = \on{f}_{i} \cdot \, b'$.
\end{defeni}


We use the identification $\{1,2,\ldots,n\} \iso \BZ/n\BZ,\, j \mapsto [j]$.
For $[j] \in \BZ/n\BZ$ we denote by $\tau_{[j]} \in S_n$ the cyclic  permutation which sends $[i]$ to $[i+j]$. We denote by $\tau_{[j]}^{\vee}$ the induced automorphism $\tau_{[j]}^{\vee}\colon P^\vee \iso P^\vee$.

Let us now define  ${\hat{\mathfrak{sl}}_n}$-crystals.
\begin{defeni}{}\label{sl_hat_cryst}
A ${\hat{\mathfrak{sl}_n}}$-crystal  is a finite set  ${\bf{B}}$ together with maps \begin{equation*}
\on{e}_{[i]},\, \on{f}_{[i]}\colon {\bf{B}} \ra {\bf{B}} \cup \{0\},\, \on{wt}\colon {\bf{B}} \ra P^\vee,\, [i] \in \BZ/n\BZ,    
\end{equation*}  
such that for every $[j] \in \BZ/n\BZ$ the operators \begin{equation*}
\on{e}_{i}:=\on{e}_{[i-j]},\, \on{f}_{i}:=\on{f}_{[i-j]},\,  {\bf{B}} \xrightarrow{\on{wt}} P^\vee \xrightarrow{\tau_{[j]}^\vee} P^\vee,\, i=1,2,\ldots,n-1
\end{equation*}
define the structure of $\mathfrak{sl}_n$-crystal on ${\bf{B}}$ that will be denoted by ${\bf{B}}^{[j]}$. 
\end{defeni}

\war{}
Note that our definition of $\hat{\mathfrak{sl}_n}$-crystal is not the standard one since we are assuming that the function $\on{wt}$ maps to $P^\vee$ that is a weight lattice of $\mathfrak{sl}_n$ (not of the whole algebra $\hat{\mathfrak{sl}_n}$), we are also considering only finite crystals.
\ewar

\ssec{}{Tensor products of  crystals} 
Let us discuss tensor products of  crystals. Given an 
$\hat{\mathfrak{sl}_n}$-crystal ${\bf{B}}$  and $[i] \in \BZ/n\BZ$ we can define maps $\varphi^+_{[i]},\, \varphi^-_{[i]}\colon {\bf{B}} \ra \BZ \cup \{\infty\}$ as follows: 
\begin{equation*}
\varphi^+_{{[i]}}(\bb) = \on{max} \Big\{m \in \mathbb{N}\,|\, (\on{e}_{[i]})^{m}(\bb) \neq 0 \Big\}, \hspace{0,1cm}  \varphi^-_{[i]}(\bb) = \on{max} \Big\{m \in \mathbb{N}\,|\, (\on{f}_{[i]})^{m}(\bb) \neq 0 \Big\}.
\end{equation*}

Given two integrable $\hat{\mathfrak{sl}_n}$-crystals $\bB,\, \bB'$ one can define an $\hat{\mathfrak{sl}_n}$-crystal structure on the set $\bB \times \bB'$ as follows:

$$\on{wt}(\bb \times \bb') = \on{wt}(\bb) + \on{wt}(\bb')$$

 \begin{equation} \on{e}_{[i]} \cdot (\bb \times \bb') =
\begin{cases}
(\on{e}_{[i]} \cdot \bb) \times \bb', \hspace{0,1cm} \on{if} \hspace{0,1cm}  \varphi^+_{[i]}(\bb) > \varphi^-_{[i]}(\bb')\\
\bb \times (\on{e}_{[i]} \cdot \bb'), \hspace{0,1cm}  \on{otherwise},
\end{cases}
\end{equation}

 \begin{equation} \on{f}_{[i]} \cdot (\bb \times \bb') =
\begin{cases}
(\on{f}_{[i]} \cdot \bb) \times \bb', \hspace{0,1cm}  \on{if} \hspace{0,1cm}  \varphi^+_{[i]}(\bb) \geqslant \varphi^-_{[i]}(\bb')\\
\bb \times (\on{f}_{[i]} \cdot \bb'), \hspace{0,1cm}  \on{otherwise}.
\end{cases}
\end{equation}



This crystal will be denoted by $\bB \otimes \bB'$. 


\ssec{}{Examples of $\mathfrak{sl}_n$-crystals}\label{ex_sl_n_cryst}
To every dominant weight $\la$ of $\mathfrak{sl}_n$ one can associate the
crystal $B_\la$ 
that should be considered as a ``discrete model" of the representation $V_\la$ (see~\cite{ka1},~\cite{ka2} for details). Explicit description of $B_\la$ was given in~\cite{kn2}, we recall the construction. Let $A(\la)$ be the Young diagram, corresponding to $\la$ i.e. if we decompose $\la=\sum_{i=1}^p \varpi_{m_i}$ with $1 \leqslant m_1 \leqslant \ldots \leqslant m_p \leqslant n$ then $A(\la)$ consists of $p$ columns of lengths $m_1,\ldots,m_p$. Let $B_\la$ be the set of semi-standard tableaux of shape $A(\la)$. Operators $\on{f}_i$ work as follows. Read the entries of the tableau from bottom to top left to right, ignoring all numbers except $i,i+1$. Replace $i+1$ by $($ and $i$ by $)$, turn $i$, corresponding to the rightmost unmatched $)$ by $i+1$. Similarly, operators $\on{e}_i$ work as follows. Read the entries of the tableau from bottom to top, left to right, ignoring all numbers except $i,i+1$. Replace $i+1$ by $($ and $i$ by $)$, turn $i+1$, corresponding to the leftmost unmatched $($ by $i$.

\begin{Ex}\label{cryst_on_fund_defe}
For $\la=\omega_l$ the crystal $B_{\omega_l}$ can be described as follows. Note that the Young diagram $A(\omega_l)$ is a column of length $l$. It follows that the set $B_{\omega_l}$ is in bijection with $l$-tuples of positive numbers $(i_1,\ldots,i_l)$ such that $1 \leqslant i_1 < \ldots < i_l \leqslant n$. Action of $\on{e}_i$ on $(i_1,\ldots,i_l)$ replaces $i$ by $i-1$, if possible,  and sends $(i_1,\ldots,i_l)$ to zero otherwise. Action of $\on{f}_i$ replaces  $i$ by $i+1$, if possible,  and sends $(i_1,\ldots,i_l)$ to  zero otherwise. The weight function $\on{wt}$ sends $(i_1,\ldots,i_l)$ to the linear function $\mathfrak{h}_0 \ni (x_1,\ldots,x_n) \mapsto x_{i_1}+\ldots+x_{i_l}$.
\end{Ex}

\rem{}
{\em{Let us give another description of the crystal structure on $B_\la$. Consider the following embedding $B_\la \subset B_{\varpi_{m_1}} \otimes \ldots \otimes B_{\varpi_{m_p}}$. Let ${{a}}$ be a semi-standard Young tableau of shape $A(\la)$. Let ${{a}}_i$ be the $i$-th column with enumeration starting from the right. We send ${{a}}$ to the element ${{a}}_1 \otimes {{a}}_2 \ldots \otimes {{a}}_p \in B_{\varpi_{m_1}} \otimes \ldots \otimes B_{\varpi_{m_p}}$.  The crystal structure on $B_\la$ is induced from the one on $B_{\varpi_{m_1}} \otimes \ldots \otimes B_{\varpi_{m_p}}$ (note that the crystal structure on $B_{\varpi_l}$ is very simple, see Example~\ref{cryst_on_fund_defe} above).}}
\erem


\begin{defeni}
A $\mathfrak{sl}_n$-crystal is normal if it is isomorphic to the disjoint union of crystals $B_\la$.
\end{defeni}

\ssec{}{Kirillov--Reshetikhin 
crystals}\label{kir_resh_cryst_def}

Let us describe the so-called Kirillov--Reshetikhin crystals. They are $\hat{\mathfrak{sl}_n}$-crystals that will be most important for us (see Remark~\ref{motiv_kr} for some motivation to restrict attention to these crystals). These $\hat{\mathfrak{sl}_n}$-crystals correspond to  representations 
$V_{l\varpi_r}$ and will be denoted ${\bf{B}}_{l\varpi_r}$. The existence of such crystals is proven in~\cite{kkmmnn} and the explicit construction is given in~\cite{shi},~\cite{kwo} and~\cite{k}, see also~\cite{oss}. 

In this section we prove that if $B_\la$ can be extended to an $\hat{\mathfrak{sl}}_n$-crystal ${\bf{B}}_\la$ such that ${\bf{B}}_\la^{[1]}$ is normal then $\la=l\varpi_r$ and such an extension is {\em{unique}}.  One of the results of this paper is a  ``geometric" construction of such an extension (see Section \ref{cryst_str_E_la}). So the present paper gives an alternative approach to the existence (and uniqueness) of  ${\bf{B}}_{l\varpi_r}$.

\rem{}\label{motiv_kr}
{\em{The (arguably) most important class of $\hat{\mathfrak{sl}}_n$-crystals are the crystals that appear as crystal graphs of (finite-dimensional) representations of $U_q(\hat{\mathfrak{sl}}_n)$. 
The following conjecture is due to Kashiwara (see~\cite[Conjecture~4.5]{shi2}): every connected 
affine crystal graph 
is isomorphic to the tensor product of Kirillov--Reshetikhin crystals. }}
\erem

We first describe the $\hat{\mathfrak{sl}}_n$-crystals, corresponding to representations $\La^l(\BC^n)=V_{\varpi_l},\, S^l(\BC^n)=V_{l\varpi_1}$.

\war{}
Note that it is not true in general that $S^l(\Lambda^r\BC^n)=V_{l\varpi_r}$. For example, for $n=4$, $l=r=2$, we have $S^2(\La^2\BC^4)=V_{2\varpi_2} \oplus \BC$.
\ewar


\begin{Ex} The crystal ${\bf{B}}_{\varpi_l}$ can be defined as follows. Pick the standard basis $\{v_{[1]},\ldots,v_{[n]}\}$ of $\BC^n$ and denote by $\{v_{[1]}^*,\ldots,v_{[n]}^*\} \subset (\BC^n)^*$ the dual basis. Consider the vector space $\La^l\BC^n$ and its projectivization $\BP(\La^l\BC^n)$. For a vector $x \in \La^l\BC^n$ we denote by $[x]$ the corresponding element of  $\BP(\La^l\BC^n)$. We define
\begin{equation*}
{\bf{B}}_{\varpi_l}=\{[v_{[i_1]}\wedge v_{[i_2]} \ldots \wedge v_{[i_l]}],\,|\,1 \leqslant i_1<i_2<\ldots <i_l \leqslant n\}.
\end{equation*}
The map $\on{e}_{[i]}$ acts via $E_{[i],[i+1]}$, the map $\on{f}_{[i]}$ acts via $E_{[i+1],[i]}$. The map $\on{wt}$ sends $[v_{[i_1]}\wedge v_{[i_2]} \ldots \wedge v_{[i_l]}]$ to the restriction of $v_{[i_1]}^*+\ldots +v_{[i_l]}^*$ to $\mathfrak{h}_0 \subset \BC^n$.
\end{Ex}

\begin{Ex} The crystal ${\bf{B}}_{l\varpi_1}$ can be described as follows. Consider the vector space $S^l\BC^n$ and its projectivization $\BP(S^l\BC^n)$, for a vector $x \in S^l\BC^n$ we denote by $[x]$ the corresponding element of  $\BP(S^l\BC^n)$.
We define 
\begin{equation*}
{\bf{B}}_{l\varpi_1}=\{[v_{[i_1]}^{p_1}\ldots v_{[i_k]}^{p_l}]\,|\,1 \leqslant i_1 < i_2 < \ldots < i_l \leqslant n,\, p_1+\ldots+p_l=l\}    
\end{equation*}
The map $\on{e}_{[i]}$ acts via $E_{[i],[i+1]}$, the map $\on{f}_{[i]}$ acts via $E_{[i+1],[i]}$. The map $\on{wt}$ sends $[v_{[i_1]}^{p_1}v_{[i_2]}^{p_2} \ldots v_{[i_l]}^{p_l}]$ to the restriction of $p_1v_{[i_1]}^*+\ldots+p_lv_{[i_l]}^*$ to $\mathfrak{h}_0 \subset \BC^n$.
\end{Ex}

\sssec{}{Sch\"utzenberger involutions and the operator $\phi$}
We start from the following definition.

\begin{defeni}{}\cite{hk}
The Sch\"utzenberger involution $\xi_{B_\la}\colon B_\la \to B_\la$ is the unique map of sets 
which satisfies 
\begin{equation*}
{\on{e}}_i(\xi_{B_\la}(b))=\xi_{B_\la}(\on{e}_{n-i}(b)),    
\end{equation*}
\begin{equation*}
{\on{f}}_i(\xi_{B_\la}(b))=\xi_{B_\la}(\on{f}_{n-i}(b)),     
\end{equation*}
\begin{equation*}
{\on{wt}}(\xi_{B_\la}(b))=w_0(\on{wt}(b)).    
\end{equation*}
Here $w_0 \in S_{n}$ is the longest element and $i \in \{1,\ldots,n-1\}$.
\end{defeni}
There is an explicit combinatorial construction of $\xi_{B_\la}$ (see, for example, \cite{llt}). If $B$ is a disjoint union of $B_\la$ then we denote by $\xi_B$ the corresponding involution of $B$.
Following \cite[Section 2.2]{hk} let $\ol{B}$ be the crystal with underlining set $\{\ol{b}\,|\, b \in B\}$ and crystal structure 
\begin{equation*}
\on{e}_i \cdot \ol{b} = \ol{\on{f}_{n-i} \cdot b},\, \on{f}_i \cdot b =\ol{\on{e}_{n-i} \cdot b},\, \on{wt}(\ol{b})= w_0(\on{wt}b).    
\end{equation*}
It follows from the definitions that the involution $\xi_{B}$ is the isomorphism of crystals $\xi_{B}\colon B \iso \ol{B}$. We set 
$
 B_\la^{(1)} := \ol{(\ol{B}_\la)|_{\mathfrak{sl}_{n-1}}}.
$
Let us describe the $\mathfrak{sl}_{n-1}$-crystal $ B_\la^{(1)}$. It follows from the definitions that   $B_\la^{(1)}=B_\la$ as a set i.e. $B_\la^{(1)}=\{b^{(1)}\,|\, b \in B_\la\}$ and the crystal structure is given by 
\begin{equation*}
\on{e}_i \cdot b^{(1)} = (\on{e}_{i+1} \cdot b)^{(1)},\, \on{f}_i \cdot b^{(1)} = (\on{f}_{i+1} \cdot b)^{(1)},\, \on{wt}(b^{(1)})=\tau_{[-1]}(\on{wt}(b)).    
\end{equation*}
Note that the composition  $\xi_{B_\la} \circ \xi_{(B_\la)|_{\mathfrak{sl}_{n-1}}}$ induces the isomorphism
\begin{equation*}
\xi_{B_\la} \circ \xi_{(B_\la)|_{\mathfrak{sl}_{n-1}}} \colon B_\la^{(1)} \iso (B_\la)|_{\mathfrak{sl}_{n-1}}
\end{equation*}
and such an isomorphism is {\em{unique}} (since $(B_\la)|_{\mathfrak{sl}_{n-1}}$ is isomorphic to the disjoint union of {\em{distinct}} irreducible $\mathfrak{sl}_{n-1}$-crystals).
We conclude that the following lemma holds.
\begin{lemma}{}\label{phi_vs_pr}
There exists the unique bijection  $\phi\colon B_\la \iso B_\la$  such that

$(1)$ $\on{wt}({\phi}(b))=\tau_{[1]}(\on{wt}(b))$


$(2)$ ${{\phi}} (\on{e}_i(b))=\on{e}_{i+1}({\phi}(b))$ and ${\phi} (\on{f}_i(b))=\on{f}_{i+1}({\phi}(b))$ for $i=1,\ldots,n-2$. 

The bijection $\phi$ is given by $\xi_{B_\la} \circ \xi_{(B_\la)|_{\mathfrak{sl}_{n-1}}}$, where $\xi_{(B_\la)|_{\mathfrak{sl}_{n-1}}}$ is the Sch\"utzenberger involution of the $\mathfrak{sl}_{n-1}$-crystal $(B_\la)|_{\mathfrak{sl}_{n-1}}$.
\end{lemma}


\sssec{}{Description of $\phi$ as a promotion operator} Let us now give a combinatorial description of the operator $\phi$ of Lemma \ref{phi_vs_pr}. Let us define the so-called Sch\"utzenberger’s promotion operator ${\bf{pr}}\colon B_{l\varpi_r} \iso B_{l\varpi_r}$. Recall that $A(l\varpi_r)$ is the Young diagram that corresponds to $l\varpi_r$. Starting from Young tableau $a \in B_{l\varpi_r}$ of shape $A(l\varpi_r)$, we denote by $H$ the horizontal strip in  $A(l\varpi_r)$ consisting of boxes of $a$ with letter $n$. Tableau ${\bf{pr}}(a)$ is constructed as follows.

Step $1$. We  remove all the letters $n$ in $a$.

Step $2$. Slide   the remaining subtableau to the southeast using the following procedure (Sch\"utzenberger’s jeu-de-taquin). 
Entering the cells of $H$
from left to right we
shift the entry above down, or the entry to the left right,
whichever is bigger (if it is same, take the vertical move) and continue this until there are
no entries above and to the left of the empty square.

Step $3$. We fill in the vacated cells with zeros.

Step $4$. Add one to each entry. The resulting tableau is ${\bf{pr}}(a)$.

It is easy to see that the operator ${\bf{pr}}(a)$ satisfies conditions $(1)$, $(2)$ of Lemma \ref{phi_vs_pr} so we must have 
\begin{equation*}
\phi=\xi_{B_\la} \circ \xi_{(B_\la)|_{\mathfrak{sl}_{n-1}}}={\bf{pr}}.   
\end{equation*}

\begin{Ex}
Let us consider the example $l=r=2$, $n=4$. We are dealing with the set $B_{2\varpi_2}$, consisting of Young tableaux of shape $\begin{Young}
&\cr
&\cr
\end{Young}$. Let us describe the operator ${\bf{pr}}=\phi$
in this particular case. 
The set $B_{2\varpi_2}$ consists of tableaux:
\begin{equation*}
\begin{Young}
1&1\cr
2&2\cr
\end{Young}, \begin{Young}
1&1\cr
2&3\cr
\end{Young}, \begin{Young}
1&1\cr
2&4\cr
\end{Young}, \begin{Young}
1&1\cr
3&3\cr
\end{Young}, \begin{Young}
1&1\cr
3&4\cr
\end{Young}, \begin{Young}
1&1\cr
4&4\cr
\end{Young}, \begin{Young}
1&2\cr
2&3\cr
\end{Young}, \begin{Young}
1&2\cr
2&4\cr
\end{Young}, \begin{Young}
1&2\cr
3&3\cr
\end{Young}, \begin{Young}
1&2\cr
3&4\cr
\end{Young},
\end{equation*}
\begin{equation*}
\begin{Young}
1&2\cr
4&4\cr
\end{Young}, 
\begin{Young}
1&3\cr
2&4\cr
\end{Young}, \begin{Young}
1&3\cr
3&4\cr
\end{Young}, \begin{Young}
1&3\cr
4&4\cr
\end{Young}, 
\begin{Young}
2&2\cr
3&3\cr
\end{Young}, \begin{Young}
2&2\cr
3&4\cr
\end{Young}, \begin{Young}
2&2\cr
4&4\cr
\end{Young}, \begin{Young}
2&3\cr
3&4\cr
\end{Young}, \begin{Young}
2&3\cr
4&4\cr
\end{Young}, \begin{Young}
3&3\cr
4&4\cr
\end{Young}.
\end{equation*}
Then ${\bf{pr}}$ acts on $B_{2\varpi_2}$ in the following way 
\begin{equation*}
\begin{Young}
1&1\cr
2&2\cr
\end{Young} \mapsto 
\begin{Young}
2&2\cr
3&3\cr
\end{Young} \mapsto 
\begin{Young}
3&3\cr
4&4\cr
\end{Young} \mapsto 
\begin{Young}
1&1\cr
4&4\cr
\end{Young} \mapsto 
\begin{Young}
1&1\cr
2&2\cr
\end{Young},~
 \begin{Young}
1&1\cr
2&3\cr
\end{Young} \mapsto 
\begin{Young}
2&2\cr
3&4\cr
\end{Young} \mapsto 
\begin{Young}
1&3\cr
3&4\cr
\end{Young} \mapsto 
\begin{Young}
1&2\cr
4&4\cr
\end{Young} \mapsto 
\begin{Young}
1&1\cr
2&3\cr
\end{Young},   
\end{equation*}
\begin{equation*}
\begin{Young}
1&1\cr
2&4\cr
\end{Young} \mapsto 
\begin{Young}
1&2\cr
2&3\cr
\end{Young} \mapsto 
\begin{Young}
2&3\cr
3&4\cr
\end{Young} \mapsto 
\begin{Young}
1&3\cr
4&4\cr
\end{Young} \mapsto 
\begin{Young}
1&1\cr
2&4\cr
\end{Young},~    
\begin{Young}
1&1\cr
3&4\cr
\end{Young} \mapsto 
\begin{Young}
1&2\cr
2&4\cr
\end{Young} \mapsto 
\begin{Young}
1&2\cr
3&3\cr
\end{Young} \mapsto 
\begin{Young}
2&3\cr
4&4\cr
\end{Young} \mapsto 
\begin{Young}
1&1\cr
3&4\cr
\end{Young},    
\end{equation*}
\begin{equation*}
\begin{Young}
1&1\cr
3&3\cr
\end{Young} \mapsto 
\begin{Young}
2&2\cr
4&4\cr
\end{Young} \mapsto 
\begin{Young}
1&1\cr
3&3\cr
\end{Young},~
\begin{Young}
1&2\cr
3&4\cr
\end{Young} \mapsto 
\begin{Young}
1&3\cr
2&4\cr
\end{Young} \mapsto 
\begin{Young}
1&2\cr
3&4\cr
\end{Young}.     
\end{equation*}

\end{Ex}

\sssec{}{Description of ${\bf{B}}_{l\varpi_r}$}
Pick a dominant $\la \in P^\vee$ and 
assume for a moment that there exists an $\hat{\mathfrak{sl}}_n$-crystal ${\mathbf{B}}_{\la}$ such that ${\mathbf{B}}_{\la}^{[0]} \simeq B_{\la}$ and ${\bf{B}}_\la^{[1]}$ is normal (we will see that this can happen only for $\la=l\varpi_r$ and the existence of such  crystal will be proven in Section \ref{cryst_str_E_la}). The goal of this section is to prove that ${\mathbf{B}}_{\la}$ is unique and can be described explicitly.

\prop{}\label{restr_cryst_class}
Let ${\mathbf{B}}_\la$ be a $\hat{\mathfrak{sl}}_n$-crystal such that ${\mathbf{B}}^{[0]}$ is isomorphic to $B_{\la}$ 
and ${\bf{B}}^{[1]}$ is normal. Then $\la=l\varpi_r$ and the operators $\on{e}_{[0]}$, $\on{f}_{[0]}$ are given by 
\begin{equation}\label{form_e_0}
\on{e}_{[0]}=(\xi_{B_\la} \circ \xi_{(B_\la)|_{\mathfrak{sl}_{n-1}}})^{-1} \circ \on{e}_1 \circ (\xi_{B_\la} \circ \xi_{(B_\la)|_{\mathfrak{sl}_{n-1}}})={\bf{pr}}^{-1} \circ {\on{e}}_1 \circ {\bf{pr}},
\end{equation}
\begin{equation}\label{form_f_0}
\on{f}_{[0]}=(\xi_{B_\la} \circ \xi_{(B_\la)|_{\mathfrak{sl}_{n-1}}})^{-1} \circ \on{f}_1 \circ (\xi_{B_\la} \circ \xi_{(B_\la)|_{\mathfrak{sl}_{n-1}}}) = {\bf{pr}}^{-1} \circ {\on{f}}_1 \circ {\bf{pr}}.  
\end{equation}
\eprop
\prf
Let us first of all note that since ${\bf{B}}^{[1]}$ is normal and has the same character as $B_\la$ we must have ${\bf{B}}^{[1]} \simeq B_\la$.
Pick isomorphisms of $\mathfrak{sl}_n$-crystals ${\bf{B}}^{[0]} \simeq B_\la$, ${\bf{B}}^{[1]} \simeq B_\la$ (note that these isomorphisms are unique). Composing one of them with the inverse of the other we obtain the 
bijection $\phi\colon B_\la \iso B_\la$ that satisfies conditions $(1)$, $(2)$ of Lemma \ref{phi_vs_pr}.
It then follows from Lemma \ref{phi_vs_pr} that $\phi=\xi_{B_\la} \circ \xi_{(B_\la)|_{\mathfrak{sl}_{n-1}}}$.
In particular, we see  that the operator $f_{[0]}$ is equal to ${\bf{pr}}^{-1} \circ \on{f}_1 \circ {\bf{pr}}$. Note  that ${\bf{pr}}^n\colon B_\la \iso B_\la$ should be the automorphism of the crystal $B_\la$ so we must have ${\bf{pr}}^n=\on{id}_{B_\la}$, hence, $\la=l\varpi_r$ (see~\cite{st} and \cite[Section~6.2]{kwo}). It now follows from the construction of ${\bf{B}}_{l\varpi_r}$ that ${\bf{B}} \simeq {\bf{B}}_{l\varpi_r}$ as $\hat{\mathfrak{sl}}_n$-crystals. 
\epr

\begin{remark}{}
In Proposition \ref{restr_cryst_class} we prove that if $B_{l\varpi_r}$ {\em{can be extended}} to an $\hat{\mathfrak{sl}}_n$-crystal then the extension is {\em{unique}} and operators ${\on{e}}_{[0]}$, ${\on{f}}_{[0]}$ {\em{must}} be given by the formulas  (\ref{form_e_0}), (\ref{form_f_0}) (compare with~\cite[Section 3]{shi}). It remains to check that the extension ${\bf{B}}_{l\varpi_r}$ indeed exists. This will be proven in Section \ref{cryst_str_E_la}. Let us recall that the existence of ${\bf{B}}_{l\varpi_r}$ is known and follows from~\cite{kkmmnn}. Our approach gives an alternative construction of  ${\bf{B}}_{l\varpi_r}$.
\end{remark}

\section{Yangians in type $A$}\label{yangian_type_A}
\ssec{}{Yangians of $\mathfrak{gl}_n,\, \mathfrak{sl}_n$, Bethe and Gaudin subalgebras}\label{type_A_Yang_Bethe}
Let us  return to Yangians.
We first of all describe the Yangian  $Y(\mathfrak{gl}_n)$. It is generated by $t_{ab}^{(r)},\, 1 \leqslant a,b \leqslant n,\, r \geqslant 1$ subject to relations
\begin{equation*}
[t_{ab}^{(r+1)},t_{cd}^{(s)}]-[t_{ab}^{(r)},t_{cd}^{(s+1)}]=t_{cb}^{(r)}t_{ad}^{(s)}-t_{cb}^{(s)}t_{ad}^{(r)},    
\end{equation*}
where $t^{(0)}_{a,b}:=\delta_{a,b}$. 
Yangian $Y(\fsl_n)$ is the subalgebra of $Y(\fgl_n)$, which is stable under all automorphisms of the form $T(u) \to f(u) T(u)$, where $f(u) \in 1 + u^{-1}\BC[[u^{-1}]]$. It follows from~\cite[Theorem~1.8.2]{m} that $Y(\mathfrak{gl}_n)=Y(\mathfrak{sl}_n) \otimes Z(Y(\mathfrak{gl}_n))$, where $Z(Y(\mathfrak{gl}_n))$ is the center of $Y(\mathfrak{gl}_n)$.

Let us also recall the description of Bethe subalgebras $\tilde{B}(C) \subset Y(\mathfrak{gl}_n)$. 
The symmetric group $S_n$ acts on $Y(\mathfrak{gl}_n)[[u^{-1}]] \otimes (\on{End}(\mathbb{C}^n))^{\otimes n}$ by permuting the tensor factors. This action factors through the embedding $S_n\hookrightarrow (\on{End}(\mathbb{C}^n))^{\otimes n}$, hence, the group algebra $\BC[S_n]$ is a subalgebra of $Y(\mathfrak{gl}_n)[[u^{-1}]] \otimes (\on{End}(\mathbb{C}^n))^{\otimes n}$. Pick $m \in \{1,2,\ldots,n\}$. Let $S_m$ be the subgroup of $S_n$ permuting the first $m$ tensor factors. Denote by $A_m$ the antisymmetrizer 
\begin{equation*}
A_m:=\frac{1}{m!}\sum_{\sigma \in S_m} (-1)^{\sigma}\sigma \in \mathbb{C}[S_m]\subset Y(\mathfrak{gl}_n)[[u^{-1}]] \otimes (\on{End}(\mathbb{C}^n))^{\otimes n}.
\end{equation*}
Note that for every $a \in \{1,\ldots,n\}$ there is an embedding
\begin{equation*}
i_a\colon Y(\mathfrak{gl}_n) \otimes \on{End}(\BC^n) \hookrightarrow Y(\mathfrak{gl}_{n}) \otimes (\on{End}(\BC^n))^{\otimes n}
\end{equation*}
which is identity on $Y(\mathfrak{gl}_n)$ and embeds $\on{End}(\BC^n)$ as the $a$-th tensor factor in $(\on{End}(\BC^n))^{\otimes n}$.
Suppose that $C \in \on{GL}_n$. 
For any $a \in \{1, \ldots ,n\}$ denote by $C_a$ the element $i_a(1\otimes C)\in Y(\mathfrak{gl}_n)[[u^{-1}]] \otimes (\operatorname{End}(\mathbb{C}^n))^{\otimes n}$. For any $a \in \{1,2,\ldots,n\}$ introduce the series with coefficients in $Y(\mathfrak{gl}_n)$ by 
\begin{equation}\label{def_tau_gl}
\tau_a(u,C) :=  \on{tr} A_a C_{1} \ldots C_{a} T_1(u) \ldots T_a(u-a+1),
\end{equation}
where we take the trace over all $a$ copies of $\on{End}(\mathbb{C}^n)$.
The algebra $\tilde{B}(C)$ is the subalgebra of $Y(\mathfrak{gl}_n)$ generated by all coefficients of the series $\tau_a(u,C)$.

\rem{}
{\em{The fact that this definition of $\tilde{B}(C)$ coincides with the one given in Section~\ref{bethe}
can be seen using that the irreducible representations of $\on{GL}_n$, corresponding to fundamental weights are wedge powers $\La^l(\BC^n)$.}}
\erem

For $C \in \on{SL}_n$ let
 $B(C) \subset Y(\mathfrak{sl}_n)$ be the corresponding Bethe subalgebra. One can show (see~\cite{imr}) that $\tilde{B}(C)=B(C) \otimes Z(Y(\mathfrak{gl}_n))$, where $Z(Y(\mathfrak{gl}_n))$ is the center of $Y(\mathfrak{gl}_n)$. 
It is easy to see that for $a \in \BC^\times$ we have $\tilde{B}(aC)=\tilde{B}(C)$ and similarly if $a^n=1$ then $B(aC)=B(C)$. It follows that actually $\tilde{B}(C),\, B(C)$ depend on $C \in \on{PGL}_n$.

Recall that $\mathfrak{h} \subset \mathfrak{gl}_n$ is the subalgebra of diagonal matrices.
Let us finally recall that to $\chi \in \mathfrak{h}$ one can associate the universal inhomogeneous Gaudin subalgebra of $U(\mathfrak{gl}_n[t])$ to be denoted $\tilde{\CA}^{\mathrm{u}}_\chi \subset U(\mathfrak{gl}_n[t])
$ (see Section~\ref{shift_univ_gaudin_def}). To the set of (distinct) points $u_1,\ldots,u_k \in \BC$ one can associate the inhomogeneous Gaudin subalgebra of $U(\mathfrak{gl}_n)^{\otimes k}$ to be denoted by $\tilde{\CA}_\chi(u_1,\ldots,u_k)$ (see Section~\ref{inhomogeneous_gaudin_via_conformal}).  Recall that $\mathfrak{h}_0 \subset \mathfrak{h}$ is the subalgebra of traceless matrices. We also consider the quotient $\mathfrak{h}/\BC$. Note now that the embedding $\mathfrak{h}_0 \subset \mathfrak{h}$ induces the identification $\mathfrak{h}_0 \iso \mathfrak{h}/\BC$ so to every $\chi \in \ol{\mathfrak{h}}$ we can associate the corresponding subalgebras $\CA^{\mathrm{u}}_\chi \subset U(\mathfrak{sl}_n[t])$,
$\CA_\chi(u_1,\ldots,u_k) \subset U(\mathfrak{sl}_n)^{\otimes k}$. We will sometimes abuse notations and denote $\mathfrak{h}/\mathbb{C}$ simply by $\mathfrak{h}$.

We will use the following lemma.
\lem{}\label{homog_comp_Gaud}
Pick $\chi \in \mathfrak{h}_0,\, s \in \BC^\times,\, c \in \BC$ and $z_1,\ldots,z_k \in \BC$. Then we have the equality $\CA_{s\chi}(z_1,\ldots,z_k)=\CA_{\chi}(sz_1+c,\ldots,sz_k+c)$. The same holds for the algebra $\tilde{\CA}_\chi(z_1,\ldots,z_k)$.
\elem
\prf
Follows from~\cite[Proposition~3]{r0} (see also~\cite[Lemma~9.2]{hkrw}). 
\epr

\ssec{}{The evaluation homomorphisms for $Y(\mathfrak{gl}_n),\, U(\mathfrak{gl}_n[t])$}\label{ev_hom_for_y_and_u_section}
Let
\begin{equation*}
E=(E_{ab})_{a,b=1,\ldots,n} \in \on{End}(\BC^n) \otimes U(\mathfrak{gl}_n)
\end{equation*}
be the matrix with coefficients in $U(\mathfrak{gl}_n)$, $E_{ab} \in \mathfrak{gl}_n$ is a matrix with $1$ on $(a,b)$-entry and zeroes on other entries. Recall that for every $z \in \BC$ we have the evaluation morphism
\begin{equation*}
 \on{\bf{ev}}_z\colon Y(\mathfrak{g}) \ra U(\mathfrak{g}),\, 
 t^{(r)}_{ab} \mapsto z^{r-1}E_{ab}.
\end{equation*}

%

Recall that we can identify $Y_0(\mathfrak{gl}_n)=\on{gr}_{2} Y(\mathfrak{gl}_n)$ with $U(\mathfrak{gl}_n[t])$ by sending $\on{gr}_{2} Y(\mathfrak{gl}_n) \ni [t_{ij}^{(r)}] \mapsto E_{ij}[r-1] \in U(\mathfrak{gl}_n[t])$.

\war{}
Note that this choice of the identification differs from the one in Proposition~\ref{iso_Y_Y} by the automorphism of $U(\mathfrak{gl}_n[t]) \iso U(\mathfrak{gl}_n[t])$ induced by the Cartan involution $\mathfrak{gl}_n \iso \mathfrak{gl}_n$ given by $x \mapsto -x^t$.
\ewar

Recall the identification $Y_{\epsilon}(\mathfrak{gl}_n) \iso Y(\mathfrak{gl}_n),\, [\hbar^m a] \mapsto \epsilon^m a$. The composition 
$Y_\epsilon(\mathfrak{gl}_n) \iso Y(\mathfrak{gl}_n) \xrightarrow{\on{\bf{ev}}_{z/\epsilon}} U(\mathfrak{gl}_n)$ is given by $[t_{ij}^{(r)}\hbar^{r-1}] \mapsto z^{r-1} E_{ij}$. We denote this composition by 
\begin{equation*}
\on{\bf{ev}}_{\epsilon;z}\colon Y_\epsilon(\mathfrak{gl}_n) \ra U(\mathfrak{gl}_n).  
\end{equation*}
One can show that $\underset{\epsilon \ra 0}{\on{lim}} \on{\bf{ev}}_{\epsilon;z}\colon U(\mathfrak{g}[t]) \ra U(\mathfrak{g})$ is the standard evaluation at $z$ homomorphism $\on{ev}_z\colon U(\mathfrak{gl}_n[t]) \ra U(\mathfrak{gl}_n)$ given by $x[l] \mapsto z^l x$ (see Lemma~\ref{compat_ev} and Corollary \ref{ext_ev_hbar} below).

More generally, we can define the evaluation homomorphism \begin{equation*}
\on{\bf{ev}}_{(z_1,\ldots,z_k)}\colon Y(\mathfrak{g}) \ra U(\mathfrak{g})^{\otimes k},\, \on{\bf{ev}}_{(z_1,\ldots,z_k)}:=(\on{\bf{ev}}_{z_1} \times \ldots \times \on{\bf{ev}}_{z_k}) \circ \Delta^k,     
\end{equation*}
where $\Delta^k\colon Y(\mathfrak{gl}_n) \ra Y(\mathfrak{gl}_n)^{\otimes k}$ is the composition $(\Delta \otimes \underbrace{1 \otimes \ldots \otimes 1}_{{k-1}}) \circ \ldots \circ  (\Delta \otimes 1) \circ \Delta $ and $\Delta\colon Y(\mathfrak{gl}_n) \ra Y(\mathfrak{gl}_n) \otimes Y(\mathfrak{gl}_n)$ is the standard comultiplication on $Y(\mathfrak{gl}_n)$ given by $\Delta(T)=T \otimes T$, where $T(u)=\Big(t_{ij}(u)\Big) \in \on{End}(\BC^n) \otimes Y(\mathfrak{g})$ and $t_{ij}(u)=\delta_{ij}+\sum_{r \geqslant 1}t_{ij}^{(r)}u^{-r}$.

Similarly, we obtain the map $\on{\bf{ev}}_{(\epsilon;z_1,\ldots,z_k)} \colon Y_\epsilon(\mathfrak{gl}_n) \ra U(\mathfrak{gl}_n)^{\otimes k}$ such that the limit $\underset{\epsilon \ra 0}{\on{lim}}\, \on{\bf{ev}}_{(\epsilon;z_1,\ldots,z_k)}$ is the standard evaluation at $z_1,\ldots,z_k$ homomorphism $\on{ev}_{z_1,\ldots,z_k} \colon U(\mathfrak{gl}_n[t]) \ra U(\mathfrak{gl}_n)^{\otimes k}$ given by $E_{ij}[r] \mapsto \sum_{a=1}^{k}E_{ij}[r]z_a^{r}$ (see Lemma~\ref{compat_ev}  and Corollary \ref{ext_ev_hbar} below). 

\rem{}
{\emph{Explicitly we have
\begin{equation*}
{\bf{ev}}_{z_1, \ldots, z_k}(T(u)) =  \left(1 + \frac{E^{(1)}}{u - z_1} \right) \ldots  \left(1 + \frac{E^{(k)}}{u - z_k} \right),
\end{equation*}
where $E^{(i)} = \on{id} \otimes \ldots \otimes \on{id} \otimes E \otimes \on{id} \otimes \ldots \otimes \on{id}$  ($E$ on the $i$-th place). It follows that 
\begin{equation}\label{ev_t_expl}
{\bf{ev}}_{z_1,\ldots,z_k}(t_{ij}^{(r)})=\sum_{l_{a_1}+\ldots+l_{a_{p}}+p=r}{E_{ij}^{(a_1)}}z_{a_1}^{l_{a_1}}\ldots {E_{ij}^{(a_p)}}z_{a_p}^{l_p}.   \end{equation}
}}
\erem

Pick distinct $z_1,\ldots,z_k \in \BC$ and $d_1,\ldots,d_k \in \BC$.
\lem{}\label{compat_ev}
Let $a(\epsilon) \in Y_{\epsilon}(\mathfrak{gl}_n),\,\epsilon \in \BC$ be an algebraic family and set $a:=a(0)$. Then $\underset{\epsilon \ra 0}{\on{lim}}\,
\on{\bf{ev}}_{(\epsilon;z_1+\epsilon d_1,\ldots,z_k+\epsilon d_k)}(a(\epsilon))$ exists 
and is equal to $\on{ev}_{\ul{z}}(a)$.
\elem
\prf
Recall that $Y_{\epsilon}(\mathfrak{gl}_n)$ is equal to the quotient $Y_{\hbar}(\mathfrak{gl}_n)/(\hbar-\epsilon)$. Let $\tilde{a} \in Y_{\hbar}(\mathfrak{gl}_n)$ be the element that corresponds to the family $a(\epsilon)$ i.e. $[\tilde{a}]_{\epsilon}=a(\epsilon)$, where $[\tilde{a}]_{\epsilon}$ is the class of $\tilde{a}$ in $Y_{\epsilon}(\mathfrak{gl}_n)$. Note that $a=[\tilde{a}]_0$. Recall now that $Y_{\hbar}(\mathfrak{gl}_n)$ is a free $\BC[\hbar]$-module.  $Y(\mathfrak{gl}_n)$ is generated by different products of elements $t^{(k)}_{mn}$ so $Y_{\hbar}(\mathfrak{gl}_n)$ as a module over $\BC[\hbar]$ is generated by the products of elements $t^{(k)}_{mn}\hbar^{k-1}$. 
We can then write $\tilde{a}$ as the linear combination of elements of the form $\hbar^{l+r_1+\ldots+r_p-p}t^{(r_1)}_{i_1j_1}\ldots t^{(r_p)}_{i_pj_p}$. Note that $[\hbar^{l+r_1+\ldots+r_p-p}t^{(r_1)}_{i_1j_1}\ldots t^{(r_p)}_{i_pj_p}]_{\epsilon}={\epsilon}^{l} [\hbar^{r_1+\ldots+r_p-p}t^{(r_1)}_{i_1j_1}\ldots t^{(r_p)}_{i_pj_p}]_{\epsilon}$ and the map $\on{\bf{ev}}_{(\epsilon;u_1+\epsilon d_1,\ldots,u_k+{\epsilon}d_k)}$ sends it to the element 
$
\epsilon^l{\bf{ev}}_{\epsilon,\ul{u}+\epsilon \ul{d}}([\hbar^{r_1-1}t_{i_1j_1}^{(r_1)}]_{\epsilon})\ldots {\bf{ev}}_{\epsilon,\ul{u}+\epsilon \ul{d}}([\hbar^{r_p-1}t_{i_pj_p}^{(r_p)}]_{\epsilon})
$
that goes to zero under the limit $\epsilon \ra 0$ if $l>0$.
So it remains to show that for every $[\hbar^{r-1}t_{ij}^{(r)}]_{\epsilon} \in Y_{\epsilon}(\mathfrak{gl}_n)$ 
$\underset{\epsilon \ra 0}{\on{lim}}\,
{\bf{ev}}_{\epsilon;\ul{u}+\epsilon\ul{d}}([\hbar^{r-1} t_{ij}^{(r)}]_{\epsilon})=\on{ev}_{\ul{u}}(E_{ij}[r-1])$.

Note that using~(\ref{ev_t_expl}) we have
\begin{equation*}
{\bf{ev}}_{\epsilon;\ul{z}+\epsilon\ul{d}}([\hbar^{r-1}t_{ij}^{(r)}]_\epsilon)=\epsilon^{r-1}{\on{ev}}_{\ul{z}/\epsilon+\ul{d}}(t_{ij}^{(r)})=\sum_{a=1}^kE_{ij}^{(a)}z_a^{r-1}+O(\epsilon).
\end{equation*}
This observation finishes the proof of the lemma since
$
\on{ev}_{\ul{z}}(E_{ij}[r-1])=\sum_{a=1}^k E_{ij}^{(a)}z_a^{r-1}.    
$
\epr

\cor{}\label{ext_ev_hbar}
There exists a homomorphism of $\BC[\hbar]$-algebras 
\begin{equation*}
{\bf{ev}}_{(\hbar;z_1+\hbar d_1,\ldots,z_k+\hbar d_k)}\colon Y_\hbar(\mathfrak{gl}_n) \ra U(\mathfrak{gl}_n)^{\otimes k}[\hbar],
\end{equation*}
whose fiber over $\hbar=\epsilon \in \BC^\times$ is ${\bf{ev}}_{(\epsilon;z_1+\epsilon d_1,\ldots,z_k+\epsilon d_k)}$ and the fiber over $\hbar=0$ is ${\on{ev}}_{\ul{z}}$. 
\ecor
\prf
Consider the natural embedding $Y_{\hbar}(\mathfrak{gl}_n) \subset Y_\hbar \otimes_{\BC[\hbar]} \BC[\hbar^{\pm 1}]$. Recall that we have the identification $Y_\hbar \otimes_{\BC[\hbar]} \BC[\hbar^{\pm 1}] \iso Y(\mathfrak{gl}_n)[\hbar^{\pm 1}]$ given by $(\hbar^ix) \otimes \hbar^l \mapsto \hbar^{l+i}x$. Consider the homomorphism ${\bf{ev}}_{\hbar^{-1}z_1+ d_1,\ldots,\hbar^{-1}z_k+d_k}\colon Y(\mathfrak{gl}_n)[\hbar^{\pm 1}] \ra U(\mathfrak{gl}_n)^{\otimes k}[\hbar^{\pm 1}]$. Composing this homomorphism with the embedding $Y_{\hbar}(\mathfrak{gl}_n) \subset Y(\mathfrak{gl}_n)[\hbar^{\pm 1}]$ we obtain the desired homomorphism (use Lemma \ref{compat_ev} to see that this homomorphism satisfies the required properties)
\begin{equation*}
{\bf{ev}}_{(\hbar;z_1+\hbar d_1,\ldots,z_k+\hbar d_k)}\colon Y_\hbar(\mathfrak{gl}_n) \ra U(\mathfrak{gl}_n)^{\otimes k}[\hbar].
\end{equation*}
 \epr

\section{Manin matrices and generators of $\tilde{\CA}^{\mathrm{u}}_\chi,\, \tilde{\CA}_\chi(z_1,\ldots,z_k)$}\label{gen_A_u_chi_exp_resid}

In this section following \cite{chemo} and \cite{tal} we introduce certain generators of the algebras $\tilde{\CA}^{\mathrm{u}}_\chi$, $\tilde{\CA}_\chi(z_1,\ldots,z_k)$, $\tilde{B}(C)$, ${\bf{ev}}_{z_1,\ldots,z_k}(\tilde{B}(C))$. 
We start from a general definition. Let $R$ be an algebra over complex numbers. Pick $M \in \on{End}(\BC^n) \otimes R$ that we can consider as $n \times n$ matrix $M=(M_{ij})$ with $M_{ij} \in R$. 

\defe{}
We say that $M$ is a Manin matrix if 
\begin{equation*}
\forall p,l,r,s=1,\ldots,n~\text{we have}~[M_{pl},M_{rs}]=[M_{rl},M_{ps}].
\end{equation*}
\edefe

Given a matrix $M \in \on{End}(\BC^n)$ we define its  column-determinant $\on{cdet}M  \in R$ by the following formula:
\begin{equation*}
\on{cdet}M:=\sum_{\sigma \in S_n} (-1)^\sigma \cdot  M_{\sigma(1)1}\ldots M_{\sigma(n)n}.	
\end{equation*}

For $i=1,\ldots,n$ we denote by $M_i \in \on{End}(\BC^n)^{\otimes n} \otimes R$ the image of $M$ under the embedding $\on{End}(\BC^n) \otimes R \hookrightarrow \on{End}(\BC^n)^{\otimes n} \otimes R$ given by 
\begin{equation*}
\on{End}(\BC^n) \otimes R \ni f \otimes x \mapsto 1 \otimes \ldots \otimes 1 \otimes \underset{i}{f} \otimes 1 \otimes \ldots \otimes 1	\otimes x \in \on{End}(\BC^n)^{\otimes n} \otimes R.
\end{equation*}

Recall that $A_n \in \BC[S_n] \subset \on{End}(\BC^n)^{\otimes n}$ is the antisymmetrizer normalized by the condition $A_n^2=A_n$. Explicitly, we have  $A_n=\frac{1}{n!}\sum_{\sigma \in S_n}(-1)^\sigma \sigma$.
\lem{}\label{prop_manin} 
Let $M \in \on{End}(\BC^n) \otimes R$ be a Manin matrix. Then 
\begin{equation*}
\on{tr}A_n M_1\ldots M_n=\on{cdet}M.
\end{equation*}
\elem
\prf
We have 
\begin{equation*}
\on{tr}A_n M_1\ldots M_n=\frac{1}{n!}\sum_{a,b \in S_n}(-1)^{ab}M_{a(1)b(1)}\ldots M_{a(n)b(n)}.
\end{equation*}
It follows from~\cite{chef} that for every $p \in S_n$ we have 
\begin{equation*}
\on{cdet}(M)=\sum_{\sigma \in S_n} \on{sgn}(\sigma) M_{\sigma(p(1))p(1)}\ldots M_{\sigma(p(n))p(n)}.
\end{equation*}
We conclude that 
\begin{equation*}
n!\on{cdet}M=\sum_{p,\sigma \in S_n}	\on{sgn}(\sigma) M_{\sigma(p(1))p(1)}\ldots M_{\sigma(p(n))p(n)}.
\end{equation*}
Replacing $\sigma \circ p$ by $a$ and $\sigma$ by $b$ we obtain that $\on{tr}A_nM_1\ldots M_n=\on{cdet}M$. 
\epr

Let us slightly modify generators of the algebra $\tilde{B}(C)$. Recall that the algebra $\tilde{B}(C)$ is generated by the elements 
\begin{equation*}
\tau_a(u,C) =  \on{tr} A_a C_{1} \ldots C_{a} T_1(u) \ldots T_a(u-a+1),\, a=1,\ldots,n.
\end{equation*}
 Recall the antipode map $\eta\colon Y(\mathfrak{gl}_n) \iso Y(\mathfrak{gl}_n)$ induced by the map $T(u) \mapsto T(u)^{-1}$. Then by~\cite[Section 1]{no} the algebra $\eta(\tilde{B}(C))=\tilde{B}(C^{-1})$ (see~\cite[Lemma~6.5]{ir0}) is generated by the elements 
\begin{equation*}
\on{tr}A_n T_1(u)\ldots T_k(u-k+1)C_{k+1}\ldots C_n.
\end{equation*}
It follows that the algebra $\tilde{B}(C)$ is generated by the elements
\begin{equation*}
\on{tr}A_nT_1(u)\ldots T_{a}(u-a+1)C_{a+1}^{-1}\ldots C_{n}^{-1},\, a=1,\ldots,n. 
\end{equation*}
These generators can be written in a generating function
\begin{equation*}
\on{tr}A_n(e^{-\partial_u}T_1(u)-C_1^{-1})\ldots (e^{-\partial_u}T_n(u)-C_n^{-1}).
\end{equation*}
Moreover, note that by~\cite[Proposition~4]{chef} matrix $e^{-\partial_u}T(u)-C^{-1}$ is Manin so we conclude from Lemma~\ref{prop_manin} that 
\begin{equation*}
\on{tr}A_n(e^{-\partial_u}T_1(u)-C_1^{-1})\ldots (e^{-\partial_u}T_n(u)-C_n^{-1})=\on{cdet}(e^{-\partial_u}T(u)-C^{-1}).
\end{equation*}

Recall the matrix $E \in \on{End}(\BC^n) \otimes U(\mathfrak{gl}_n)$.
Set $E(u):=\sum_{r \geqslant 0}E[r]u^{-r-1} \in \on{End}(\BC^n) \otimes U(\mathfrak{gl}_n[t])[[u^{-1}]]$.
\lem{}
We have 
\begin{equation*}
[E_{ij}(u),E_{kl}(v)]=\frac{1}{u-v}\Big(E_{il}(v)\delta_{jk}-E_{il}(u)\delta_{jk}-E_{kj}(v)\delta_{li}+E_{kj}(u)\delta_{li}\Big)	
\end{equation*}
i.e. $E(u)$ is a Lax matrix of $\mathfrak{gl}_n$-Gaudin type (see~\cite[Definition~5]{chef}).
\elem
\prf
Follows from the relation
$
[E_{ij},E_{kl}]=E_{il}\delta_{jk}-E_{kj}\delta_{li},
$.
\epr

\cor{}
The matrix $\partial_u-E(u) 
$ is a Manin matrix.
\ecor
\prf
Follows from~\cite[Proposition~2]{chef}, compare with~\cite[Equation~(18)]{chef}.
\epr

We are now ready to give an explicit description of the generators of the algebra $\tilde{\CA}_{\chi}^{\mathrm{u}}$. The following proposition is essentially due to Talalaev (\cite{tal}). 
\prop{}\label{expl_gen_A_u_chi}
The subalgebra $\tilde{\CA}^{\mathrm{u}}_\chi \subset U(\mathfrak{gl}_n[t])$ is generated by the coefficients in front of $u^r\partial_u^l$, $l \in \BZ_{\geqslant 0}$, $r \in \BZ_{\leqslant 0}$ of
$
\on{cdet}(E(u)-\partial_u+\chi).
$
\eprop
\prf
Recall that the generators of $\tilde{B}(C)$ are the coefficients of $\on{cdet}(e^{-\partial_u}T(u)-C^{-1})$. 
Consider the family $\on{cdet}\epsilon^{-1}(e^{-\epsilon \partial_u}T(u/\epsilon)-\on{exp}(-\epsilon \chi)) $ of elements of $\tilde{B}(\on{exp}(\epsilon \chi))$.
After the identification $\tilde{B}(\on{exp}(\epsilon\chi)) \simeq \tilde{B}_\epsilon(\on{exp}(\epsilon\chi))$ the element  $\on{cdet}\epsilon^{-1}(e^{-\epsilon\partial_u}T(u/\epsilon)-\on{exp}(-\epsilon\chi))$ becomes 
\begin{equation*}
\on{cdet}\epsilon^{-1}\Big(e^{-\epsilon\partial_u}\Big(1+\epsilon\sum_{r \geqslant 1}[\hbar^{r-1}t_{ij}^{(r)}]u^{-r}\Big)-\on{exp}(-\epsilon \chi)\Big) \in \tilde{B}_\epsilon(\on{exp}(\epsilon \chi)).
\end{equation*}
We have 
\begin{multline*}
e^{-\epsilon\partial_u}\Big(1+\epsilon\sum_{r \geqslant 1}[\hbar^{r-1}t_{ij}^{(r)}]u^{-r}\Big)-\on{exp}(-\epsilon \chi)=\\(1-\epsilon \partial_u)\Big(1+ \epsilon\sum_{r \geqslant 1}[\hbar^{r-1}t_{ij}^{(r)}]u^{-r}\Big)-1+\epsilon\chi+O(\epsilon^2)=\\=\epsilon\Big(\sum_{r \geqslant 1}[\hbar^{r-1}t_{ij}^{(r)}]u^{-r}-\partial_u+\chi\Big)+O(\epsilon^2).
\end{multline*}
It follows that 
\begin{equation*}
\underset{\epsilon \ra 0}{\on{lim}}\,\on{cdet}\epsilon^{-1}(e^{-\epsilon \partial_u}T(u/\epsilon)-\on{exp}(-\epsilon \chi))=\on{cdet}(E(u)-\partial_u+\chi).
\end{equation*}
It remains to recall that by Theorem~\ref{lim_Bethe} we have $\tilde{\CA}^{\mathrm{u}}_\chi=\underset{\epsilon \ra 0}{\on{lim}}\,\tilde{B}_\epsilon(\on{exp}(\epsilon \chi))$. The fact that the coefficients of $\on{cdet}(E(u)-\partial_u+\chi)$ generate $\tilde{\CA}^{\mathrm{u}}_\chi$ follow from the similar statement for the classical limit $\tilde{\ol{\CA}}{}^{\mathrm{u}}_\chi \subset S^\bullet(\mathfrak{gl}_n[t])$ together with the fact that $\on{gr}_{PBW}\tilde{\CA}^{\mathrm{u}}_\chi = \on{gr}_{PBW} \tilde{\ol{\CA}}{}^{\mathrm{u}}_\chi$. 
\epr

\rem{}
{\em{For $\chi=0$ the results of Proposition~\ref{expl_gen_A_u_chi} should be compared with~\cite[Theorem~3.1]{chemo}.}} 
\erem

Let us now pick distinct points $z_1,\ldots,z_k \in \BC$ and recall the corresponding standard Lax matrix for the Gaudin system
\begin{equation*}
L_{\ul{z}}(u):=\sum_{i=1}^k \frac{E^{(i)}}{u-z_i}.	
\end{equation*}

\cor{}\label{expl_gen_A_chi_z}
For every $\chi \in \mathfrak{gl}_n$ the algebra $\tilde{\CA}_\chi(z_1,\ldots,z_k)$ is generated by the coefficients in front of $u^r\partial_u^l$ of $\on{cdet}(L(u)-\partial_u-\chi)$. 
\ecor
\prf
Recall that $\CA_\chi(z_1,\ldots,z_k)=\on{ev}_{z_1,\ldots,z_k}(\CA^{\mathrm{u}}_{-\chi})$. 
Now the claim follows from the fact that 
\begin{equation*}
\on{ev}_{z_1,\ldots,z_k}(E(u))=\sum_{i=1}^k\sum_{r=0}^{\infty}Ez_i^{r}u^{-r-1}=\sum_{i=1}^k \frac{E^{(i)}}{u-z_i} =L_{\ul{z}}(u).	
\end{equation*}

\epr


Pick $\xi \in \mathfrak{h}$.
Our goal for now is to introduce certain compatible families of generators of the algebras ${\bf{ev}}_{\ul{z}/\varepsilon+\ul{d}}\tilde{B}(\on{exp}(\epsilon\chi))$, $\tilde{\CA}_{\chi}(\ul{z})$.

Consider the following two-parameter family of 
differential operators, here $\epsilon,\varepsilon \in \BC^\times$ 
\begin{equation*}
\on{cdet}\epsilon^{-1}(e^{-\epsilon\partial_u}(1+L_{\ul{z}/\varepsilon+\ul{d}}(u/\epsilon))-\on{exp}(-\epsilon\chi))=\sum_{k \geqslant 0}b_{k,\epsilon,\varepsilon}(u)\partial_u^k.
\end{equation*}
Since $e^{-\partial_u}u^r=u^{r-1}e^{-\partial_u}$ for every $r \in \BZ$ it follows that $b_{k,\epsilon,\varepsilon}(u)$ are rational functions with poles at the points $\frac{\epsilon}{\varepsilon}z_i+\epsilon j+\epsilon d_i$, $i,j=1,\ldots,k$. 
Consider then the following elements
\begin{equation*}
r_{k,z_i,l,\epsilon,\varepsilon}:=\sum_{j=1}^k\on{res}_{u=\frac{\epsilon}{\varepsilon}z_i+\epsilon j+\epsilon d_i}(u-z_i)^lb_{k,\epsilon,\varepsilon}(u)du \in {\bf{ev}}_{\ul{z}/\varepsilon+\ul{d}}(\tilde{B}(\on{exp}(\epsilon \chi))),\, l=0,\ldots,k.
\end{equation*}

We can also decompose 
\begin{equation*}
\on{cdet}(L_{\ul{z}}(u)-\partial_u-\xi)=\sum_{k \geqslant 0}b^0_{k,\ul{z},\xi}(u)\partial_u^k
\end{equation*}
and consider the elements
\begin{equation*}
r_{k,z_i,l,\xi}^0:=\on{res}_{u=z_i}(u-z_i)^lb^0_{k,\ul{z},\xi}(u)du,\, l=0,\ldots,k.
\end{equation*}

The following proposition is an immediate consequence of the results above.
\prop{}\label{gen_A_chi_z_via_res}
The algebra ${\bf{ev}}_{\ul{z}/\varepsilon+\ul{d}}(\tilde{B}(\on{exp}(\epsilon \chi)))$ is generated by elements $r_{k,z_i,l,\epsilon,\varepsilon}$. The algebra $\tilde{\CA}_{\xi}(\ul{z})$ is generated by the elements $r_{k,z_i,l,\xi}^0$
\eprop

\rem{}\label{shift_inv_ev_Bethe}
{\em{Note that for every $c \in \BC$ the elements  $r_{k,z_i,l,\epsilon,\varepsilon}$, $r^0_{k,z_i,l,\xi}$ do not change under the simultaneous shift $(z_1,\ldots,z_k) \mapsto (z_1+c,\ldots,z_n+c)$ so in particular 
\begin{equation*}
{\bf{ev}}_{((z_1/\varepsilon)+d_1,\ldots,(z_k/\varepsilon)+d_k)}(\tilde{B}(\on{exp}(\epsilon \chi)))=
{\bf{ev}}_{((z_1+c)/\varepsilon+d_1,\ldots,(z_k+c)/\varepsilon+d_k)}(\tilde{B}(\on{exp}(\epsilon \chi))),
\end{equation*}
\begin{equation*}
\tilde{\CA}_\xi(z_1,\ldots,z_n)=\tilde{\CA}_\xi(z_1+c,\ldots,z_k+c).
\end{equation*} 
}}
\erem

\section{Alcoves}\label{alcoves} 
Let $S \subset T$ be the compact torus i.e. $S:=T \cap U(n)$, where $U(n) \subset \on{GL}_n$ is the group of unitary matrices. Abusing notations, we denote by the same symbol $S$ the quotient $(T \cap U(n))/U(1) \subset \on{PGL}_n$, here $U(1) \subset U(n)$ is the center. Recall the Cartan subalgebra $\mathfrak{h} \subset \mathfrak{gl}_n$ consisting of diagonal matrices. We use the natural bijection $\{1,2,\ldots,n\} \iso \BZ/n\BZ,\, j \mapsto [j]$ to parametrize coordinates of $\mathfrak{h}$ by elements of $\BZ/n\BZ$. 
Again abusing notations we will also denote by $\mathfrak{h}$ the quotient of the diagonal matrices in $\mathfrak{gl}_n$ by $\mathbb{C}$.
We have a covering:  
\begin{equation*}
\on{exp}\colon \mathfrak{h}_{\BR} \ra S,\, \chi \mapsto \on{exp}(2\pi i \chi).
\end{equation*}
We denote by $S^{\mathrm{reg}} \subset S$ the subset of regular elements. Note that the preimage $\,\on{exp}^{-1}(S \setminus S^{\mathrm{reg}})$ consists of the points $(a_{[1]},\ldots,a_{[n]}) \in \mathfrak{h}_{\BR}$ such that $a_{[i]}-a_{[j]} \in \BZ$ for some \([i] \neq [j]\). For $[i] \neq [j]$ and $k \in \BZ$ we set
\begin{equation*}
H_{[i],[j]}^k=\{(a_{[1]},\ldots,a_{[n]}) \in \mathfrak{h}_{\BR}\,|\, a_{[i]}-a_{[j]}=k\}. \end{equation*}
We will call the hyperplanes $H_{[i],[j]}^k \subset \mathfrak{h}_{\BR}$ the {\em{affine walls}}, and the connected components of the complement to the union of all $H_{[i],[j]}^k$ the \emph{alcoves}. Let $\La$ be the coweight lattice of $\mathfrak{gl}_n$. We have the natural identification $\La=\BZ^n$. Abusing notations, we will denote by the same symbol $\Lambda$ the quotient $\La/\BZ(1,1,\ldots,1)$.
Consider the standard action $S_n \curvearrowright \La$ via the permutation of the basis elements. 
Note that we have the left action of $S_{n} \ltimes \La $ on $\mathfrak{h}$ given by 
\begin{equation*}
(\sigma;m_1,\ldots,m_n) \cdot (a_1,\ldots,a_n)=(a_{\sigma^{-1}(1)}+m_{\sigma^{-1}(1)},\ldots,a_{\sigma^{-1}(n)}+m_{\sigma^{-1}(n)}).
\end{equation*}
\rem{}
{\em{Recall that we have $(\sigma_1,\la_1)\cdot (\sigma_2,\la_2)=(\sigma_1\sigma_2,\sigma_2^{-1}(\la_1)+\la_2)$ for $\sigma_1,\,\sigma_2 \in S_n,\, \la_1,\,\la_2 \in \La$.}}
\erem
We define the \emph{extended} affine Weyl group as $\widehat{W}^{\mathrm{ext}}:=S_n \ltimes \La$. The natural action of $S_n \curvearrowright \mathfrak{h}$ and the action of $\La \curvearrowright \mathfrak{h}$ via translations induce the action  $\widehat{W}^{\mathrm{ext}} \curvearrowright \mathfrak{h}$. 

The coroot lattice $\La_r \subset \La$ is generated by the vectors $E_{ii}-E_{i+1i+1} \in \mathfrak{h}$, $i=1,\ldots,n-1$. 
Note that $[\La:\La_r]=n$. The affine Weyl group $\widehat{W}:=S_n \ltimes \La_r$ is an index $n$ subgroup of the extended affine Weyl group $\widehat{W}^{\mathrm{ext}}$.

The $\widehat{W}^{\mathrm{ext}}$-action on $\mathfrak{h}_{\BR}$ preserves the affine walls, so it acts on the set alcoves. It is easy to see that the alcoves
are fundamental domains for the action $\widehat{W} \curvearrowright \mathfrak{h}_{\BR}$. By choosing the ``base'' alcove 
\begin{equation*}
Q=\{(a_{[1]},\ldots,a_{[n]}) \in \mathfrak{h}_{\BR},\,|\, a_{[n]}+1 \geqslant a_{[1]} \geqslant a_{[2]} \geqslant \ldots \geqslant a_{[n]}\},
\end{equation*}
we obtain a bijection between the alcoves and the elements of $\widehat{W}$ and a surjective $n:1$ map from $\widehat{W}^{\mathrm{ext}}$ to the set of alcoves. More precisely,  
to $\hat{w}=(\sigma,m_{[1]},\ldots,m_{[n]}) \in \widehat{W}^{\mathrm{ext}}$ we associate 
\begin{multline}\label{defeni_of_alcove_corr_w}
Q_{\hat{w}}:=\hat{w}Q=\\
=\{(a_{[1]},\ldots,a_{[n]}) \in \mathfrak{h}_{\BR}\,|\, a_{\sigma([n])}-m_{[n]}+1 \geqslant a_{\sigma([1])}-m_{[1]} \geqslant \ldots \geqslant a_{\sigma([n])}-m_{
[n]}\}.    
\end{multline}

\war{}
Let us point out once again that the $\widehat{W}^{\mathrm{ext}}$-action on the set of alcoves is not free, so it may happen that $Q_{\hat{w}}=Q_{\hat{w}'}$ for distinct $\hat{w},\, \hat{w}' \in \widehat{W}^{\mathrm{ext}}$. This will be important in the proof of Proposition~\ref{dif_fin_res}.
\ewar

The base alcove $Q$ is separated by the walls 
\begin{equation*}
H^{0}_{[1],[2]},\ldots,H_{[n-1],[n]}^{0},H^{-1}_{[n],[1]}    
\end{equation*}
that we denote by $H_1,H_2,\ldots,H_n$ respectively.
We denote by $Q^{\mathrm{reg}}_{\hat{w}}$ the 
interior of $Q_{\hat{w}}$. The alcove $Q_{\hat{w}}$ is separated by the walls 
\begin{equation*}
H_{\sigma([1]),\sigma([2])}^{m_{[1]}-m_{[2]}},\ldots,H_{\sigma([n-1]),\sigma([n])}^{m_{[n-1]}-m_{[n]}}, H_{\sigma([n]),\sigma([1])}^{m_{[n]}-m_{[1]}-1}
\end{equation*}
that we enumerate by numbers $1,2,\ldots,n$ and denote by $H^{\hat{w}}_1,\ldots,H^{\hat{w}}_n$. Note that $H^{\hat{w}}_i=\hat{w}(H_i)$. 

\section{Limits to the affine wall}\label{limits_to_the_wall}
In this section we describe the limits of different families of subalgebras depending on $\chi \in \mathfrak{h}$ to a generic point of an (affine) wall $H^k_{[i],[j]}$. 

We pick $\hat{w} \in \widehat{W}^{\mathrm{ext}}$ and consider the corresponding alcove $Q_{\hat{w}}$ (see~(\ref{defeni_of_alcove_corr_w})).
Pick $\chi_0 \in Q_{\hat{w}}$ lying on some $H_{[i],[j]}^{k}$ but not on any other affine wall (such $\chi_0 \in H_{[i],[j]}^k$ will be called \emph{subregular}). Set $C_0:=\on{exp}(2\pi i \chi_0)$. Pick also $\chi \in Q^{\mathrm{reg}}_{\hat{w}}$ and consider the following family $C(\varepsilon):=C_0\on{exp}(2\pi i\varepsilon \chi)=\on{exp}(2\pi i (\chi_0+\varepsilon \chi))$, where $\varepsilon \in \BC$ (one can also consider the formal version of this family).
Set $B(C_0,\chi):=\underset{\varepsilon \ra 0}{\on{lim}} B(C(\varepsilon))$. Recall that we have the identification $Y(\mathfrak{g}) \simeq Y_{\epsilon}(\mathfrak{g})$. We denote by $B_\epsilon(C_0,\chi)$ the image of $B(C_0,\chi)$ in $Y_{\epsilon}(\mathfrak{g})$. We also set $\CA^{\mathrm{u}}_{(\chi_0,\chi)}:=\underset{\varepsilon \ra 0}{\on{lim}}\,\CA^{\mathrm{u}}_{\chi_0+\varepsilon \chi}$.

\rem{}
{\em{The algebra $B(C_0,\chi)$ is the limit subalgebra for the family $\{B(C),\, C \in T^{\mathrm{reg}}\}$. It follows from~\cite{ir0} that the space of all possible limit subalgebras is classified by the Deligne-Mumford space $\ol{M_{0,n+2}}$ of stable rational curves with $n+2$ marked points.}}
\erem

We will use the following lemma in the proof of Proposition~\ref{lim_B} below.
\lem{}\label{cart_Bethe}
Pick $C \in T^{\mathrm{reg}}$ then $\mathfrak{h}_0 \subset B(C)$ (we use the embedding $U(\mathfrak{g}) \subset Y(\mathfrak{g})$). 
\elem
\prf
Recall the filtration $F_2$ on $Y(\mathfrak{g})$ (see Section~\ref{two_filtr_yang}).
It follows from~\cite{ir2} that $\on{gr}_{2} B(C)=U(\mathfrak{h}_0[t]) \subset U(\mathfrak{sl}_n[t])=\on{gr}_{2} Y(\mathfrak{g})$ so $\mathfrak{h}_0 \subset \on{gr}_2 B(C)$. Note that $\on{deg}_{F_2}\mathfrak{h}_0=0$ and the only element of degree $-1$ in $Y(\mathfrak{g})$ is zero.
The claim follows. 
\epr

The following proposition should be compared  with~\cite[Theorem~B]{ir2}.
\prop{}\label{lim_B}
We have 
\begin{equation*}
B(C_0,\chi)=B(C_0) \otimes_{Z(U(\mathfrak{z}_{\mathfrak{g}}(\chi_0)))} \CA_{\chi}(\mathfrak{z}_{\mathfrak{g}}(\chi_0)),\, \CA^{\mathrm{u}}_{(\chi_0,\chi)}=\CA^{\mathrm{u}}_{\chi_0} \otimes_{Z(U(\mathfrak{z}_{\mathfrak{g}}(\chi_0)))} \CA_{\chi}(\mathfrak{z}_{\mathfrak{g}}(\chi_0)).
\end{equation*}
Moreover, the subalgebra  $B(C_0,\chi) \subset Y(\mathfrak{g})$ is generated by $B(C_0)$ and the element $h=h_{ij}=E_{ii}-E_{jj} \in U(\mathfrak{g}) \subset Y(\mathfrak{g})$, the subalgebra $\CA^{\mathrm{u}}_{(\chi_0,\chi)} \subset U(\mathfrak{g}[t])$ is generated by $\CA^{\mathrm{u}}_{\chi_0}$ and the element $h_{ij} \in U(\mathfrak{g}) \subset U(\mathfrak{g}[t])$. 
\eprop
\prf
We have $B(C_0) \subset \underset{\epsilon \ra 0}{\on{lim}} B(C(\epsilon))=B(C_0,\chi)$. We claim that $h_{ij}$ lies in the limit. Indeed, this follows from the fact that  $C=C(\epsilon)$ is regular for $\epsilon \neq 0$ so by Lemma~\ref{cart_Bethe} we have $\mathfrak{h}_0 \subset B(C(\epsilon))$.
Note now that the derived subalgebra  $\mathfrak{z}_{\mathfrak{g}}(\chi_0)^{\mathrm{der}}$ of 
$\mathfrak{z}_{\mathfrak{g}}(\chi_0)$
 is equal to 
$\on{Span}_{\BC}(E_{ij},h_{ij},E_{ji})\simeq \mathfrak{sl}_2$ so the algebra $\CA_{\chi}(\mathfrak{z}_{\mathfrak{g}}(\chi_0))$ is equal to $\BC[h_{ij}] \otimes Z(U(\mathfrak{z}_{\mathfrak{g}}(\chi_0)))$ i.e. is (freely) generated by $h_{ij}$ over the center $Z(U(\mathfrak{z}_{\mathfrak{g}}(\chi_0)))$.
We conclude that 
$B(C_0) \otimes_{Z(U(\mathfrak{z}_{\mathfrak{g}}(\chi_0)))} \CA_{\chi}(\mathfrak{z}_{\mathfrak{g}}(\chi_0)) \subset 
B(C_0,\chi)$.
Note that $Z(U(\mathfrak{z}_{\mathfrak{g}}(\chi_0))) \subset B(C_0) \subset B(C_0,\chi)$ so $B(C_0,\chi)$ is contained in the centralizer $Z_{Y(\mathfrak{g})}(Z(U(\mathfrak{z}_{\mathfrak{g}}(\chi_0))))$. It follows from~\cite[Lemma~6.14]{ir2} that   
\begin{equation*}
Z_{Y(\mathfrak{g})}(Z(U(\mathfrak{z}_{\mathfrak{g}}(\chi_0))))=Y(\mathfrak{g})^{\mathfrak{z}_{\mathfrak{g}}(\chi_0)} \otimes_{Z(U(\mathfrak{z}_{\mathfrak{g}}(\chi_0)))} U(\mathfrak{z}_{\mathfrak{g}}(\chi_0)).
\end{equation*}
Algebra $B(C_0) \otimes_{Z(U(\mathfrak{z}_{\mathfrak{g}}(\chi_0)))} \CA_{\chi}(\mathfrak{z}_{\mathfrak{g}}(\chi_0))$ is a maximal commutative subalgebra of $Y(\mathfrak{g})^{\mathfrak{z}_{\mathfrak{g}}(\chi_0)} \otimes_{Z(U(\mathfrak{z}_{\mathfrak{g}}(\chi_0)))} U(\mathfrak{z}_{\mathfrak{g}}(\chi_0))$ (compare with the proof of Proposition~\ref{A_chi_max}) so we must have the equality $B(C_0) \otimes_{Z(U(\mathfrak{z}_{\mathfrak{g}}(\chi_0)))} \CA_{\chi}(\mathfrak{z}_{\mathfrak{g}}(\chi_0)) = B(C_0,\chi)$. 

The proof for $\CA^{\mathrm{u}}_{(\chi_0,\chi)}$ is similar.
\epr

Recall that $\CA^{\mathrm{u}}_{(\chi_0,\chi)}=\CA^{\mathrm{u}}_{\chi_0} \otimes_{Z(U(\mathfrak{z}_{\mathfrak{g}}(\chi_0)))} \CA_{\chi}(\mathfrak{z}_{\mathfrak{g}}(\chi_0))$. We set  
\begin{equation*}
\ol{\CA}{}^{\mathrm{u}}_{(\chi_0,\chi)}:=\ol{\CA}{}^{\mathrm{u}}_{\chi_0} \otimes_{Z(S^\bullet(\mathfrak{z}_{\mathfrak{g}}(\chi_0)))} \ol{\CA}_{\chi}(\mathfrak{z}_{\mathfrak{g}}(\chi_0))
\end{equation*}
that can be considered as a classical limit of  $\CA^{\mathrm{u}}_{(\chi_0,\chi)}$. 

\rem{}
{\em{In the same way as in the proof of Proposition~\ref{lim_B} one can show that $\ol{\CA}{}^{\mathrm{u}}_{(\chi_0,\chi)}=\underset{\epsilon \ra 0}{\on{lim}}\,\ol{\CA}{}^{\mathrm{u}}_{\chi_0+\epsilon\chi}$ and $\ol{\CA}{}^{\mathrm{u}}_{(\chi_0,\chi)}$ is generated by $\ol{\CA}{}^{\mathrm{u}}_{\chi_0}$ and $h$.
}} 
\erem

\lem{}\label{central_chi_chi0}
We have 
\begin{equation*}
\ol{\CA}{}^{\mathrm{u}}_{(\chi_0,\chi)}=Z_{S^{\bullet}(\mathfrak{g}[t])}(\Omega_{\chi_0},h),\, 
\CA^{\mathrm{u}}_{(\chi_0,\chi)}=Z_{U(\mathfrak{g}[t])}(\tilde{\Omega}_{\chi_0},h).  
\end{equation*}
\elem
\prf
The proof is  similar to the proof of Proposition~\ref{class_gaud_centr}. 
Considering the family $Z_{S^{\bullet}(\mathfrak{g}[t])}(\Omega_{\kappa\chi_0},h)$ and taking the limit $\kappa \ra \infty$ we conclude (as in the proof of Proposition~\ref{class_gaud_centr}) that 
\begin{equation*}
\underset{\kappa \ra \infty}{\on{lim}}Z_{S^\bullet(\mathfrak{g}[t])}(\Omega_{\kappa\chi_0},h) \subset Z_{S^\bullet(\mathfrak{z}_{\mathfrak{g}}(\chi_0)[t])}(\Omega_{\mathfrak{z}_{\mathfrak{g}}(\chi_0)},h).
\end{equation*}
It remains to show that $\ol{\CA}{}^{\mathrm{u}}_{0}(\mathfrak{z}_{\mathfrak{g}(\chi_0)}) \cdot h=Z_{S^\bullet(\mathfrak{z}_{\mathfrak{g}}(\chi_0)[t])}(\Omega_{\mathfrak{z}_{\mathfrak{g}}(\chi_0)},h)$. Note that   $\mathfrak{z}_{\mathfrak{g}}(\chi_0)$ is a reductive Lie algebra with semisimple part being isomorphic to $\mathfrak{sl}_2$. Now the claim reduces to the $\mathfrak{sl}_2$-case when this is clear. 

So we have shown that $\ol{\CA}{}^{\mathrm{u}}_{(\chi_0,\chi)}=Z_{S^{\bullet}(\mathfrak{g}[t])}(\Omega_{\chi_0},h)$. The proof of the equality  $\CA^{\mathrm{u}}_{(\chi_0,\chi)}=Z_{U(\mathfrak{g}[t])}(\tilde{\Omega}_{\chi_0},h)$ follows from the equality $\ol{\CA}^{\mathrm{u}}_{(\chi_0,\chi)}=Z_{S^\bullet(\mathfrak{g}[t])}({\Omega}_{\chi_0},h)$ in the same way as in the proof of Proposition~\ref{ca_as_centr}.
\epr

\lem{}\label{extend_C2}
The two-parameter family $B_{\epsilon}(\on{exp}(\epsilon(\chi_0+\varepsilon \chi)))$ depending on $(\epsilon,\varepsilon) \in (\BC^\times)^2$ extends to the continuous family on the whole $\BC^2$ by setting 
\begin{equation*}
(\epsilon,0) \mapsto B_{\epsilon}(\on{exp}(\epsilon \chi_0),\chi),\, (0,\varepsilon) \mapsto \CA^{\mathrm{u}}_{\chi_0+\varepsilon\chi},\, (0,0) \mapsto \CA^{\mathrm{u}}_{(\chi_0,\chi)}.
\end{equation*}
\elem
\prf
As in the proof of Proposition~\ref{W_lim} we consider the family 
\begin{equation*}
X_{\epsilon}(\chi_0+\varepsilon \chi) := \frac{\epsilon}{c_V}\psi_\epsilon\Big(\on{tr}_V \on{exp}\rho_V(\epsilon(\chi_0+\varepsilon \chi))T^{(3)}\big)
 \in B_{\epsilon}(\on{exp}(\epsilon(\chi_0+\varepsilon \chi))),\, \epsilon \in \BC^\times,\, \varepsilon \in \BC.
\end{equation*}

It follows from the proof of Proposition~\ref{W_lim} that the family $X_{\epsilon}(\chi_0+\varepsilon \chi)$ extends to $\epsilon=0$ via $X_0(\chi_0+\varepsilon \chi)=\tilde{\Omega}_{\chi_0+\varepsilon \chi}$.
We conclude that  every limit algebra to $(0,\varepsilon)$ contains $\tilde{\Omega}_{\chi_0+\varepsilon\chi}$ (here $\varepsilon$ may be equal to zero). Note also that $h_{ij} \in B_{\epsilon}(\on{exp}(\epsilon(\chi_0+\varepsilon \chi)))$ so every limit algebra contains $h_{ij}$. It follows that every limit algebra to $(0,\varepsilon)$ ($\varepsilon \in \BC$) contains $\tilde{\Omega}_{\chi_0+\varepsilon\chi}, h_{ij}$ so is contained in the centralizer $Z_{U(\mathfrak{g}[t])}(\tilde{\Omega}_{\chi_0+\varepsilon\chi}, h_{ij})$ that is equal to $\CA^{\mathrm{u}}_{\chi_0+\varepsilon\chi}$ for $\varepsilon \neq 0$ (see Proposition \ref{ca_as_centr}) and to $\CA^{\mathrm{u}}_{(\chi_0,\chi)}$ for $\varepsilon=0$ (see Lemma \ref{central_chi_chi0}).  From the dimension estimates (same as in the proof of Theorem \ref{lim_Bethe}) we conclude that every limit algebra to $(0,\varepsilon)$ ($\varepsilon \in \BC$) is $\CA^{\mathrm{u}}_{\chi_0+\varepsilon\chi}$ for $\varepsilon \neq 0$ and  $\CA^{\mathrm{u}}_{(\chi_0,\chi)}$ for $\varepsilon=0$. We also note that every limit algebra to $(\epsilon,0)$ ($\epsilon \neq 0$) contains $B_{\epsilon}(\on{exp}(\epsilon \chi_0)),\,h_{ij}$, hence, by  Proposition~\ref{lim_B} 
and the maximality of $B_\epsilon(\on{exp}(\epsilon \chi_0,\chi))$ (see \cite[Main Theorem of Section 1.2]{ir0}) it must be equal to $B_\epsilon(\on{exp}(\epsilon \chi_0,\chi))$.

To finish the proof we need to check that the family that we have constructed is indeed continuous. Let $\pi$ be our two-parameter family and let $\pi_i\colon \BC^2 \ra \on{Gr}(d_i,Q_i)$ be the corresponding maps to Grassmannians. We need to check that if we have some convergent sequence $p_n \ra p$ on $\BC^2$ then the sequence of images $\pi_i(p_n)$ converges to $\pi_i(p)$. It is enough to show that every convergent subsequence of the sequence $\pi_i(p_n)$ converges to $\pi_i(p)$. Pick a convergent subsequence $\pi_i(p_{j^{(i)}_n})$ of the sequence $\pi_i(p_n)$ and assume that it converges to some $P_i$. Pick a   subsequence $p_{j^{(i)}_{j^{(i+1)}_n}}$ of $p_{j_n}$ such that $\pi_i(p_{j^{(i)}_{j^{(i+1)}_n}})$ converges to some $P_{i+1}$. Continuing in this way we obtain a sequence of vector spaces (use Lemma \ref{cont_ver_cor}) $P_i \subset P_{i+1} \subset \ldots $ such that $P:=\bigcup_{j \geqslant i}P_i$ is a commutative algebra (use Lemma \ref{cont_ver_cor}), its dimension can only jump. It now follows from the observations above together with Lemma \ref{cont_ver_cor} that the algebra $P$ must be $\pi(p)$. We conclude that $P_i=\pi_i(p)$.

\epr

\rem{}
{\em{Note that there is an alternative proof of Lemma \ref{extend_C2} that uses explicit generators described in Section \ref{gen_A_u_chi_exp_resid} (compare with the proof of Proposition \ref{expl_gen_A_u_chi}). Note also that we are not claiming that the family constructed in Lemma \ref{extend_C2} is holomorphic (this should be true but we do not need this).}}
\erem

Pick distinct points $z_1,\ldots,z_k \in \BC$.  We set $\CA_{(\chi_0,\chi)}(\ul{z}):=\underset{\epsilon \ra 0}{\on{lim}}\,\CA_{\chi_0+\epsilon \chi}(\ul{z}) \subset U(\mathfrak{g})^{\otimes k}$. Let us show that the algebra $\CA_{(\chi_0,\chi)}(\ul{z})$ is equal to $\on{ev}_{\ul{z}}(\CA^{\mathrm{u}}_{(\chi_0,\chi)})$ (see Proposition \ref{ev_of_A_univ_chi_0} below). We start from the following lemma.
\lem{}\label{gen_A_z_chi_0_chi}          
We have $\CA_{(\chi_0,\chi)}(\ul{z})=\CA_{\chi_0}(\ul{z}) \otimes_{Z(U(\mathfrak{z}_{\mathfrak{gl}_n}(\chi)))} \CA_{\chi}(\mathfrak{z}_{\mathfrak{gl}_n}(\chi_0))$.    
The algebra     
$\CA_{(\chi_0,\chi)}(\ul{z})$ is  generated by  $\CA_{\chi_0}(\ul{z})$
and the element $\Delta^k(h_{ij})=h_{ij}^{(1)}+\ldots+h_{ij}^{(k)}$.
\elem
\prf
Same proof as of Proposition~\ref{lim_B}.
\epr

\prop{}\label{ev_of_A_univ_chi_0}
We have $\on{ev}_{\ul{z}}(\CA^{\mathrm{u}}_{(\chi_0,\chi)})=\CA_{(-\chi_0,-\chi)}(\ul{z})$.
\eprop
\prf
Use Proposition~
\ref{realiz_our_via_conf}
together with Proposition~\ref{lim_B}, Lemma~\ref{gen_A_z_chi_0_chi} and the maximality of $\CA_{(-\chi_0,-\chi)}(\ul{z}) \subset U(\mathfrak{g})^{\otimes k}$ (follows from the fact that $\CA_{-\chi_0} \subset U(\mathfrak{g})^{\mathfrak{z}_{\mathfrak{g}}(\chi_0)}$
is maximal that was proved in~\cite[Proposition~9.12]{hkrw}). 
\epr

\section{Simplicity of spectra of $B(X)$ on a tensor product of KR modules}\label{simpl_spec_kr}
In this section we recall the main results of~\cite{imr} that will be used later. 

 Pick $k \geqslant 1$ and let $V_{a_1\varpi_{b_1}},\ldots,V_{a_k\varpi_{b_k}}$ be irreducible polynomial representations of $\mathfrak{gl}_n$, corresponding to rectangular Young diagrams. 
Let $\ul{z}=(z_1,\ldots,z_k)$ be a $k$-tuple of (distinct) complex numbers $z_i \in \BC$, later we will put some restrictions on $\on{Re}(z_j)$. We denote by $V_{a_j\varpi_j}(z_j)$ the corresponding irreducible representations of $Y({\mathfrak{gl}_n})$ (via ${\bf{ev}}_{z_j}$). The representations  $V_{a_j\varpi_j}(z_j)$ are called Kirillov--Reshetikhin modules.
Recall that to every $C \in G$, we can associate Bethe subalgebra $B(C) \subset Y(\mathfrak{sl}_n)$ (see Section~\ref{type_A_Yang_Bethe}) so we obtain the action $B(C) \curvearrowright V_1(z_1) \otimes \ldots \otimes V_k(z_k)$ (via $(\on{\bf{ev}}_{z_1} \otimes \ldots \otimes \on{\bf{ev}}_{z_k}) \circ \Delta^{k}$).

The following proposition was proved in~\cite{imr}.
\prop{}\label{crit_norm}
For $C \in S$ and $\on{Re}(z_j)=-b_j-n+a_j$ the algebra $B(C)$ acts on the tensor product $V_{a_1\varpi_{b_1}}(z_1) \otimes \ldots \otimes V_{a_k\varpi_{b_k}}(z_k)$ via normal operators.
\eprop

The following proposition holds by~\cite[Section~4.3]{imr}.
\prop{}\label{cyc_fin}
Fix distinct $z_1,\ldots,z_k \in \BC$ and $d_1,\ldots,d_k \in \BC$. For almost every $s \in \BC$  the action of $B(X)$ on $V_{a_1\varpi_{b_1}}(sz_1+d_1) \otimes \ldots \otimes V_{a_k\varpi_{b_k}}(sz_k+d_k)$ has a cyclic vector for $X=C \in T^{\mathrm{reg}}$ and $X=(C_0,\chi)$.
\eprop

\rem{}
{\em{Proposition~\ref{cyc_fin} remains true for any  $X \in \ol{M_{0,n+2}}$.}}
\erem

As a corollary of Proposition~\ref{crit_norm} and Proposition~\ref{cyc_fin}, we obtain the following proposition.
\prop{}\label{fin_bad}
Fix distinct  $z_1,\ldots,z_k \in i\mathbb{R}$ and let $d_j:=-b_j-n+a_j$. Then the set of $s \in \mathbb{R}$ such that the action of $B(X)$ on $V_{a_1\varpi_{b_1}}(sz_1+d_1) \otimes \ldots \otimes V_{a_k\varpi_{b_k}}(sz_k+d_k)$ does not have simple spectrum for some $X=C \in S^{\mathrm{reg}}$ or $X=(C_0,\chi)$ is finite. 
\eprop

\rem{}
{\em{In the Proposition~\ref{fin_bad} one can replace $X$ by any  $X \in \ol{M_{0,n+2}}$.}}
\erem

From now on we fix distinct $z_1,\ldots,z_k \in i\BR$. Using Proposition~\ref{fin_bad} we get \cor{}\label{cond_simp_spec!}
There exists a real number $N \in \BR$ such that for every $s \in \BR_{\geqslant N}$ and $X \in S^{\mathrm{reg}}$ or $X=(C_0,\chi)$ the action of $B(X)$ on $V_{a_1\varpi_{b_1}}(sz_1+d_1) \otimes \ldots \otimes V_{a_k\varpi_{b_k}}(sz_k+d_k)$ has simple spectrum.
\ecor

\section{Further properties of the two-parameter family $B_{\epsilon}(\on{exp}(\epsilon(\chi_0+\varepsilon \chi)))$}\label{further_prop_of_family}


In this section we construct certain families relating images in $\on{End}(V_{\la_1}\otimes \ldots \otimes V_{\la_k})$ of Bethe algebras and inhomogeneous Gaudin algebras.
We consider $T$ as a subset of $\mathfrak{h}$.
\lem{}\label{ev_of_Bethe} 
We have ${\bf{ev}}_{z}(\tilde{B}(C))=\tilde{\CA}_{C^{-1}}$.
\elem
\prf
Follows from~\cite[Section~2]{no}, see also~\cite{imr}.
\epr

\prop{}\label{lim_B_var_epsilon_compl}
For $\chi \in \mathfrak{h}_{0}^{\mathrm{reg}}$ we have 
\begin{equation*}
\underset{\varepsilon \ra 0}{\on{lim}}\,B_{\varepsilon}(\on{exp}(\chi))=U(\mathfrak{h}_0[t]).    
\end{equation*}

Moreover, for Weil generic $\chi \in \mathfrak{h}_{0}^{\mathrm{reg}}$ (i.e. for $\chi$ in the complement to a certain
countable union of Zariski closed subsets of $\mathfrak{h}_{0}^{\mathrm{reg}}$) the two-parameter family $B_{\varepsilon}(\on{exp}(\epsilon \chi))$ depending on $\epsilon,\varepsilon \in \mathbb{C}^\times$ extends to the continuous two-parameter family on the blowup $\on{Bl}_{(0,0)}\mathbb{C}^2$ as follows (here $\epsilon \in \mathbb{C}$, $\varepsilon,c \in \mathbb{C}^\times$):
\begin{equation*}
((\epsilon,0),[1:0]) \mapsto U(\mathfrak{h}_0[t]),\, ((0,\varepsilon),[0:1]) \mapsto B_{\varepsilon}(1) \cdot \CA_{\chi},     
\end{equation*}
\begin{equation*}
((0,0),[c:1]) \mapsto \CA^{\mathrm{u}}_{c\chi}, \,
((0,0),[0:1]) \mapsto \CA^{\mathrm{u}}_0 \cdot \CA_{\chi}.
\end{equation*}
\eprop
\prf
It follows from Lemma~\ref{cart_Bethe} that $\mathfrak{h}_0 \subset B(\on{exp}(\epsilon \chi))$. Since $\on{deg}_2\mathfrak{h}_0=0$ we conclude that $\mathfrak{h}_0$ lies in every limit algebra of the family $B_{\varepsilon}(\on{exp}(\epsilon \chi))$.
Let us also note that according to~\cite{i}  the algebra $B(\on{exp}(\epsilon \chi))$ contains elements 
\begin{equation*}
\sigma_{i}(\epsilon):=2J(t_{\varpi_i})-\sum_{\al \in \Delta_+} \frac{e^\al(\on{exp}(\epsilon \chi))+1}{e^{\al}(\on{exp}(\epsilon \chi))-1}(\al,\al_i)x_\al^+ x_\al^-.    
\end{equation*}
We denote by $\sigma_{i,\varepsilon}(\epsilon)$ the corresponding elements of $B_\varepsilon(\on{exp}(\epsilon \chi))$.
We see that the family $\sigma_{i,\varepsilon}(\epsilon)$ extends to  $\on{Bl}_{(0,0)}\BC^2$ as follows:
\begin{equation*}
((\epsilon,0),[1:0]) \mapsto 2 t_{\varpi_i}[1],\, ((0,\varepsilon),[0:1]) \mapsto -\sum_{\al \in \Delta_+}\frac{(\al,\al_i)}{(\al,\chi)}x_{\al}^+x_{\al}^-,  
\end{equation*}
\begin{equation*}
((0,0),[c:1]) \mapsto   2t_{\varpi_i}[1]-\sum_{\al \in \Delta_+}\frac{(\al,\al_i)}{c(\al,\chi)}x_{\al}^+x_{\al}^-,\,  
\end{equation*}
We conclude that every limit algebra at the point $((\epsilon,0),[1:0]),\,\epsilon \in \BC$  contains $\mathfrak{h}_0[1]$ so lies in $Z_{U(\mathfrak{g}[t])}(\mathfrak{h}_0[1])=U(\mathfrak{h}_0[t])$, hence, coincides with $U(\mathfrak{h}_0[t])$ since $\on{dim}_{F_1} U(\mathfrak{h}_0[t])=\on{dim}_{F_1}B_{\varepsilon}(\on{exp}(\epsilon \chi))$. Similarly, every limit algebra at the point $((0,\varepsilon),[0:1])$ clearly contains $B(1)$ and also contains the linear combinations of elements $\sum_{\al \in \Delta_+}\frac{(\al,\al_i)}{(\al,\chi)}x_{\al}x_{\al}^-$ i.e. the quadratic part of $\CA_\chi$ (see~\cite{v}).
Note also that the centralizer of $U(\mathfrak{g})^{\mathfrak{g}}$ in $Y(\mathfrak{g})$ is equal to $Y(\mathfrak{g})^{\mathfrak{g}} \cdot U(\mathfrak{g})=Y(\mathfrak{g})^{\mathfrak{g}} \otimes_{U(\mathfrak{g})^{\mathfrak{g}}} U(\mathfrak{g})$ (see~\cite[Lemma~6.14]{ir2}).
Since for Weil generic $\chi \in \mathfrak{h}^{\mathrm{reg}}_{0}$ the algebra $\CA_\chi$ is equal to the centralizer of its quadratic part (see~\cite{r3}) we then conclude that  for Weil generic $\chi$ every limit algebra at the point $((0,\varepsilon),[0:1])$ is contained in $B_{\varepsilon}(1) \cdot \CA_\chi$, hence, is equal to $B_{\varepsilon}(1) \cdot \CA_\chi$ since $\on{dim}_{F_1}B_{\varepsilon}(\on{exp}(\epsilon \chi))=\on{dim}_{F_1}(B_{\varepsilon}(1) \cdot \CA_\chi)$ (compare with the proof of~\cite[Theorem~6.13]{ir2}).
Similarly, every limit algebra at the point $((0,0),[0:1])$ contains the elements $\sum_{\al \in \Delta_+}\frac{(\al,\al_i)}{(\al,\chi)}x_{\al}x_{\al}^-$ and also the element $\tilde{\omega}_0$ so coincides with $\CA^{\mathrm{u}}_0 \cdot \CA_\chi$ (here we use the maximality of $\CA^{\mathrm{u}}_0 \cdot \CA_\chi \subset U(\mathfrak{g}[t])$ which follows from Proposition~\ref{gr_to_univ}).
Finally, every limit algebra at the point $((0,0),[c:1])$ contains the element $\tilde{\omega}_{c\chi}$ therefore is contained in the centralizer $Z_{U(\mathfrak{g}[t])}(\tilde{\omega}_{c\chi})$ that is equal to $\CA^{\mathrm{u}}_{c\chi}$ by Remark~\ref{cent_til_w_chi_reg_chi}. The argument as in the end of the proof of Lemma \ref{extend_C2} finishes the proof.
\epr

\rem{}
{\em{
We expect that Proposition \ref{lim_B_var_epsilon_compl} is true for every $\chi \in \mathfrak{h}_0^{\mathrm{reg}}$. To prove this it is enough to show that the algebra $\CA_\chi$ coincides with the centralizer of its quadratic part for every $\chi \in \mathfrak{h}_0^{\mathrm{reg}}$.   Note also that we are not claiming that the family constructed in  Proposition \ref{lim_B_var_epsilon_compl} is holomorphic (this should be true but we do not need this for our purposes).
}}
\erem

\cor{}\label{lim_B_var_epsilon}
For Weil generic $\chi \in i\mathfrak{h}_{0,\BR}^{\mathrm{reg}}$ (i.e., for $\chi$ in the complement to a certain
countable union of Zariski closed subsets in $i\mathfrak{h}_{0,\mathbb{R}}^{\mathrm{reg}}$) the two-parameter family $B_{\varepsilon}(\on{exp}(\epsilon \chi))$ depending on $\epsilon,\varepsilon \in \mathbb{R}^\times$ extends to the two-parameter family on the blowup $\on{Bl}_{(0,0)}\mathbb{R}^2$ as follows (here $\epsilon \in \mathbb{R}$, $\varepsilon,c \in \mathbb{R}^\times$):
\begin{equation*}
((\epsilon,0),[1:0]) \mapsto U(\mathfrak{h}_0[t]),\, ((0,\varepsilon),[0:1]) \mapsto B_{\varepsilon}(1) \cdot \CA_{\chi},     
\end{equation*}
\begin{equation*}
((0,0),[c:1]) \mapsto \CA^{\mathrm{u}}_{c\chi}, \,
((0,0),[0:1]) \mapsto \CA^{\mathrm{u}}_0 \cdot \CA_{\chi}.
\end{equation*}
\ecor

Let us now introduce the notation. Recall that we have 
\begin{equation*}
\on{Bl}_{(0,0)}\BR^2=\{((a,b),[c:d]) \in \BR^2 \times \BR\BP^1\,|\, ad=bc\}.
\end{equation*}
We consider $\on{Bl}_{(0,0)}\BR^2$ as a topological space (with the standard induced topology on it).
We denote by $\on{Bl}_{(0,0)}(\BR^2 \setminus (\{0\} \times \BR^\times))$ the following topological space:
\begin{equation*}
\on{Bl}_{(0,0)}(\BR^2 \setminus (\{0\} \times \BR^\times)):=
\{((a,b),[c:d]) \in (\BR^2 \setminus (\{0\} \times \BR^\times)) \times \BR\BP^1\,|\, ad=bc\}
\end{equation*}
 
\rem{}
{\em{Note that the topological space $\on{Bl}_{(0,0)}(\BR^2 \setminus (\{0\} \times \BR^\times))$ is not an algebraic variety (it is a constructible set).}}
\erem 
 
We also consider the following subset of $\on{Bl}_{(0,0)}(\BR^2 \setminus (\{0\} \times \BR^\times))$:
\begin{multline}\label{main_bl_param}
\on{Bl}_{(0,0)}(\BR \times [0,1/N) \setminus (\{0\} \times (0,1/N)))):=\\
=\{((a,b),[c:d]) \in (\BR \times [0,1/N) \setminus (\{0\} \times (0,1/N))) \times \mathbb{R}\BP^1\,|\, ad=bc\}.
\end{multline}
To simplify notation, we set $K:=\BR \times [0,1/N) \setminus (\{0\} \times (0,1/N))$ and will denote the space~(\ref{main_bl_param})  by $\on{Bl}_{(0,0)}K$.

We are now ready to prove the proposition that describes a partial compactification of the family $\on{Im}({\bf{ev}}_{\ul{z}/\varepsilon+\ul{d}}B(\on{exp}(\chi)) \ra \on{End}(V_{1} \otimes \ldots \otimes V_{k}))$.
Note that if representations $V_1,\ldots,V_k$ have no multiplicities for the action of $\mathfrak{h}_0$ then the proof of Proposition~\ref{two_stupid_limits} below simplifies. Indeed, if this is the case then $\on{Im}(U(\mathfrak{h}_0)^{\otimes k} \ra \on{End}(V_1 \otimes \ldots \otimes V_k))$ is already the maximal commutative subalgebra of $\on{End}(V_1 \otimes \ldots \otimes V_k)$ so the statement of Proposition~\ref{two_stupid_limits} can be deduced from Corollary~\ref{lim_B_var_epsilon}.

\prop{}\label{two_stupid_limits}
For $z_1,\ldots,z_k \in i\mathbb{R}$ and $\chi \in i\mathfrak{h}_{0,\BR}^{\mathrm{reg}}$
we have $(\varepsilon \in (0,1/N))$
\begin{equation*}
\underset{\varepsilon \ra 0}{\on{lim}}\, \on{Im}({\bf{ev}}_{\ul{z}/\varepsilon+\ul{d}}B(\on{exp}(\chi)) \ra \on{End}(V_{1} \otimes \ldots \otimes V_{k}))=\on{Im}(\CA_{\on{exp}(-\chi)}^{\otimes k} \ra  \on{End}(V_{1} \otimes \ldots \otimes V_{k})).
\end{equation*}
Moreover, for Weil generic $\chi \in i\mathfrak{h}^{\mathrm{reg}}_{0,\mathbb{R}}$  
the two-parameter family \begin{equation*}
\on{Im}({\bf{ev}}_{\ul{z}/\varepsilon+\ul{d}}B(\on{exp}({\epsilon}\chi)) \ra \on{End}(V_{1} \otimes \ldots \otimes V_k))
\end{equation*}
depending on $\epsilon \in \mathbb{R}^\times,$ $\varepsilon \in (0,1/N)$ extends to the (continuous) two-parameter family on  $\on{Bl}_{(0,0)}K$ as follows (here $c \in \mathbb{R}^\times,\, \epsilon \in (0,1/N)$):
\begin{equation*}
((\epsilon,0),[1:0]) \mapsto 
\on{Im}(\CA_{\on{exp}(-\epsilon\chi)}^{\otimes k} \ra \on{End}(V_{1} \otimes \ldots \otimes V_{k})), 
\end{equation*}
\begin{equation*}
((0,0),[1:0]) \mapsto \on{Im}(\CA_{\chi}^{\otimes k} \ra \on{End}(V_1 \otimes \ldots \otimes V_k)),
\end{equation*}
\begin{equation*}
((0,0),[c:1]) \mapsto \on{Im}(\CA_{-c\chi}(\ul{z})  \ra \on{End}(V_{1} \otimes \ldots \otimes V_{k})),
\end{equation*}
\begin{equation*}
((0,0),[0:1]) \mapsto \on{Im}(\CA_{0}(\ul{z}) 
\cdot \Delta^k(\CA_{\chi}) \ra \on{End}(V_{1} \otimes \ldots \otimes V_{k})).
\end{equation*}
\eprop
\prf
Consider the topological space 
\begin{multline*}
\on{Bl}_{(0,0)}(\BR^2 \setminus (\BR^\times \times \{0\} \cup \{0\} \times \BR^\times)):=\\=	
\{((a,b),[c:d]) \in (\BR^2 \setminus (\BR^\times \times \{0\} \cup \{0\} \times \BR^\times)) \times \BR\BP^1\,|\, ad=bc\}.
\end{multline*}

Let us first of all note that ${\bf{ev}}_{\ul{z}/\varepsilon+\ul{d}}B(\on{exp}(\epsilon \chi))={\bf{ev}}_{\varepsilon;\ul{z}+\varepsilon\ul{d}}B_{\varepsilon}(\on{exp}(\epsilon \chi))$. Note also that by our assumptions on $V_i$ the algebra $\on{Im}({\bf{ev}}_{\ul{z}/\varepsilon+\ul{d}}B(\on{exp}({\epsilon}\chi)) \ra \on{End}(V_{1} \otimes \ldots \otimes V_k))
$ is maximal commutative subalgebra of $\on{End}(V_1 \otimes \ldots V_k)$ of dimension $\on{dim}(V_1 \otimes \ldots \otimes V_k)$.
Recall that the family $B_{\varepsilon}(\on{exp}(\epsilon \chi))$ extends to $\on{Bl}_{(0,0)}(\BR^2)$. We claim that the family of morphisms ${\bf{ev}}_{\varepsilon;\ul{z}+\varepsilon\ul{d}}$ also extends to  $\on{Bl}_{(0,0)}(\BR^2)$ by $\on{ev}_{\ul{z}}$. Indeed, note that ${\bf{ev}}_{\varepsilon;\ul{z}+\varepsilon\ul{d}}$ does not depend on $\epsilon$ and $B_\varepsilon (\on{exp}(\epsilon \chi)) \subset Y_{\varepsilon}(\mathfrak{g})$ so the claim easily follows from Lemma~\ref{compat_ev} (see Corollary \ref{ext_ev_hbar}).
Note that the images of $\CA_{-c\chi}(\ul{z}),\, (\CA_{0}(\ul{z})\cdot \Delta^k(\CA_{\chi}))$ in $\on{End}(V_{1} \otimes \ldots \otimes V_{k})$ are maximal commutative subalgebras (follows from~\cite[Corollary~11.9]{hkrw} using Lemma~\ref{homog_comp_Gaud}).

Taking the image in $\on{End}(V_1 \otimes \ldots \otimes V_k)$ of the family considered in   Corollary \ref{lim_B_var_epsilon} and using that $\on{Bl}_{(0,0)}\BR^2$ is Noetherian (so the image of the whole family coincides with the image of a large enough filtration term of it)  we conclude that the two-parameter family 
\begin{equation*}
\on{Im}({\bf{ev}}_{\ul{z}/\varepsilon+\ul{d}}B(\on{exp}({\epsilon}\chi)) \ra \on{End}(V_{1} \otimes \ldots \otimes V_k))
\end{equation*}
depending on $(\epsilon,\varepsilon) \in (\mathbb{R}^\times)^2$ extends to the two-parameter family on  $\on{Bl}_{(0,0)}(\mathbb{R}^2 \setminus (\BR^\times \times \{0\} \cup \{0\} \times \mathbb{R}^\times)) \setminus \{((0,0),[1:0])\}$ as follows (here $c \in \mathbb{R}^\times,\, \epsilon \in \mathbb{R}$):
\begin{equation*}
((0,0),[c:1]) \mapsto \on{Im}(\CA_{-c\chi}(\ul{z})  \ra \on{End}(V_{1} \otimes \ldots \otimes V_{k})),
\end{equation*}
\begin{equation*}
((0,0),[0:1]) \mapsto \on{Im}((\CA_{0}(\ul{z}) 
\cdot \Delta^k(\CA_{\chi})) \ra \on{End}(V_{1} \otimes \ldots \otimes V_{k}).
\end{equation*}
It remains to show that this family extends to $\on{Bl}_{(0,0)}K$ by setting  
\begin{equation*}
((\epsilon,0),[1:0]) \mapsto 
\on{Im}(\CA_{\on{exp}(-\epsilon\chi)}^{\otimes k} \ra \on{End}(V_{1} \otimes \ldots \otimes V_{k})), 
\end{equation*}
\begin{equation*}
((0,0),[1:0]) \mapsto \on{Im}(\CA_{\chi}^{\otimes k} \ra \on{End}(V_1 \otimes \ldots \otimes V_k)).
\end{equation*}


Note that the images of algebras ${\bf{ev}}_{\ul{z}/\varepsilon+\ul{d}}(B(\on{exp}(\epsilon \chi)))$, $\CA^{\otimes k}_{\on{exp}(-\epsilon \chi)}$, $\CA_\chi^{\otimes k}$, $\CA_{-c\chi}(\ul{z})$, $(\CA_0(\ul{z}) \cdot \Delta^k(\CA_\chi))$ in $\on{End}(V_1 \otimes \ldots \otimes V_k)$ are maximal commutative so we can replace algebras above by the corresponding algebras for $\mathfrak{gl}_n$. 
 Consider the following open subset $U \subset \on{Bl}_{(0,0)}K$:
\begin{equation*}
U:=\{((\epsilon,c\epsilon),[1:c])\} \subset \on{Bl}_{(0,0)}K.\end{equation*}
Recall now the generators $r_{k,z_i,l,\epsilon,\varepsilon}$ of the algebra ${\bf{ev}}_{\ul{z}/\varepsilon+\ul{d}}(\tilde{B}(\on{exp}(\epsilon \chi)))$ (see Section~\ref{gen_A_u_chi_exp_resid}). Recall also the generators $r^0_{k,z_i,l,\xi}$ of the algebra $\tilde{\CA}_{\xi}(z_1,\ldots,z_k)$, here $\chi \in \mathfrak{h}_0$.
Let us first of all note that the family $r_{k,u_i,l,\epsilon,\varepsilon}$ extends continuously to the family on $U \setminus \{((\epsilon,0),[1:0])\}$ by setting
\begin{equation*}
((0,c\epsilon),[1:c])  \mapsto r^0_{k,z_i,l,-c^{-1}\chi} \in \CA_{-\chi}(c^{-1}z_1,\ldots,c^{-1}z_k).
\end{equation*}
To see this, recall that $\varepsilon=c\epsilon$ and note that we have the equality (compare with the proof of Proposition~\ref{expl_gen_A_u_chi}): 
\begin{equation}\label{nicest_limit!}
\underset{\epsilon \ra 0}{\on{lim}}\,\epsilon^{-n}\on{cdet}(e^{-\epsilon \partial_u}(1+L_{\ul{z}/\varepsilon+\ul{d}}(u/\epsilon))-\on{exp}(-\epsilon \chi))=\on{cdet}(L_{c^{-1}\ul{z}}(u)-\partial_u+\chi).
\end{equation}
The equality~(\ref{nicest_limit!}) implies that $\underset{\epsilon \ra 0}{\on{lim}}\,b_{k,l,\epsilon,\varepsilon}(u)=b^0_{k,l,c^{-1}\ul{z},-\chi}(u)$, hence, 
\begin{multline*}
\underset{\epsilon \ra 0}{\on{lim}}\,r_{k,z_i,l,\epsilon,\varepsilon}=\underset{\epsilon \ra 0}{\on{lim}}\,\sum_{j=1}^k \on{res}_{u=c^{-1}z_i+\epsilon j +\epsilon d_i}(u-z_i)^lb_{k,\epsilon,\varepsilon}(u)du=\\
=\on{res}_{u=c^{-1}z_i}(u-z_i)^lb^0_{k,l,c^{-1}\ul{z},-\chi}=r^0_{k,c^{-1}z_i,l,-\chi}.
\end{multline*}

It remains to show that the families $r_{k,z_i,l,\epsilon,\varepsilon}$ extend to the points $((\epsilon,0),[1:0])$ and generate algebras $\tilde{\CA}_{\on{exp}(-\epsilon \chi)}^{\otimes k}$, $\tilde{\CA}_{\chi}^{\otimes k}$. Indeed, recall that the elements $r_{k,z_i,l,\epsilon,\varepsilon}$ are sums of residues of the elements $(u-z_i)^lb_{k,\epsilon,\varepsilon}$ at the points $c^{-1}z_i+\epsilon j +\epsilon d_i$. Recall also that the generating function for  $b_{k,\epsilon,\varepsilon}$ can be written in the following form:
\begin{multline*}
\epsilon^{-n}\on{tr}A_n(e^{-\epsilon\partial_u}(1+L_{\ul{z}/\varepsilon+\ul{d}}(u))_1-\on{exp}(-\epsilon \chi)_1)\ldots (e^{-\epsilon\partial_u}(1+L_{\ul{z}/\varepsilon+\ul{d}}(u))_n-\on{exp}(-\epsilon \chi)_n) = \\ =\epsilon^{-n}\sum_{a=0,\ldots,n}(-1)^{n-a}{n\choose{a}}\on{tr}A_{n}
(1+\frac{\epsilon E^{(1)}_1}{u-c^{-1}z_1-\epsilon d_1})\ldots (1+\frac{\epsilon E^{(1)}_a}{u-c^{-1}z_1-\epsilon(d_1+a-1)}) \ldots \\
\ldots  (1+\frac{\epsilon E^{(k)}_1}{u-c^{-1}z_k-\epsilon d_k}) \ldots (1+\frac{\epsilon E^{(k)}_a}{u-c^{-1}z_k-\epsilon(d_k+a-1)})(\on{exp}(\epsilon \chi))_{a+1} \ldots (\on{exp}(\epsilon \chi))_n e^{-a\epsilon\partial_u}
\\=\sum_{k \geqslant 0}b_{k,\epsilon,\varepsilon}(u)\partial_u^k.
\end{multline*}
Since $c^{-1} \ra \infty$ and $\epsilon$ is clearly bounded we see that at the point $((\epsilon,0),[1:0])$ the family $r_{k,z_i,l,\epsilon,\varepsilon}$ converges to the same element as the family $r_{1,z_i,l,\epsilon,\varepsilon}^{(i)}=1\otimes \ldots \otimes 1 \otimes \underset{i}{r_{1,z_i,l,\epsilon,\varepsilon}} \otimes 1 \otimes  \ldots \otimes 1$. So we can assume that $k=1$ and then we are dealing with the family ${\bf{ev}}_{z_{i}/\varepsilon+d_i}(\tilde{B}(\on{exp}(\epsilon\chi)))$ that by Lemma~\ref{ev_of_Bethe} is exactly $\tilde{\CA}_{\on{exp}(-\epsilon \chi)}$. It remains to show that the limits as $c \ra 0$ of the elements $r_{1,z_i,l,\epsilon,\varepsilon}$ generate $\tilde{\CA}_{\on{exp}(-\epsilon \chi)}$ for $\epsilon \neq 0$ and generate $\CA_{-\chi}$ as $\epsilon \ra 0$. Indeed, directly from the definitions the elements $r_{1,z_i,l,\epsilon,\varepsilon}$ do not depend on $c$, $d_i$ (compare with Remark~\ref{shift_inv_ev_Bethe}). Note also that as $\epsilon \ra 0$ we have 
\begin{equation*}
\underset{\epsilon \ra 0}{\on{lim}}\,\on{cdet} \epsilon^{-1}(e^{-\epsilon \partial_u}(1+L_{0}(u))-\on{exp}(\epsilon\chi))	=L_0(u)-\partial_u+\chi.
\end{equation*}
This observation finishes the proof. 


\epr

\rem{}
{\em{We expect that for every $\chi \in \mathfrak{h}^{\mathrm{reg}}_0$ the family ${\bf{ev}}_{\ul{z}/\varepsilon +\ul{d}} B(\on{exp}(\epsilon \chi))$ extends to the blowup $\on{Bl}_{(0,0)}\BC^2$ by setting (here $\epsilon, \varepsilon, c \in \BC^\times$)
\begin{equation*}
((\epsilon,0),[1:0]) \mapsto \CA_{\on{exp}(-\epsilon \chi)}^{\otimes k}    ,\, ((0,\varepsilon),[1:0]) \mapsto {\bf{ev}}_{\ul{z}/\varepsilon+\ul{d}}(B(1)) \cdot \CA_{\chi}
\end{equation*}
\begin{equation*}
((0,0),[c:1]) \mapsto \CA_{-c\chi}(\ul{z}),\, ((0,0),[0:1]) \mapsto \CA_0(\ul{z}) \cdot \Delta^k(\CA_{\chi}),\, ((0,0),[1:0]) \mapsto \CA_{\chi}^{\otimes k}.    
\end{equation*}
For example, it follows from our proof of Proposition~\ref{two_stupid_limits}  
that for a fixed $\epsilon \in \BC^\times$ we have $\underset{\epsilon \rightarrow 0}{\on{lim}}\,{\bf{ev}}_{\ul{z}/\varepsilon+\ul{d}}(B(\on{exp}(\epsilon \chi)))=\CA_{\on{exp}(-\epsilon \chi)}^{\otimes k}$.
It is also easy to show that $\underset{\epsilon \ra 0}{\on{lim}}\,{\bf{ev}}_{\ul{z}/\varepsilon+\ul{d}}B(\on{exp}(\epsilon \chi)) \supset {\bf{ev}}_{\ul{z}/\varepsilon+\ul{d}}(B(1))\cdot \CA_\chi$ but the equality is not clear. We expect that ${\bf{ev}}_{\ul{z}/\varepsilon+\ul{d}}(B(1)) \cdot \CA_\chi \subset U(\mathfrak{g})^{\otimes k}$ is the maximal commutative subalgebra and then the equality above will follow. Another problem with extending the family ${\bf{ev}}_{\ul{z}/\varepsilon+\ul{d}}B(\on{exp}(\epsilon \chi))$ to the blowup $\on{Bl}_{(0,0)}\BC^2$ is that we do not know Poincar\'e series of the algebras ${\bf{ev}}_{\ul{z}/\varepsilon+\ul{d}}B(\on{exp}(\epsilon \chi))$ so the method similar to the one used in the proof of Proposition~\ref{two_stupid_limits} (see also the proof of \cite[Proposition~10.17]{hkrw}) cannot be applied.}}
\erem

\section{Crystal structure on $\CE_C(\ul{\la})$}\label{cryst_str_E_la}

Let us now pick $k \geqslant 1$ and  let $V_{1},\ldots,V_{k}$ be irreducible polynomial representations of $\mathfrak{gl}_n$ with highest weights $\la_1,\ldots,\la_k$ respectively.  Let us also fix distinct $z_1,\ldots,z_k \in i\BR$. We are not assuming  that $V_j$ correspond to rectangular Young diagrams here. For our purposes it is enough to assume that Corollary~\ref{cond_simp_spec!} holds for tensor products of $V_j$. More precisely, we assume the following.
\ass{}\label{our_assumption}
We assume that there exist $d_1,\ldots,d_k \in \BC$ and $N \in \BR$ such that for every $s \in \BR_{\geqslant N}$, $X \in S^{\mathrm{reg}}$ or $X=(C_0,\chi)$ and $1 \leqslant m_1<\ldots<m_p \leqslant k$ the action of $B(X)$ on $V_{m_1}(sz_{m_1}+d_{m_1}) \otimes \ldots \otimes V_{m_p}(sz_{m_p}+d_{m_p})$ has simple spectrum.
\eass

\rem{}
{\em{By the results of Section~\ref{simpl_spec_kr} this assumption holds automatically for $V_j$, corresponding to rectangular Young diagrams. We will see in Remark~\ref{simple_spec_impl_kr} that  this assumption actually implies that $V_j$ correspond to rectangular Young diagrams.}}
\erem

Recall that we have fixed $z_i,\,d_i$. For simplicity, let us denote the representation $V_{1}(sz_1+d_1) \otimes \ldots \otimes V_{k}(sz_k+d_k)$ of $Y(\mathfrak{g})$ by  $\ul{V}(s),\, s \in \BR_{\geqslant N}$. 
For $X=C \in S^{\mathrm{reg}}$ or $X=(C_0,\chi)$ let $\mathcal{E}_X(\ul{\la},s)$ be the set of eigenlines of $B(X)$ acting on $\ul{V}(s)$. Note that the set $\BR_{\geqslant N}$ is  contractible so simply-connected, hence, $\mathcal{E}_X(\ul{\la},s)$ are canonically identified for all $s \in \BR_{\geqslant N}$ and we denote this set simply by $\CE_X(\ul{\la})$.

Our goal for now is to define on $\CE_C(\ul{\la})$ the structure of  $\hat{\mathfrak{sl}_n}$-crystal (see Definition~\ref{sl_hat_cryst}).

We pick $\hat{w} \in \widehat{W}^{\mathrm{ext}}$ and consider the corresponding alcove $Q_{\hat{w}} \supset Q^{\mathrm{reg}}_{\hat{w}}$. Taking $\chi \in Q^{\mathrm{reg}}_{\hat{w}}$ such that $C=\on{exp}(2\pi i \chi)$ we define a crystal structure on $\mathcal{E}_{\on{exp}(2\pi i \chi)}(\ul{\la})$ as follows. Pick $j=1,\ldots,n$ and consider the corresponding wall $H^{\hat{w}}_j$ and
pick a generic element $\chi_0$ 
lying on the intersection of this wall with $Q_{\hat{w}}$. 
Set $C:=\on{exp}(2\pi i \chi)$, $C_0:=\on{exp}(2\pi i \chi_0)$ and
consider the algebra $B(C_0,\chi)$ acting on $\ul{V}(s)$. Note that we have the identification $\mathcal{E}_C(\ul{\la}) \iso \mathcal{E}_{(C_0,\chi)}(\ul{\la})$. Recall  that by Proposition~\ref{lim_B} the algebra $B(C_0,\chi)$ 
is generated by $B(C_0)$ together with the element 
$
h$, corresponding to the wall $H^{\hat{w}}_j$. We can then decompose $\ul{V}(s)$ as $B(C_0)$-module in the direct sum of weight spaces $\ul{V}(s)=\bigoplus_{\eta\colon B(C_0) \ra \BC} \ul{V}(s)^{\eta}$. The action of  $B(C_0,\chi)$ on $\ul{V}(s)$ has simple spectrum so the action $h \curvearrowright \ul{V}(s)^{\eta}$ has simple spectrum i.e. $\ul{V}(s)^{\eta}=\bigoplus_{i \in \BZ}\ul{V}(s)^{\eta}_i$ with $\ul{V}(s)^{\eta}_i$ being one-dimensional. Moreover, note that $\mathfrak{h}_0 \subset B(C_0,\chi)$ so $\mathfrak{h}_0$ acts on each $\ul{V}(s)^{\eta}_i$ via some weight $\mu \in \mathfrak{h}_0^*$. Since $\ul{V}(s)^{\eta}_i$ is one-dimensional it  can be considered as an element of $\CE_C(\ul{\la})$ and we can define 
\begin{equation}\label{cryst_str_e_la}
\on{e}_{[j]}(\ul{V}(s)^\eta_i):=\ul{V}(s)^\eta_{i+1},\, \on{f}_{[j]}(\ul{V}(s)^{\eta}_i):=\ul{V}(s)^{\eta}_{i-1},\, \on{wt}(\ul{V}(s)^{\eta}_i)=\mu,
\end{equation}
where $\on{e}_{[j]}(\ul{V}(s)^\eta_i):=0$ if $\ul{V}(s)^\eta_{i+1}=0$, similarly, $\on{f}_{[j]}(\ul{V}(s)^{\eta}_i):=0$ if $\ul{V}(s)^\eta_{i-1}=0$.

Let us now note that the definition of $\on{e}_{[j]},\on{f}_{[j]}$ does not depend on the choice of $\chi_0$ because the subset of $Q_{\hat{w}}$, consisting of regular elements and subregular elements of $H_j^{\widehat{w}}$ is simply-connected  (we also use the Main Theorem of  \cite{ir0}, where the parametrizing space of limits of Bethe algebras is described).
We will denote by $\CE^{\hat{w}}_{\on{exp}(2\pi i \chi)}(\ul{\la})$ the set $\CE_{\on{exp}(2\pi i \chi)}(\ul{\la})$ with a structure (\ref{cryst_str_e_la}) defined as above.



\lem{}\label{shift_equiv}
Pick $\sigma \in S_n$ and $(m_1,\ldots,m_n),\, (m_1',\ldots,m_n') \in \La$. Set $\hat{w}=(\sigma;m_1,\ldots,m_n)$, $\hat{w}'=(\sigma;m_1',\ldots,m_n')$. Pick also $\chi \in Q^{\mathrm{reg}}_{\hat{w}},\, \chi' \in Q^{\mathrm{reg}}_{\hat{w}'}$. Then there is a bijection of sets $\CE^{\hat{w}_1}_{\on{exp}(2\pi i \chi_1)}(\ul{\la}) \iso \CE^{\hat{w}_2}_{\on{exp}(2\pi i \chi_2)}(\ul{\la})$ that commutes with $e_{[j]}, f_{[j]}, \on{wt}$. \elem
\prf
The map $\on{exp}\colon \mathfrak{h}_{\BR} \ra S$ sends alcoves $Q_{\hat{w}},\, Q_{\hat{w}'}$ to the same subset of $S$. The image of $H^{\hat{w}}_j$ coincides with the image of $H^{\hat{w}'}_j$. The claim follows.
\epr

Let $\CE^{\hat{w},\mathrm{fin}}_{\on{exp}(2\pi i \chi)}(\ul{\la})$ be the set $\CE^{\hat{w}}_{\on{exp}(2\pi i \chi)}(\ul{\la})$ together with operations  $\on{e}_{[j]},\on{f}_{[j]},\, j=1,\ldots,n-1$, $\on{wt}$ (i.e. we forget $\on{e}_{[n]},\on{f}_{[n]}$).
\prop{}\label{iso_fin_cr_res}
Pick $\sigma \in S_n$.
There is a canonical isomorphism  $\CE^{\sigma,\mathrm{fin}}_{\on{exp}(2\pi i \chi)}(\ul{\la}) \iso \mathcal{E}_{-\chi}(\ul{\la})$ compatible with $\on{e}_j,\on{f}_j,\on{wt}$.
\eprop
\prf  
Consider the family $B_{\epsilon}(\on{exp}(2\pi i \epsilon \chi))$, $\epsilon \in (0,1/N)$. Note that by Theorem~\ref{lim_Bethe} the limit of this family to $\epsilon=0$ is $\CA^{\mathrm{u}}_{2\pi i\chi}$. 
It follows from Lemma \ref{compat_ev} and Corollary \ref{ext_ev_hbar} that the family of homomorphisms ${\bf{ev}}_{\epsilon;z_1+\epsilon d_1,\ldots,z_k+\epsilon d_k}$ ($\epsilon \in \BR$) extends to $\epsilon=0$ via ${\on{ev}}_{\ul{z}}$. Taking the image of the family above in $\on{End}(V_1 \otimes \ldots \otimes V_k)$, we obtain the family of commutative subalgebras of $\on{End}(V_1 \otimes \ldots \otimes V_k)$. Every algebra of our family above acts on $V_1 \otimes \ldots \otimes V_k$ with simple spectrum (for $\epsilon \in (0,1/N)$ this follows from the Assumption \ref{our_assumption}, for $\epsilon=0$ this follows from Proposition \ref{realiz_our_via_conf} together with \cite{ffr}, see also \cite[Secion 11]{hkrw}).
Using  Lemma~\ref{homog_comp_Gaud}, we  obtain the identification of sets
$\mathcal{E}^{\sigma}_{\on{exp}(2\pi i \chi)}(\ul{\la}) \iso \mathcal{E}_{-\chi}(\ul{\la})$.

We now pick $[j]=[1],\ldots,[n-1]$ and our goal is to show that the operators $\on{e}_{[i]},\on{f}_{[i]}, \on{wt}$ are compatible with the identification above. To see that, we can consider the two-parameter family $B_{\epsilon}(\on{exp}(\epsilon(\chi_0+\varepsilon \chi))),\, \epsilon \in (0,1/N),\,\varepsilon \in \BR^\times$ that by Lemma~\ref{extend_C2} extends to a family depending on  $\epsilon \in [0,1/N),\,\varepsilon \in \BR$ as follows: 
\begin{equation*}
(\epsilon,0) \mapsto B_{\epsilon}(\on{exp}(\epsilon \chi_0),\chi),\, (0,\varepsilon) \mapsto \CA^{\mathrm{u}}_{\chi_0+\varepsilon\chi},\, (0,0) \mapsto \CA^{\mathrm{u}}_{(\chi_0,\chi)}.
\end{equation*}
Recall that the family of homomorphisms ${\bf{ev}}_{\epsilon;z_1+\epsilon d_1,\ldots,z_k+\epsilon d_k}$ ($\epsilon \in (0,1/N)$) extends to $\epsilon=0$ via ${\on{ev}}_{\ul{z}}$  so we can consider the image of the family above (with $\epsilon \in [0,1/N)$, $\varepsilon \in \BR$) in $\on{End}(V_1 \otimes \ldots \otimes V_k)$. Every algebra in the image acts on $V_1 \otimes \ldots \otimes V_k$ with simple spectrum: for $(\varepsilon,\epsilon) \in  \BR \times (0,1/N)$ this follows from Assumption \ref{our_assumption}, for $\epsilon=0$ this follows from \cite[Section 11]{hkrw} together with Propositions \ref{realiz_our_via_conf}, \ref{ev_of_A_univ_chi_0}.
The claim follows since every two paths  on $\BR \times [0,1/N)$ from $(1,1)$ to $(0,0)$ are homotopic.
\epr

\rem{}
{\em{Recall that for $x \in Y(\mathfrak{g})$ we have ${\bf{ev}}_{(\epsilon;z_1+\epsilon d_1,\ldots,z_r+\epsilon d_r)}(\psi_\epsilon(x))={\on{ev}}_{z_1/\epsilon+d_1,\ldots,z_r/\epsilon+d_r}(x)$. We use this in the proof of Proposition \ref{iso_fin_cr_res}.}}
\erem

For $[j] \in \BZ/n\BZ$ we denote by $\CE^{\hat{w},[j]}_{\on{exp}(2\pi i \chi)}(\ul{\la})$ the set $\CE^{\hat{w}}_{\on{exp}(2\pi i \chi)}(\ul{\la})$ together with operations $e_{[1]},\ldots,e_{[j-1]},e_{[j+1]},\ldots,e_{[n]},\, f_{[1]},\ldots,f_{[j-1]},f_{[j+1]},\ldots,f_{[n]},\, \on{wt}$ (we forget the action of $\on{e}_{[j]},\on{f}_{[j]}$). Recall that $\tau_{[j]} \in S_n$ is an  element that sends $[i]$ to $[i+j]$.

\prop{}\label{dif_fin_res}
For $\hat{w}=(\sigma;m_1,\ldots,m_n)$ we have an isomorphism  $\CE^{\hat{w},[j]}_{\on{exp}(2\pi i \chi)}(\ul{\la}) \simeq \CE_{-\tau_{[j]}(\chi)}(\ul{\la})$, that identifies  $\on{e}_{[i]}$ with $\on{e}_{i-j}$ and $\on{f}_{[i]}$ with $\on{f}_{i-j}$. 
\eprop
\prf
We have $\hat{w}=(\sigma;m_1,\ldots,m_n)$.
Recall that $Q_{\hat{w}}$ consists of $(a_{[1]},\ldots,a_{[n]})$ such that
\begin{equation*}
a_{\sigma([n])}-m_{[n]}+1 \geqslant a_{\sigma([1])}-m_{[1]} \geqslant \ldots \geqslant a_{\sigma([n])}-m_{[n]}.   
\end{equation*}
Consider the element $\hat{w}':=\hat{w} \cdot (\tau_{[j]},[0,\ldots,0,\underbrace{1,\ldots,1}_{j}])=(\sigma \circ \tau_{[j]},m_1',\ldots,m_n')$. Note that $Q_{\hat{w}}=Q_{\hat{w}'}$ and $H^{\hat{w}}_{[i]}=H^{\hat{w}'}_{[i-j]}$. We obtain the identification $\CE^{\hat{w}}_{\on{exp}(2\pi i \chi)}(\ul{\la}) \iso \CE^{\hat{w}'}_{\on{exp}(2\pi i \chi)}(\ul{\la})$ such that the action of $\on{e}_{[l]},\on{f}_{[l]}$ on $\CE^{\hat{w}}_{\on{exp}(2\pi i \chi)}(\ul{\la})$ corresponds to the action of $\on{e}_{[l-j]},\on{f}_{[l-j]}$ on $\CE^{\hat{w}'}_{\on{exp}(2\pi i \chi)}(\ul{\la})$. So we obtain the identification $\CE^{\hat{w},[j]}_{\on{exp}(2\pi i \chi)}(\ul{\la}) \iso \CE^{\hat{w}',\mathrm{fin}}_{\on{exp}(2\pi i \chi)}$. 
It  follows from Lemma~\ref{shift_equiv} that we have an identification $\CE^{\hat{w}',\mathrm{fin}}_{\on{exp}(2\pi i \chi)} \iso \CE^{\sigma \circ \tau_{[j]},\mathrm{fin}}_{\on{exp}(2\pi i \tau_{[j]}(\chi))}$.
Note now that by Proposition~\ref{iso_fin_cr_res} there is an isomorphism of crystals  $\CE^{\sigma \circ \tau_{[j]},\mathrm{fin}}_{\on{exp}(2\pi i \tau_{[j]}(\chi))}(\ul{\la}) \iso \CE_{-\tau_{[j]}(\chi)}(\ul{\la})$. After composing the identifications above we conclude that we have an isomorphism $\CE^{\hat{w},[j]}_{\on{exp}(2\pi i \chi)}(\ul{\la}) \iso \CE_{-\tau_{[j]}(\chi)}(\ul{\la})$ that sends $\on{e}_{[i]} \mapsto \on{e}_{[i-j]},\, \on{f}_{[i]} \mapsto \on{f}_{[i-j]}$.

\epr

\cor{}\label{restr_ident}
We have a canonical isomorphism $\CE^{[j]}_{\on{exp}(2\pi i\chi)}(\ul{\la}) \simeq \CE_{-\tau_{[j]}(\chi)}(\ul{\la})$ that becomes an isomorphism of crystals if we identify $\on{e}_{[i]}$ with $\on{e}_{i-j}$ and $\on{f}_{[i]}$ with $f_{i-j}$. 
\ecor

\section{Identification of $\CE_C(\ul{\la})$ with tensor product of KR crystals}\label{E_la_via_tens_prod}

Pick $\chi \in Q_1$, $k \in \BZ_{\geqslant 1}$. Recall that $\ul{\la}=\la_1,\ldots,\la_k$ is a collection of dominant weights (of $\mathfrak{gl}_n$), $z_1,\ldots,z_k \in i\BR$ and the Assumption \ref{our_assumption} holds.
\prop{}\label{r_1_crystal}
Assume that $k=1$. Then $\la=l\varpi_r$ and the crystal $\mathcal{E}_{\on{exp}(2\pi i \chi)}(\la)$ coincides with the Kirillov--Reshetikhin crystal ${\bf{B}}_\la$, corresponding to the irreducible representation $V_\la$ (see Section~\ref{cryst_main}). 
\eprop
\prf
Follows from  Proposition~\ref{restr_cryst_class}, Corollary \ref{restr_ident} and \cite[Theorem 1.8]{hkrw}.
\epr

\rem{}\label{simple_spec_impl_kr}
{\em{Let us point out that Proposition~\ref{r_1_crystal},
in particular, implies that if the Assumption~\ref{our_assumption} holds for $V_{1},\ldots,V_k$ then each $V_i$ is the irreducible representation, corresponding to a rectangular diagram. 
}}
\erem

Recall that in Section~\ref{cryst_main} the tensor product of  $\hat{\mathfrak{sl}_n}$-crystals was defined.
The goal of the rest of the Section is to construct the isomorphism of crystals $\CE_{\on{exp}(2\pi i \chi)}(\ul{\la}) \iso \CE_{\on{exp}(2\pi i \chi)}(\la_1) \otimes \ldots \otimes  \CE_{\on{exp}(2\pi i \chi)}(\la_k)
$ (see Proposition~\ref{iso_with_tensor}) and then to conclude (using Proposition~\ref{r_1_crystal}) that we have an isomorphism $\CE_{\on{exp}(2\pi i \chi)}(\ul{\la}) \iso {\bf{B}}_{\la_1} \otimes \ldots \otimes {\bf{B}}_{\la_k}$ (see Theorem~\ref{KR_tensor_via_geometry}). 
The similar claims for the $\mathfrak{sl}_n$-crystals $\CE_{\chi}(\ul{\la})$ follow from the results of~\cite{hkrw}.

\prop{}\label{iso_with_tensor}
For $\on{Im}z_1 \gg \on{Im}z_2 \gg \ldots \gg \on{Im}z_k$  we have the canonical isomorphism of crystals 
\begin{multline*}
\CE_{\on{exp}(2\pi i \chi)}(\ul{\la}) \iso \\ 
\CE_{\on{exp}(2\pi i \chi)}(\la_1) \otimes \CE_{\on{exp}(2\pi i \chi)}(\la_2) \otimes \ldots \otimes \CE_{\on{exp}(2\pi i \chi)}(\la_{k-1}) \otimes  \CE_{\on{exp}(2\pi i \chi)}(\la_k).
\end{multline*}
\eprop
\prf
Let us construct an isomorphism of sets  
\begin{equation*}
\CE_{\on{exp}(2\pi i \chi)}(\ul{\la}) \iso \CE_{\on{exp}(2\pi i \chi)}(\la_1) \otimes \ldots \otimes  \CE_{\on{exp}(2\pi i \chi)}(\la_k).
\end{equation*}
Consider the family 
\begin{equation*}
\on{Im}(\on{\bf{ev}}_{(\ul{z}/\varepsilon+\ul{d})}(B(\on{exp}(2\pi i \chi))) \ra \on{End}(V_{1} \otimes \ldots \otimes V_k))
\end{equation*}
with $z_l \in i\BR,\, \varepsilon \in (0,\frac{1}{N})$  (see Assumption~\ref{our_assumption})
and note that by Proposition~\ref{two_stupid_limits} the limit as $\varepsilon \ra 0$ is $\on{Im}(\CA_{\on{exp}(-2\pi i \chi)}^{\otimes k} \ra \on{End}(V_1 \otimes \ldots \otimes V_k))$ 
that is exactly $\on{Im}({\bf{ev}}_{z_1}(B(\on{exp}(2\pi i \chi))) \otimes \ldots \otimes {\bf{ev}}_{z_k}(B(\on{exp}(2\pi i \chi))) \ra \on{End}(V_1 \otimes \ldots V_k))$.
Since the set of eigenlines of the action ${\bf{ev}}_{z_1}(B(\on{exp}(2\pi i \chi))) \otimes \ldots \otimes {\bf{ev}}_{z_k}(B(\on{exp}(2\pi i \chi)))$ is exactly $\CE_{\on{exp}(2\pi i \chi)}(\la_1) \otimes \ldots \otimes \CE_{\on{exp}(2\pi i \chi)}(\la_k)$ we obtain the identification $\CE_{\on{exp}(2\pi i \chi)}(\ul{\la}) \iso \CE_{\on{exp}(2\pi i \chi)}(\la_1) \otimes \ldots \otimes  \CE_{\on{exp}(2\pi i \chi)}(\la_k)$. We claim that this identification is an isomorphism of crystals. To check this, it is enough to show that for every $j=1,\ldots,n$ the identification above induces the isomorphism of crystals $\CE^{[j]}_{\on{exp}(2\pi i \chi)}(\ul{\la}) \iso \CE^{[j]}_{\on{exp}(2\pi i \chi)}(\la_1) \otimes \ldots \otimes  \CE^{[j]}_{\on{exp}(2\pi i \chi)}(\la_k)$. 
Consider the element $\hat{w}_j:=(\tau_{[j]};0,\ldots,0,\underbrace{1,\ldots,1}_j) \in \hat{W}^{\mathrm{ext}}$ and let us also fix an element $\chi_j \in Q_{\tau_{[j]}}^{\mathrm{reg}}$.
It follows from Lemma~\ref{shift_equiv} and the proof of Proposition~\ref{dif_fin_res} that we have identifications $\CE^{[j]}_{\on{exp}(2\pi i \chi)}(\ul{\la}) \iso \CE^{\hat{w}_j,\mathrm{fin}}_{\on{exp}(2\pi i \chi)}(\ul{\la}) \iso \CE^{\tau_{[j]},\mathrm{fin}}_{\on{exp}(2\pi i \chi_j)}(\ul{\la})$, 
 $\CE^{[j]}_{\on{exp}(2\pi i \chi)}(\la_l) \iso \CE^{\hat{w}_j,\mathrm{fin}}_{\on{exp}(2\pi i \chi)}(\la_l) \iso \CE^{\tau_{[j]},\mathrm{fin}}_{\on{exp}(2\pi i \chi_j)}(\la_l)$, $l=1,\ldots,k$ that are compatible with the identifications above. Therefore it is enough to show that the bijection 
 $\CE^{\tau_{[j]},\mathrm{fin}}_{\on{exp}(2\pi i \chi_j)}(\ul{\la}) \iso \CE^{\tau_{[j]},\mathrm{fin}}_{\on{exp}(2\pi i \chi_j)}(\la_1) \times \ldots \times \CE^{\tau_{[j]},\mathrm{fin}}_{\on{exp}(2\pi i \chi_j)}(\la_k)$ is an isomorphism of crystals (with the tensor product crystal structure on the right-hand side). 
 By Proposition~\ref{iso_fin_cr_res} we have the isomorphism of crystals $\CE^{\tau_{[j]},\mathrm{fin}}_{\on{exp}(2\pi i \chi_j)}(\la_l) \iso 
 \CE_{-\chi_j}(\la_l)$ given by the monodromy along the path connecting $\epsilon=0$ and $\epsilon=1/2N$ in the family 
 \begin{equation*}
 \epsilon \mapsto \on{Im}({\bf{ev}}_{\ul{z}/\epsilon+\ul{d}}(B(\on{exp}(2\pi i \epsilon \chi_j))) \ra \on{End}(V_1 \otimes \ldots \otimes V_k)),\, \epsilon \in (0,1/N), 
 \end{equation*}
 \begin{equation*}
 0 \mapsto \on{Im}(\CA_{-\chi_j}(2\pi i \cdot \ul{z}) \ra \on{End}(V_1 \otimes \ldots \otimes V_k))
 \end{equation*}
  so we just need to check that the bijection 
 \begin{equation}\label{iso_deal_with}
 \CE^{\tau_{[j]},\mathrm{fin}}_{\on{exp}(2\pi i \chi_j)}(\ul{\la}) \iso \CE_{-\chi_j}(\la_1) \times \ldots \times \CE_{-\chi_j}(\la_k)
 \end{equation} 
 given by the composition of the product of monodromies above and the monodromy along the path connecting $\varepsilon=0$ and $\varepsilon=1$ of the family 
 \begin{equation*}
 \varepsilon \mapsto \on{Im}({\bf{ev}}_{\ul{z}/\varepsilon+\ul{d}}B(\on{exp}(2\pi i\chi_j)) \ra \on{End}(V_{1} \otimes \ldots \otimes V_{k})),\, \varepsilon \in \BR^\times,     
 \end{equation*}
 \begin{equation*}
 0 \mapsto \on{Im}(\CA_{\on{exp}(-2\pi i\chi_j)}^{\otimes k} \ra  \on{End}(V_{1} \otimes \ldots \otimes V_{k}))  
 \end{equation*} is an isomorphism of crystals. 

Consider the family $\on{Im}({\bf{ev}}_{\ul{z}/\varepsilon+\ul{d}}(B(\on{exp}(2\pi i \epsilon \chi_j))) \ra \on{End}(V_1 \otimes \ldots \otimes V_k))$, $\epsilon, \varepsilon \in \mathbb{R}$. We can assume that $\chi_j$ is Weil generic so it follows from Proposition~\ref{two_stupid_limits} that this family extends to the family parametrized by $\on{Bl}_{(0,0)}K$.  
Our goal is to describe the path $p \colon [0,1] \ra \on{Bl}_{(0,0)}K$ that induces the isomorphism~(\ref{iso_deal_with}) and then  replace it by a homotopy equivalent path inside $\on{Bl}_{(0,0)}K$.

We start from describing the path $p \colon [0,1] \ra \on{Bl}_{(0,0)}K$. Note that $p(0)=((1,1/2N),2N:1])$, $p(1)=((0,0),[1:0])$.
Recall that the isomorphism~(\ref{iso_deal_with}) is the composition of two isomorphisms: one is the isomorphism 
\begin{equation*}
\CE^{\tau_{[j]},\mathrm{fin}}_{\on{exp}(2\pi i \chi_j)}(\ul{\la}) \iso \CE^{\tau_{[j]},\mathrm{fin}}_{\on{exp}(2\pi i \chi_j)}(\la_1) \times \ldots \times \CE^{\tau_{[j]},\mathrm{fin}}_{\on{exp}(2\pi i \chi_j)}(\la_k)
\end{equation*}
and the second one is the isomorphism 
\begin{equation*}
\CE^{\tau_{[j]},\mathrm{fin}}_{\on{exp}(2\pi i \chi_j)}(\la_1) \times \ldots \times \CE^{\tau_{[j]},\mathrm{fin}}_{\on{exp}(2\pi i \chi_j)}(\la_k) \iso \CE_{-\chi_j}(\la_1) \times \ldots \times \CE_{-\chi_j}(\la_k).
\end{equation*}
The first isomorphism is induced by the path $p_1\colon [0,1] \ra \on{Bl}_{(0,0)}K$ given by $c \mapsto ((\frac{1}{2N},\frac{1-c}{2N}),[1:1-c])$ and the second isomorphism is induced by the path $p_2\colon [0,1] \ra \on{Bl}_{(0,0)}K$ given by $c \mapsto ((\frac{1-c}{2N},0),[1:0])$.  So we get $p=p_1 * p_2$ i.e. $p$ is obtained by gluing $p_1,p_2$.


Consider now the following path $q \colon [0,1] \ra \on{Bl}_{(0,0)}K$ from $((0,0),[1:0])$ to $((1,1),[1:1])$. Path $q$ will be the gluing of two different paths $q_1,q_2$, $q=q_1*q_2$. Path $q_1\colon [0,1] \ra \on{Bl}_{(0,0)}K$ is given by  $c \mapsto ((0,0),[1:c])$, path $q_2$ is given by $c \mapsto ((\frac{c}{2N},\frac{c}{2N}),[1:1])$.

Consider the composition $p * q$. Note that $p * q$ is a cycle and it is easy to see that this cycle is homotopic to zero: indeed, note that we have a family of continuous maps $\gamma_t \colon \on{Bl}_{(0,0)}K \ra \on{Bl}_{(0,0)}K,\, ((a,b),[x:y]) \mapsto ((ta,tb),[x:y])$ that retracts $\on{Bl}_{(0,0)}K$ on $\mathbb{R}\BP^1$. The cycles $\gamma_t(p * q)$ are homotopic and $\gamma_1(p*q)=p*q$. It remains to note that $\gamma_0(p * q)$ is homotopic to zero. 

We conclude that the isomorphism induced by $p$ is equal to the isomorphism $\CE^{\tau_{[j]},\mathrm{fin}}_{\on{exp}(2\pi i \chi_j)}(\ul{\la}) \iso \CE_{-\chi_j}(\la_1) \times \ldots \times \CE_{-\chi_j}(\la_k)$ induced by $q^{-1}$. 
Note now that the isomorphism induced by $q^{-1}=q_2^{-1}* q_1^{-1}$ is the composition of the isomorphism $\CE^{\tau_{[j]},\mathrm{fin}}_{\on{exp}(2\pi i \chi_j)}(\ul{\la}) \iso \CE_{-\chi_j}(\ul{\la})$ induced by $q_1^{-1}$ and the isomorphism 
$\CE_{-\chi_j}(\ul{\la}) \iso \CE_{-\chi_j}(\la_1) \times \ldots \times \CE_{-\chi_j}(\la_k)$ 
induced by $q_2^{-1}$. It follows from Proposition~\ref{iso_fin_cr_res} that the first isomorphism is the isomorphism of crystals, it follows from~\cite{hkrw} that the second isomorphism is the isomorphism of crystals, where the crystal structure on $\CE_{-\chi_j}(\la_1) \times \ldots \times \CE_{-\chi_j}(\la_k)$ is given by $\CE_{-\chi_j}(\la_1) \otimes \ldots \otimes  \CE_{-\chi_j}(\la_k)$.

\epr

We are now ready to prove the main theorem of this section.
\th{}\label{KR_tensor_via_geometry}
For $C \in S^{\mathrm{reg}}$ and $\on{Im}z_1 \gg \on{Im}z_2 \gg \ldots \gg \on{Im}z_k$, the crystal $\mathcal{E}_C(\ul{\la})$ coincides with the Kirillov--Reshetikhin crystal ${\bf{B}}_{\la_1} \otimes {\bf{B}}_{\la_2} \otimes \ldots \otimes {\bf{B}}_{\la_k}$, corresponding to the tensor product of irreducible representations $V_{\la_1} \otimes V_{\la_2} \otimes \ldots \otimes V_{\la_k}$.
\eth
\prf
Follows from Propositions~\ref{r_1_crystal},~\ref{iso_with_tensor}.

\epr

\section{Monodromy}\label{monodromy}

\subsection{Spaces $\ol{M_{0,n+2}}$, $\ol{S^{\mathrm{reg}}}$}
\subsubsection{Space $\ol{M_{0,n+2}}$} Recall that $\ol{M_{0,n+2}}$ is the Deligne-Mumford space of stable rational curves with $n+2$ marked points. The points of $\ol{M_{0,n+2}}$ are isomorphism classes of curves of genus $0$, with $n+2$ ordered marked points and possibly with nodes, such that each component has at least $3$ distinguished points (either marked points or nodes). 

The space $\ol{M_{0,n+2}}$ is a compactification of $M_{0,n+2}$ (configuration space of $n+2$ ordered pairwise distinct points $(z_{1},\ldots,z_{n+2})$ of $\BP^1$ modulo $\on{PGL}_2$). We have 
\begin{equation}\label{descr_M_0}
M_{0,n+2} \simeq \{(z_1,\ldots,z_n) \in (\BC^\times)^n\,|\, z_i \neq z_j\}/\BC^\times = T^{\mathrm{reg}},
\end{equation}
 where $T^{\mathrm{reg}} \subset \on{PGL}_n$ is the subalgebra of diagonal matrices.

Space $\ol{M_{0,n+2}}$ has a natural stratification indexed by the combinatorial type of the curve with marked points. Namely, the strata are indexed by trees whose leaves correspond to the marked points $0,z_1,\ldots,z_n,\infty$ where the inner vertices represent the connected components and the edges correspond to the nodal points.


\subsubsection{Space $\ol{S^{\mathrm{reg}}}$} We will be interested in a certain real form of  $\ol{M_{0,n+2}}$ that we will denote by $\ol{S^{\mathrm{reg}}}$ since it is a compactification of $S^{\mathrm{reg}}$. Namely, consider the (anti-holomorphic) involution on $\ol{M_{0,n+2}}$ sending $(C,z_0,z_1,\ldots,z_n,z_{n+1})$ to  $(C,\ol{z}_{n+1},\ol{z}_1,\ldots,\ol{z}_n,\ol{z}_{0})$.
The space $\ol{S^{\mathrm{reg}}}$ is the fixed-point locus of this involution. Note that after the identification (\ref{descr_M_0}) the involution on $M_{0,n+2}$ is given by $(z_1,\ldots,z_n) \mapsto (\ol{z}_1^{-1},\ldots,\ol{z}_n^{-1})$ so its fixed points are precisely $S^{\mathrm{reg}}$. The whole space $\overline{S^{\mathrm{reg}}}$ is the moduli space of curves where all marked points are on the unit circle, and the points $0$, $\infty$ belong to the same component. 
We refer  the reader to \cite[Section 7.2.4]{ikpr} and \cite{ce} for details. 

Stratification of $\ol{M_{0,n+2}}$ induces the stratification of $\ol{S^{\mathrm{reg}}}$.

\subsubsection{Extended affine cactus group}
See \cite[Section 10.3]{ikpr} for the details. 

We recall that the affine cactus group $AC_n$ is the group with generators $s_{ij}$ for $1 \leqslant i \neq j \leqslant n$ and the following relations. We consider the natural cyclic order on the unit circle $S^1$ and use the embedding  $\{1,\ldots,n\} \hookrightarrow S^1,\, i \mapsto e^{\frac{2\pi i}{n}}$ to obtain the cyclic order on $\{1,2,\ldots,n\}$. Then, for every pair $i \neq  j$ of elements of $\{1,\ldots,n\}$ it makes sense to talk about the segment $[i,j] \subset S^1$. The defining relations of $AC_n$ are: 

$(1)$ $s_{ij}^2=1$,

$(2)$ $s_{ij}s_{kl}=s_{kl}s_{ij}$ if $[i,j] \cap [k,l]=\varnothing$,

$(3)$ $s_{ij}s_{kl}=s_{w_{ij}(l)w_{ij}(k)}s_{kl}$ if $[k,l] \subset [i,j]$,
where $w_{ij}$ is the element which reverses $[i,j]$ and leaves invariant the elements outside the interval.


We define an action of $Z/n\BZ$ on $AC_n$ by $r \cdot s_{ij}=s_{i+1\,j+1}$. The extended affine cactus group $\widetilde{AC}_n$ is the semi-direct product $(\BZ/n\BZ) \ltimes AC_n$.

\subsubsection{Equivariant fundamental group of $\ol{S^{\mathrm{reg}}}$} 
Group $S_n$ acts naturally on the space $\overline{S^{\mathrm{reg}}}$.
We will be interested in the $S_n$-equivariant fundamental group of $\overline{S^{\mathrm{reg}}}$ which is defined as follows. Pick a basepoint $x \in \overline{S^{\mathrm{reg}}}$ which corresponds to a configuration of evenly spaced points of $U(1)$, namely $z_k=e^{\frac{2\pi i k}{n}}$.
\begin{equation*}
\pi_1^{S_n}(\overline{S^{\mathrm{reg}}}) := \{(g,p)\,|\, g \in S_n,~\text{$p$ is a homotopy class of paths from $x$ to $gx$}\}.
\end{equation*}
The multiplication in $\pi_1^{S_n}(\overline{S^{\mathrm{reg}}})$ is defined as follows:
\begin{equation*}
(g_1,p_1) \cdot (g_2,p_2) = (g_1g_2,p_1 * g_1(p_2)).
\end{equation*}



Abusing notations for $1 \leqslant i \neq j \leqslant n$ let $(w_{ij},s_{ij}) \in \pi^{S_n}_1(\overline{S^{\mathrm{reg}}})$ be the element corresponding to the shortest path $\gamma_{ij}$ connecting the base point $x$ with $w_{ij}x$.

\begin{propos}{}{\em{(\cite[Theorem 8.3]{ce} and \cite[Theorem 11.12]{ikpr})}} There is an isomorphism $\widetilde{AC}_n \iso \pi_1^{S_n}(\overline{S^{\mathrm{reg}}})$ that sends 
$
s_{ij} \in \widetilde{AC}_n$  to $(w_{ij},s_{ij}) \in \pi^{S_n}_1(\overline{S^{\mathrm{reg}}})$
and $[1] \in \BZ/n\BZ$  to the pair consisting of the cyclic permutation $(12\ldots n)$ and the shortest path connecting $x$ with $(12\ldots n)x$.
\end{propos}

\subsection{Description of the monodromy}
The Deligne-Mumford space $\overline{M_{0,n+2}}$, being a compactification of $T^{\mathrm{reg}}$, parametrizes all possible limits of the Bethe subalgebras $B(C)$ in the Yangian (see \cite{ir0}). According to \cite{imr}, any subalgebra $B(X)$ corresponding to any point $X\in \overline{M_{0,n+2}}$ acts on $\bigotimes V_{\lambda_i}(u_i)$ with a cyclic vector. Restricting to $\overline{S^{\mathrm{reg}}}\subset\overline{M_{0,n+2}}$, and requiring the evaluation parameters $u_i$ to be as above, we guarantee that the action of the corresponding subalgebra on the corresponding module is Hermitian hence semisimple. So we can regard $\CE_{C}(\ul{\la})$ as a $S_n$-equivariant covering of the space $\overline{S^{\mathrm{reg}}}\subset\overline{M_{0,n+2}}$. Our goal in this section is to confirm \cite[Conjecture 7.1]{imr} (see also \cite[Conjecture 1.6]{ikr}) describing the monodromy of this covering in terms of the above affine crystal structure (compare with \cite{hkrw}).

We claim that the monodromy of this covering can be expressed in terms of the partial Sch\"utzenberger involutions for the  Kirillov--Reshetikhin crystals. Namely, for any proper subdiagram in the affine diagram $\widetilde{A_{n-1}}$ (i.e. for any proper subset in the set of affine simple roots), we can decompose our KR crystal into connected components with respect to the corresponding (finite-type) Levi subalgebra and apply the Sch\"utzenberger involution assigned to this Levi to each of the components.


Recall that if $B$ is a normal $\mathfrak{g}$-crystal then the Sch\"utzenberger involution is the {\em{unique}} isomorphism $\xi_B\colon B \iso B$ that preserves connected components of $B$ and such that:
\begin{equation*}
e_i(\xi_B(b))=\xi_B(f_{\theta(i)}(b)),
\end{equation*}
\begin{equation*}
f_i(\xi_B(b))=\xi_B(e_{\theta(i)}(b)),
\end{equation*}
\begin{equation*}
\on{wt}(\xi_B(b))=w_0\on{wt}(b),
\end{equation*}
where $w_0 \in W$ is the longest element and $\theta$ is the involution of the set of simple roots induced by $-w_0$.

Recall that to $C \in S^{\mathrm{reg}}$ and $w \in \widehat{W}$ we can associate the Kirillov--Reshetikhin crystal $\CE_{C}^{w}$. Abusing the notations, we denote by $\CE$ the KR crystal corresponding to $C \in Q^{\mathrm{reg}}$, we can identify it with crystals corresponding to other alcoves using the action of $\widehat{W}$.

\begin{thm}{}\label{th:monodromy}
The generators $s_{ij}$ of the monodromy group of $\overline{S^{\mathrm{reg}}}$ act on $\mathcal{E}$ as a partial Sch\"utzenberger involution with respect to the corresponding subdiagram in the  Dynkin diagram of $\widetilde{A_{n-1}}$, see \cite[Conjecture~7.1]{imr} and \cite[Conjecture 1.6]{ikr}. 
\end{thm}
\begin{proof}
Let us deal with the case of $s_{1n}$, the other cases are similar. We only need to check that the bijection induced by the monodromy along the path $\gamma_{1n}$ is compatible with crystal structures and preserves connected components. The path $\gamma_{1n}$ can be replaced by the following path $\gamma$. Let $\chi=\on{diag}(u_1,\ldots,u_n)$ be an element of $\mathbb{R}^n$ such that $w_0(\chi)=-\chi$ (i.e., $u_j=-u_{n+1-j}$) and all $u_i$ are distinct. Consider the path $\gamma^{+}(t)=\on{exp}(2\pi i t\chi) \in S^{\mathrm{reg}}$, $t \in (0,1]$. We can extend it to $0$ by setting $\gamma^{+}(0):=\underset{t \rightarrow 0}{\on{lim}}\gamma^+(t) \in \overline{S^{\mathrm{reg}}}$. Set $\gamma^{-}\colon [0,1] \rightarrow \overline{S^{\mathrm{reg}}}$, $\gamma^{-}(t)=w_0(\gamma^{+}(1-t))$ (note that for $t \neq 1$, we have $\gamma^-(t)=\operatorname{exp}(2\pi i (t-1)\chi)$). The path $\gamma\colon [0,1] \rightarrow \overline{S^{\mathrm{reg}}}$ is the concatenation of the paths $\gamma^+$ and $\gamma^-$.
Note now that $B(\gamma(\frac{1}{2}))$ is generated by $B(1)$ together with the image of $\mathcal{A}_\chi$.
Let us now decompose our KR-module as the direct sum of weight spaces of $B(1)$. Each weight space will be a $\mathfrak{g}$-module without multiplicities since otherwise, the action of $\mathcal{A}_\chi$ on this component will not be cyclic (and so the action of $B(1) \cdot \mathcal{A}_\chi$ cannot have simple spectrum on the whole module). It then follows that our bijection preserves components. 
Using Corollary \ref{restr_ident} we see that the crystal structure on the spectrum of the algebra $B(1) \cdot \mathcal{A}_{\chi}$ after restricting to the appropriate Levi (determined by the alcove of $\chi$) coincides with the one defined in \cite[Section 12.1]{hkrw}. Now, the compatibility with the crystal structure follows from  \cite[Proposition 13.5]{hkrw}. 
\end{proof}

\begin{remark}{}
Theorem \ref{th:monodromy} confirms \cite[Conjecture 7.1]{ikr} (see also \cite[Conjecture 1.6]{ikr}). We expect that Theorem \ref{th:monodromy} also follows from
\cite[Theorem 1.3(3)]{kar} on the monodromy of spectra of trigonometric Gaudin model in terms of crystal commutors, 
combined with \cite[Theorem~5.1]{fru} on Howe duality that switches trigonometric Gaudin model with $XXX$ Heisenberg model in type $A$ (and goes back to \cite{mtv}). 
\end{remark}

\subsection{Application: monodromy of the spectrum of $QH^*_{\mathsf{A} \times \BC^\times}(\mathfrak{M})$ and cross-walling bijections}
Let $\mathfrak{M}=\mathfrak{M}(\mathsf{w})$ be a disjoint union of type $A$ quiver varieties $\mathfrak{M}({\mathsf{v}},{\mathsf{w}})$ with fixed framing vector ${\mathsf{w}}=(w_i)$ and arbitrary dimension vector ${\mathsf{v}}=(v_i)$ ($w_i, v_i$ are collections of non-negative numbers, labeled by the vertices of the quiver). 
Let $QH^*_{\mathsf{A} \times \BC^\times}(\mathfrak{M})$ be the algebra of {\emph{quantum cohomology}} of $\mathfrak{M}$, one can think about this algebra as about the {\em{family}} of subalgebras in $\on{End}(H^*_{\mathsf{A} \times \BC^\times}(\mathfrak{M}))$ depending (in logarithmic coordinates) on an element of $H^2(\mathfrak{M},\BC)$. Recall now that by the results  of Varagnolo \cite{varagnolo_yangians_on_quivers} (following \cite{nakajima_reps_of_quantum}), and Maulik--Okounkov (\cite{mo}) (see \cite{mcbreen_thesis} for the comparison of the two approaches) 
\begin{equation*}
H^*_{\mathsf{A} \times \BC^\times}(\mathfrak{M}) \simeq V_{\la_1}(z_1) \otimes \ldots \otimes V_{\la_k}(z_k)
\end{equation*}
and the action of the Yangian on it has a  {\em{geometric realization}} (here $\la_j$ are fundamental weights that can be constructed from the framing $(w_i)$). After these identifications, the action of  $QH^*_{\mathsf{A} \times \BC^\times}(\mathfrak{M})$ on $V_{\la_1}(z_1) \otimes \ldots \otimes V_{\la_k}(z_k)$ identifies with the action of the family of Bethe subalgebras of the Yangian (our parameter $C$ above naturally corresponds to an element of $H^2(\mathfrak{M},\BC)$), see \cite{mo}. 
So, the following theorem is a corollary of our results.  
\begin{thm}{}\label{thm_monodromy_qh}
The monodromy operators of eigenvalues of multiplication by elements of the quantum cohomology ring $QH^*_{\mathsf{A} \times \BC^\times}(\mathfrak{M})$ are generated by partial Sch\"utzenberger involutions. 
\end{thm}

\begin{remark}{}
The result of Theorem \ref{thm_monodromy_qh} can be also proven for an arbitrary ADE quiver variety such that $w_i \neq 0$ only for {\emph{minuscule}} fundamental weights $\omega_i$. 
\end{remark}

In the joint work in progress \cite{rep_th_gies_II} with Ivan Losev, the first author studies the cross-walling bijections for ADHM spaces (also known as Gieseker varieties) and for type $A$ quiver varieties.  

The space of cocharacters $\nu$ naturally splits into {\em{chambers}} separated by {\em{walls}}. It turns out that when $\nu, \nu'$ lie in adjacent chambers that share a codimension
one face, then the corresponding categories $\CO_\nu(\CA_{c})$, $\CO_{\nu'}(\CA_c)$ are {\em{derived equivalent}}. The derived equivalence is called the {\em{cross-walling functor}} (it is expected to be {\em{symplectically dual}} to the well-known {\em{wall-crossing}} functor). Losev proved in  \cite{losev} that the associated graded to cross-walling functor (w.r.t. a certain filtration) is ``perverse'' in the sense of Rouquier so it induces a bijection $\on{Irr}\CO_\nu(\CA_c) \iso \on{Irr}\CO_{\nu'}(\CA_c)$ between irreducible objects of our categories $\CO$. This bijection is called a {\em{cross-walling bijection}}. We have the natural identifications of both $\on{Irr}\CO_\nu(\CA_c)$ and $\on{Irr}\CO_{\nu'}(\CA_c)$ with $X^T$ so the cross-walling bijection is an automorphism of the set $X^T$.
In \cite{rep_th_gies_II}, we prove the following theorem. 

\begin{thm}{}\label{thm_cw_bij_ADHM}
Let $\mathfrak{M}({\mathsf{v}},{\mathsf{w}})$ be a type $A$ quiver variety with $\mu:=\sum_{i}w_i\omega_i-\sum_iv_i\alpha_i$ being dominant. Then, after appropriate identifications, cross-walling bijections for $\mathfrak{M}({\mathsf{v}},{\mathsf{w}})$ can be described explicitly: they are given by partial Sch\"utzenberger involutions. 
\end{thm}

Combining Theorems \ref{thm_cw_bij_ADHM}, \ref{thm_monodromy_qh} we see that the monodromy for $QH^*_{\mathsf{A} \times \BC^\times}(\mathfrak{M})$ coincides (after appropriate identifications) with the cross-walling bijection for the symplectically dual variety to $\mathfrak{M}$ (which in the case when $\mu$ is dominant is again a type $A$ quiver variety). Cross-walling bijections should be interchanged with wall-crossing bijections via the symplectic duality (cf. \cite[Section 10.1]{blpw}).
So, our computation supports the general conjecture formulated in the last paragraph of \cite[Section 1.13]{hkrw} being applied to the type $A$ Nakajima quiver varieties.

\appendix\label{app}

\section{Notation}

Here is a list of the notation used in the paper.

\ssec{}{Lie groups, Lie algebras, roots, weights, representations}
\begin{itemize}
    \item $\mathfrak{g}$, simple finite-dimensional Lie algebra
    \item $G$, adjoint Lie group with Lie algebra 
    $\mathfrak{g}$
    \item $\mathfrak{h}$, Cartan subalgebra of $\mathfrak{g}$
    \item $T$, Cartan subgroup of $G$
    \item $\varpi_i$, fundamental weight of $\mathfrak{g}$
    \item $\Delta_+$, positive roots of $\mathfrak{g}$
    \item $\{\al_i\}_{i=1,\ldots,\on{rk}\mathfrak{g}}$, the set of simple roots of $\mathfrak{g}$
\item $P^\vee$, weight lattice of $\mathfrak{sl}_n$
\item $\La$, coweight lattice of $\mathfrak{gl}_n$ or $\mathfrak{pgl}_n$ 
\item $\La_r$, the coroot lattice of $\mathfrak{gl}_n$ or $\mathfrak{pgl}_n$ (depending on the context)
\item $\mathfrak{h}_0$, subalgebra of traceless diagonal matrices
\item $\mathfrak{h}^{\mathrm{reg}}_0$, traceless diagonal matrices with distinct eigenvalues 
\item $T$, diagonal matrices in $\on{GL}_n$ or $\on{PGL}_n$ (depending on the context)
\item $T^{\mathrm{reg}}$, diagonal matrices in $\operatorname{GL}_n$ or  $\on{PGL}_n$ with distinct eigenvalues
\item $S \subset U(n)$, compact subtorus of diagonal matrices
\item $S^{\mathrm{reg}} \subset S$, regular elements of $S$
\item $V_\la$, irreducible finite-dimensional representation of $\mathfrak{g}$, corresponding to $\la$
\item $\ul{V}(s)$, the representation $V_{1}(sz_1+d_1) \otimes \ldots \otimes V_{k}(sz_k+d_k)$ of $Y(\mathfrak{g})$
\item $\mathfrak{z}_{\mathfrak{g}}(\chi)$, centralizer of $\chi \in \mathfrak{g}$
\item $\hat{\mathfrak{g}}$, affine Lie algebra, corresponding to $\mathfrak{g}$
\end{itemize}

\ssec{}{General algebra and geometry}
\begin{itemize}
    \item $\on{dim}_{Q} D$ for a filtration $Q^{\bullet}$ on a vector space $U$ and $D \subset U$ we have $\on{dim}_Q D:=\sum_{i \geqslant 0} \on{dim}\left(\frac{D \cap Q^iU}{D \cap Q^{i-1}U}\right) q^i \in \BZ[q].    
$
    \item $Rees(A)$, Rees algebra of a filtered algebra $A$
    \item If $a \in R$ is an element of some algebra $R$ then $a^{(i)}=1 \otimes \ldots \otimes 1 \otimes \underset{i}{a} \otimes 1 \otimes \ldots \otimes 1$
    \item $\ol{M_{0,n+2}}$, the Deligne-Mumford space of stable genus $0$ curves with $n+2$ marked points 
    \item $\mathfrak{M}=\mathfrak{M}(\mathsf{w})=\bigsqcup_{\mathsf{v}}\mathfrak{M}(\mathsf{v},\mathsf{w})$, disjoint union of type $A$ Nakajima quiver varieties with fixed framing $\mathsf{w}$
    \item ${\mathsf{A}} \times \mathbb{C}^{\times}$, torus naturally acting on $\mathfrak{M}$
\end{itemize}

\ssec{}{Specific elements and operators}
\begin{itemize}
    \item $\omega'=\frac{1}{2}\sum_{\al \in \Delta^+}x_{\al}^{\pm} \otimes x_{\al}^{\mp}+\sum_{j=1}^{\on{rk}\mathfrak{g}} \frac{1}{(\al_j,\al_j)}h_j \otimes t_{\varpi_j} \in Y(\mathfrak{g})$
    \item $\omega_\chi = \frac{1}{2}\sum_{a}x_a[0]^2+ \chi[1] \in S^\bullet(\mathfrak{g}[t])$
    \item $\tilde{\omega}_\chi = \frac{1}{2}\sum_{a}x_a[0]^2+ \chi[1] \in U(\mathfrak{g}[t])$
    \item $\Omega_\chi= \sum_{a}x_a[0]x_a[1]+\chi[2]  \in S^\bullet(\mathfrak{g}[t])$
    \item $\tilde{\Omega}_\chi = \sum_{a}x_a[0]x_a[1]+\chi[2]\in U(\mathfrak{g}[t])$
    \item $\tilde{\Omega}_{\chi,\epsilon}=\sum_a \hbar^2 \epsilon^{-3} x_ax_{a}[1]+\epsilon^{-3}\hbar \chi[2] \in U_\epsilon(\mathfrak{g}[t])$
    \item $\tau_{[j]}\colon \BZ/n\BZ \ra \BZ/n\BZ$, the permutation that sends $[i]$ to $[i+j]$
    \item $\xi_B$, Sch\"utzenberger involution of a normal $\mathfrak{sl}_n$-crystal $B$
    \item  
${\bf{pr}}$, promotion operator
    \item 
$L_{\ul{z}}=\sum_{i=1}^k \frac{E^{(i)}}{u-z_i}$
\item
$E \in \on{End}(\BC^n) \otimes U(\mathfrak{gl}_n)$ matrix $(E_{ij})_{i,j=1,\ldots,n}$
\item $A_m=\frac{1}{m!}\sum_{\sigma \in S_m}\sigma \subset \on{End}(\BC^n)^{\otimes n}$
\item $c_{\mathfrak{g}}$, the scalar by which $\sum_a x_a^2$ acts on $\mathfrak{g}$
\item $(C_0,\chi)$, the pair of $C_0=\on{exp}(2\pi i \chi_0)$, $\chi \in Q_{\hat{w}}^{\mathrm{reg}}$, here $\chi_0 \in Q_{\hat{w}}$ is a subregular element
\end{itemize}

\ssec{}{Conformal blocks}
\begin{itemize}
\item $\tilde{\mathfrak{g}}(z)=\mathfrak{g}((t-z))$
\item $\tilde{\mathfrak{g}}(\ul{z})=\bigoplus_{i=1}^k \tilde{\mathfrak{g}}(z_i)$
\item $\hat{\mathfrak{g}}(z)$, one-dimensional central extension of $\mathfrak{g}(z)$
\item $\hat{\mathfrak{g}}(\ul{z})$, one-dimensional central extension of $\tilde{\mathfrak{g}}(\ul{z})$
\item $\mathfrak{g}_{\ul{z}}$, rational $\mathfrak{g}$-valued functions on $\BP^1 \setminus \{z_1,\ldots,z_k\}$
\item $\mathbb{V}_{0,z}$, the induced module $\on{Ind}^{\hat{\mathfrak{g}}}_{\mathfrak{g}(z) \oplus \BC \mathsf{K}}\BC$
\item $\mathbb{M}_z$, the induced module $\on{Ind}^{\hat{\mathfrak{g}}(z)}_{\mathfrak{g}[[t]] \oplus \BC \mathsf{K}} M$
 \item $H(\mathbb{M}_{z_1},\ldots,\mathbb{M}_{z_k})$, space of conformal blocks
 \item $\mathring{D}=\on{Spec}\BC((t))$, $\mathring{D}_z=\on{Spec}\BC((t-z))$
 \item $\on{Op}(X)$, moduli space of $G^\vee$-opers on $X$
\end{itemize}

\ssec{}{Commutative subalgebras}
\begin{itemize}
\item $\CA^{\mathrm{u}}_\chi$, inhomogeneous universal  Gaudin subalgebra of $U(\mathfrak{g}[t])$
\item $\tilde{\ol{\CA}}{}^{\mathrm{u}}_\chi$, inhomogeneous universal  Gaudin subalgebra of $U(\mathfrak{gl}_n[t])$
\item $\ol{\CA}{}^{\mathrm{u}}_\chi$, classical inhomogeneous universal  Gaudin subalgebra of $S^\bullet(\mathfrak{g}[t])$
\item $\tilde{\ol{\CA}}{}^{\mathrm{u}}_\chi$, inhomogeneous universal  Gaudin subalgebra of $U(\mathfrak{gl}_n[t])$
\item $\CA_{(\chi_0,\chi)}^{\mathrm{u}}$, the limit $\underset{\epsilon \ra 0}{\on{lim}}\,\CA^{\mathrm{u}}_{\chi_0+\epsilon\chi}$
\item $\ol{\CA}{}^{\mathrm{u}}_{(\chi_0,\chi)}$, the limit $\underset{\epsilon \ra 0}{\on{lim}}\,\ol{\CA}{}^{\mathrm{u}}_{\chi_0+\epsilon\chi}$
\item $B(C)$, Bethe subalgebra of $Y(\mathfrak{g})$, corresponding to $C \in G$
\item $\tilde{B}(C)$, Bethe subalgebra of $Y(\mathfrak{gl}_n)$, corresponding to $C \in \on{GL}_n$
\item $\CA_\chi(u_1,\ldots,u_k)$, inhomogeneous Gaudin subalgebra of $U(\mathfrak{g})^{\otimes k}$
\item $\CA_{(\chi_0,\chi)}(u_1,\ldots,u_k)$, the limit $\underset{\epsilon \ra 0}{\on{lim}}\,\CA_{\chi_0+\epsilon \chi}(u_1,\ldots,u_k)$
\item $B(C_0,\chi)$, the limit $\underset{\epsilon \ra 0}{\on{lim}}\,B(C_0\on{exp}(2\pi \epsilon \chi))$ 
\item $\CA(M_1,\ldots,M_k)$, commutative subalgebra of $\on{End}(M_1 \otimes \ldots \otimes M_k)$ defined via conformal blocks
\item $Z$, Poisson center of the completion $\widetilde{S}^\bullet(\mathfrak{g}((t^{-1})))$
\item $\CZ$, center of the completion $\widetilde{U}(\hat{\mathfrak{g}})_{-1/2}$
\end{itemize}

\ssec{}{Chambers and alcoves}
\begin{itemize}
\item $W$, Weyl group of $\mathfrak{g}$
\item $O_w$, closed chamber (w.r.t. $W$), corresponding to $w \in W$
\item $O^{\mathrm{reg}}_w$, open chamber (w.r.t. $W$), corresponding to $w \in W$
\item $H_\al$, wall  for the action $W \curvearrowright \mathfrak{h}$
    \item $\widehat{W}$, semidirect product $S_n \ltimes \ol{\La}_r$
    \item $\widehat{W}^{\mathrm{ext}}$, semidirect product $S_n \ltimes \ol{\La}$
    \item $Q_{\hat{w}}$, closed alcove (w.r.t. $\widehat{W}$), corresponding to $\hat{w} \in \widehat{W}$
       \item $Q_{\hat{w}}^{\mathrm{reg}}$, open alcove (w.r.t. $\widehat{W}$), corresponding to $\hat{w} \in \widehat{W}$\item $H_{[i],[j]}^k$, affine wall for the action $\widehat{W} \curvearrowright \ol{\mathfrak{h}}$
\end{itemize}

\ssec{}{Yangians}
\begin{itemize}
\item $Y(\mathfrak{g})$, Yangian of $\mathfrak{g}$ in Drinfeld's realization
\item $Y_V(\mathfrak{g})$, Yangian of $\mathfrak{g}$ in $RTT=TTR$ realization 
\item $\CR$, universal $R$-matrix for $Y(\mathfrak{g})$ 
\item ${\bf{ev}}_{z_1,\ldots,z_k}\colon Y(\mathfrak{g}) \ra U(\mathfrak{g})^{\otimes k}$, $\on{ev}_{u_1,\ldots,u_k}\colon U(\mathfrak{g}[t]) \ra U(\mathfrak{g})^{\otimes k}$, $\on{ev}'_{u_1,\ldots,u_k}\colon U(t^{-1}\mathfrak{g}[t^{-1}]) \ra U(\mathfrak{g})^{\otimes k}$, $\on{ev}'_{\infty}\colon U(t^{-1}\mathfrak{g}[t^{-1}]) \ra S^\bullet(\mathfrak{g})$, $\on{ev}'_{u_1,\ldots,u_k,\infty}\colon U(t^{-1}\mathfrak{g}[t^{-1}]) \ra U(\mathfrak{g})^{\otimes k} \otimes S^\bullet(\mathfrak{g})$, $\on{ev}'_{u_1,\ldots,u_k,\chi}\colon U(t^{-1}\mathfrak{g}[t^{-1}]) \ra U(\mathfrak{g})^{\otimes k}$, evaluation homomorphisms
\item  $\Delta^k\colon Y(\mathfrak{g}) \ra Y(\mathfrak{g})^{\otimes k}$, comultiplication homomorphism
\item $\eta\colon Y(\mathfrak{g}) \ra Y(\mathfrak{g})$, antipode
\item $V(\la,0)$, irreducible representation of $Y(\mathfrak{g})$, corresponding to dominant weight $\la$
\end{itemize}

\ssec{}{Crystals}
\begin{itemize}
\item For a $\mathfrak{sl}_n$-crystal $B$, the Kashiwara operators $\on{e}_i,\on{f}_i\colon B \ra B \cup \{0\}$, $i=1,\ldots,n-1$, $\on{wt}\colon B \ra P^\vee$
\item For a $\hat{\mathfrak{sl}}_n$-crystal ${\bf{B}}$, the Kashiwara operators $\on{e}_{[i]},\on{f}_{[i]}\colon B \ra B \cup \{0\}$, $[i] \in \BZ/n\BZ$, $\on{wt}\colon {\bf{B}} \ra P^\vee$.
 
    \item $B_\la$, the $\mathfrak{sl}_n$-crystal of the representation $V_\la$
    \item ${\bf{B}}_{l\varpi_r}$, the Kirillov--Reshetikhin crystal of $V_{l\varpi_r}$
    \item $\CE_C(\ul{\la})$, the set of eigenlines of $B(C) \curvearrowright V_{\la_1} \otimes \ldots \otimes V_{\la_k}$
    \item $\CE_{\chi}(\ul{\la})$, the set of eigenlines of $\CA_\chi(z_1,\ldots,z_k) \curvearrowright V_{\la_1} \otimes \ldots \otimes V_{\la_k}$
  \end{itemize}

\end{document}